\documentclass[a4paper,12pt]{article}
\usepackage{longtable}
\usepackage{setspace}
\usepackage{multirow}
\usepackage{epsfig}
\usepackage{amsmath}
\usepackage{amssymb}
\usepackage{tabularx}
\usepackage{enumerate}
\usepackage{graphicx}
\usepackage{subfigure}
\usepackage{cases}
\usepackage[compress]{cite}
\usepackage{lineno}
\usepackage{color}
 \usepackage{epstopdf}
\usepackage{geometry}
\newtheorem{thm}{Theorem}[section]
\newtheorem{lem}[thm]{Lemma}
\newtheorem{cor}[thm]{Corollary}

\newtheorem{rmk}[thm]{Remark}
\newtheorem{defi}{Definition}[section]
\newtheorem{pppp}{Proof}

\newcommand{\qed}{\hspace{1em}\mbox{\raisebox{0.65ex}{\fbox{}}}}

\numberwithin{equation}{section}

\newcommand{\be}{\begin{equation}}
\newcommand{\ee}{\end{equation}}
\newcommand\bes{\begin{eqnarray}} \newcommand\ees{\end{eqnarray}}
\newcommand{\bess}{\begin{eqnarray*}}
\newcommand{\eess}{\end{eqnarray*}}
\newcommand{\bpf}{{\bf Proof:\ \ }}
\newcommand{\epf}{\mbox{}\hfill $\Box$}

\begin{document}

\thispagestyle{empty}
\title{A West Nile virus nonlocal model with free boundaries and seasonal
succession\thanks{The work is partially supported by the NNSF of China (Grant No. 11771381), Graduate
Research and Innovation Projects of Jiangsu Province (XKYCX20-011).}}
\date{\empty}
\author{Liqiong Pu$^{a}$,\  Zhigui Lin$^a$ and Yuan Lou$^{b, c}$\\
{\small $^a$ School of Mathematical Science, Yangzhou University, Yangzhou 225002, China}\\
{\small $^b$ School of Mathematical Sciences, Shanghai Jiao Tong University, Shanghai 200240, China}\\
{\small $^c$ Department of Mathematics, Ohio State University, Columbus, OH 43210, USA}
}
\maketitle
\begin{quote}
\noindent
{\bf Abstract.} {
 \small
 The paper deals with a West Nile virus (WNv) model, where the nonlocal diffusion
  characterizes a long-range dispersal, the free boundary
  describes the spreading front, and seasonal succession accounts for the effect of the warm and cold seasons. The well-posedness of the model is
 established, and its long-term dynamical behaviours,
 which depend on the generalized eigenvalues of the corresponding linear operator,
 are investigated. For both spatially independent and nonlocal WNv models with seasonal successions, the generalized eigenvalues are determined.
We develop the indexes to the case with the free boundary and
apply these indexes to determine whether spreading or vanishing happens. Our new criteria extends
previous results for the case with the nonlocal diffusion and the case with the free boundary.
The generalized eigenvalues reveal that there exists positive correlation between the duration of the warm season and the risk of infection.
Moreover, the initial infection length, the initial infection scale and the spreading ability to the new area play an important role for the long time behavior of the solution.

}


\end{quote}

\section{Introduction}
\setcounter{equation}{0}

West Nile virus (WNv) is an emerging mosquito-borne virus that can cause a severe, life-threatening neurological disease in humans and horses, and it is widely distributed throughout the world with considerable impact on both public health and animal health \cite{BJL2013}.

For nearly two decades, many mathematical models for WNv have been proposed and studied.
However, most models are focused on the  non-spatial transmission dynamics \cite{BGD, WBL, WZ}.
In fact, what we should actually do is to consider the spatial spreading, which is an important factor to affect the persistence and eradication of WNv.
To utilize the cooperative characteristic of cross-infection dynamics and estimate the spatial spread rate of infection, Lewis et al. in \cite{LRD}  proposed the following simplified spatial-dependent WNv model
\be\left\{
\begin{array}{ll}
u_{1t}=d_1\Delta u_1 +\alpha_{b}\beta_{b}\frac{(N_b-u_1)}{N_b} u_2 -\gamma_{b}u_1, & t>0,\ x\in\Omega, \\[2mm]
u_{2t}=d_2\Delta u_2 +\alpha_{m}\beta_b \frac{(A_m-u_2)}{N_b} u_1 -d_{m} u_2,  & t>0,\ x\in\Omega,
\end{array} \right.
\label{a01} \ee
where $u_1(t,x)$ and $u_2(t,x)$ represent the population
densities of infected birds and mosquitoes at the location $x$ and time $t\geq 0$, respectively, and $0<u_1(0,x)\leq N_b, 0<u_2(0,x)\leq A_m$. The total population of birds $N_b$ and mosquitoes $A_m$ are 
assumed to be positive constants, $d_1$ and $d_2$ represent diffusion coefficients for birds and mosquitoes, respectively. $\alpha_{b}$ and $\alpha_{m}$ account for WNv transmission probability per bite to birds and mosquitoes. $\beta_{b}$ is the biting rate of mosquitoes on birds, $\gamma_b $ is the recovery rate of birds from WNv, and $d_m$ is adult mosquitos' death rate.

For convenience of mathematical analysis, we denote the parameters in \eqref{a01} by
$$a_1=\frac{\alpha_{b}\beta_{b}}{N_b},\, e_1=N_b, b_1=\gamma_b,\, a_2=\frac{\alpha_{m}\beta_b}{N_b},\, e_2=A_m,\, b_2=d_m.$$
Under the new notations, system \eqref{a01} is rewritten as
\be\left\{
\begin{array}{ll}
u_{1t}=d_1\Delta u_1 +a_1(e_1-u_1)u_2 -b_1 u_1, &t>0,\  x\in\Omega, \\[2mm]
u_{2t}=d_2\Delta u_2 +a_2(e_2-u_2) u_1 -b_2 u_2,  &t>0,\ x\in\Omega.
\end{array} \right.
\label{a02} \ee
For the ODE version of \eqref{a02}, the authors in \cite{LRD}  derived the basic reproduction number
\begin{equation}
R_0=\sqrt{\frac{a_1a_2e_1e_2}{b_1b_2}}
\label{a03}
\end{equation}
by the next generation matrix method \cite{VW} and illustrated that the virus vanishes for $R_0<1$, while for $R_0>1$, the disease-endemic equilibrium stabilizes.
They further considered the existence of traveling waves of \eqref{a02} and deduced that the spread rate, which is determined by a linearized system, is equivalent to the minimal wave speed of the non-linear one.

Although \eqref{a02} can be used to estimate the speed of disease transmission, it can not be used to understand the spreading front of the infected areas. Recently, Lin and Zhu \cite{LZ2017} investigated the following improved version of \eqref{a02} under our notations, in which the spreading fronts are explicitly described as free boundaries:
\begin{equation}\left\{
\begin{array}{ll}
u_{1t}=d_1u_{1xx}+a_1(e_1-u_1)u_2 -b_1 u_1, &t>0,\  g(t)<x< h(t), \\[2mm]
u_{2t}=d_2u_{2xx}+a_2(e_2-u_2) u_1 -b_2 u_2,  &t>0,\  g(t)<x< h(t),\\[2mm]
u_1(t,x)=u_2(t,x)=0,  &t>0,\   x\in\{g(t), h(t)\},\\[2mm]
h(0)=h_0, \, h'(t)=-\mu u_{1x}(t, h(t)), &t>0,\\[2mm]
g(0)=-h_0, \, g'(t)=-\mu u_{1x}(t, g(t)), &t>0,\\[2mm]
u_1(0,x)=u_{1,0}(x), u_2(0,x)=u_{2,0}(x),    &-h_0\leq x\leq h_0,
\end{array} \right.
\label{a04}
\end{equation}
where $x=h(t)$ and $x=g(t)$ are the moving boundaries to be determined together with $(u_1, u_2)$, $(g(t), h(t))$ is the infected interval, and the initial functions satisfy
\bess\left\{
\begin{array}{ll}
u_{i,0}(x)\in C^2([-h_0, h_0]), &u_{i,0}(-h_0)=u_{i,0}(h_0)=0,\\[2mm]
0<u_{i,0}(x)\leq e_i, & x\in(-h_0, h_0), i=1, 2.
\end{array} \right.
\eess
It is proved in \cite{LZ2017} that problem \eqref{a04} has a unique solution which is defined for all $t > 0$, and when $R_0\leq 1$, the virus always vanishes eventually, and if $R_0 > 1$, they  proved that spreading-vanishing dichotomy holds. Subsequently, the asymptotic spreading speed of \eqref{a04} was determined in \cite{WND2019} when spreading occurs. Some free boundary problems similar to \eqref{a04} have been extensively studied over the past decade, see \cite{Wang2019, Wang2021, Liu2019} and references therein.

We notice that in both problems \eqref{a02} and \eqref{a04}, the movement of birds and mosquitoes in space is approximately regarded as random diffusion,  which is called local diffusion, expressed as
$d_1u_{1xx}$ and $d_2u_{2xx}$, respectively. However, the classical Laplace operator can not describe all diffusion processes in nature. Murray pointed out in \cite{M1998} that the reaction-diffusion equation in the form of \eqref{a02} or \eqref{a04} can only be used to describe the model with sparse density. But in the embryonic development model, the cell density involved is large, and the Laplace operator can not describe their accurate diffusion process, while the convolution operator
\be
(J*u-u)(t,x):=\int_\mathbb{R}J(x-y)u(t,y)dy-u(t,x)
\label{a11}
\ee
can overcome this problem. We call \eqref{a11}  nonlocal diffusion operator.

In the recent years, the problems with local diffusion have been extensively studied in the literature. However,  it has been increasingly recognized that the nonlocal diffusion as a long-range process can better model some natural phenomena. Cao et al. \cite{Cao2019} introduced
the free boundary model with nonlocal diffusion, which is a natural extension of
a free boundary model with local diffusion in \cite{Du2010}, and they showed that for the spreading-vanishing criteria,
the nonlocal diffusion model has quite different characteristics
in comparison to the corresponding local diffusion model. Subsequently, the spreading speed of the problem in \cite{Cao2019} when the expansion occurs was solved by Du et al. \cite{DL2019}. Inspired by the single species free boundary model with nonlocal diffusion in \cite{Cao2019}, two species nonlocal diffusion systems with free boundaries were investigated; see recent work \cite{DW2019 , Wang20201, Wang20202} for competition  and predator-prey models, and \cite{ZL2020, Zhao2020} for  epidemic models.

Recently, to better describe the dispersal (especially long-range dispersal) of birds and mosquitoes, Du and Nie \cite{Du2020} discussed the nonlocal version of \eqref{a04} as follows:
\be\left\{
\begin{array}{ll}
u_{1t}=d_1\mathcal{L}_1[u_1]+a_1(e_1-u_1)u_2 -b_1 u_1, &t>0,\ g(t)<x< h(t), \\[2mm]
u_{2t}=d_2\mathcal{L}_2[u_2]+a_2(e_2-u_2) u_1 -b_2 u_2,  &t>0,\ g(t)<x< h(t),\\[2mm]
u_1(t,x)=u_2(t,x)=0,  &t>0,\ x\in\{g(t), h(t)\},\\[2mm]
h'(t)=\mu\int_{g(t)}^{h(t)}\int_{h(t)}^{+\infty}J_1(x-y)u_1(t,x)dydx, &t>0,\\[2mm]
g'(t)=-\mu\int_{g(t)}^{h(t)}\int_{-\infty}^{g(t)}J_1(x-y)u_1(t,x)dydx, &t>0,\\[2mm]
u_i(0,x)=u_{i,0}(x),\,i=1,2,   &-h_0\leq x\leq h_0,
\end{array} \right.
\label{a05} \ee
where
\bess
\mathcal{L}_i[u]=\mathcal{L}_i[u;g,h](t,x)=\int_{g(t)}^{h(t)}J_i(x-y)u(t,y)dy-u(t,x).
\eess
They 
showed that problem \eqref{a05} is well-posed and the spreading-vanishing dichotomy holds. At the same time, they gave the criteria when spreading and vanishing can occur.
A significant difference was observed between problem \eqref{a05} and problem \eqref{a04}, which was studied in \cite{LZ2017}.
The most marked difference
is that the spreading may have infinite speed (or accelerated spreading) in the nonlocal model. In fact, Wang et al.\cite{WND2019} showed that the spreading speed of the local model is finite whenever spreading occurs.

On the other hand, as discussed in \cite{HZ}, the temporal changes of the environment affect the growth and development of species, and the environmental variations
due to alternating seasons in nature not only affect the growth of species but also impact on the composition of communities. For example,
in temperate lakes,
phytoplankton and zooplankton have a growing season in warm months, after which species die or  enter the dormant period in winter. To investigate the effects of seasonal succession on the dynamical behaviour of the population, Steiner et al.\cite{S2009} and Hu and Tessier \cite{HT1995} did a lot of experiments to get some data on the effect of seasonal succession on phytoplankton competition, and Klausmeier \cite{K2010} considered a well-known Rosenzweig-McArthur model to study the effect of seasonal alternation on the dynamic behavior of the model. Later, the complete dynamic behavior of a class of Lotka-Volterra competitive models with seasonal succession have been investigated 
by Hsu and Zhao \cite{HZ}, and so on.   Introducing seasonal succession into the diffusion logistic model with a free boundary, Peng and Zhao \cite{PZ2013} studied the following problem
\be\left\{
\begin{array}{ll}
u_{t}=-k u, &m\omega<t\leq m\omega+(1-\tau)\omega, 0<x <h(t), \\[2mm]
u_{t}-du_{xx}=u(a-bu), &m\omega+(1-\tau)\omega<t\leq (m+1)\omega, 0<x <h(t), \\[2mm]
u_x(t,0)=0,  &t>0,\\[2mm]
u(t, h(t))=0,&t>0,\\[2mm]
h(t)=h(m\omega), &m\omega<t\leq m\omega+(1-\tau)\omega,\\[2mm]
h'(t)=-\mu u_x(t, h(t)), &m\omega+(1-\tau)\omega<t\leq (m+1)\omega,\\[2mm]
h(0)=h_0,\\[2mm]
u(0,x)=u_0(x),   &0\leq x\leq h_0,
\end{array} \right.
\label{a06} \ee
where $0<\tau<1 $, and $\omega, k$ are the positive constants. The parameters $\omega$ and $\tau$ account for the period of seasonal succession and the duration ratio of the warm season, respectively.  They gave the criteria of spreading and vanishing, which determined whether the species $u$ spatially spreads to infinity or vanishes in a limited space interval, and
illustrated the influence of the duration of the warm season and cold season on the dynamical behavior. Furthermore, the spreading speed of the species when the spreading occurs was derived.

Inspired by the work of \cite{Cao2019,Du2020,PZ2013}, we consider a West Nile virus nonlocal model with free boundaries and seasonal succession, which reads as follows:
{\footnotesize
\be\left\{
\begin{array}{ll}
u_{1t}=d_1\mathcal{L}_1[u_1]+a_1(e_1-u_1)u_2 -b_1 u_1, &m\omega<t\leq m\omega+(1-\delta)\omega, g(t)<x< h(t), \\[2mm]
u_{2t}=d_2\mathcal{L}_2[u_2]+a_2(e_2-u_2) u_1 -b_2 u_2,  &m\omega<t\leq m\omega+(1-\delta)\omega, g(t)<x< h(t),\\[2mm]
u_{1t}=d_1\mathcal{L}_1[u_1]-b_1 u_1, &m\omega+(1-\delta)\omega<t\leq (m+1)\omega, g(t)<x< h(t), \\[2mm]
u_{2t}=-k u_2,  &m\omega+(1-\delta)\omega<t\leq (m+1)\omega, g(t)<x< h(t), \\[2mm]
u_i(t,x)=0,  &t>0,\ x\in\{g(t), h(t)\}, i=1,2\\[2mm]
h'(t)=\sum\limits_{i=1}^{2}\mu_i\int_{g(t)}^{h(t)}\int_{h(t)}^{+\infty}J_i(x-y)u_i(t,x)dydx, &m\omega<t\leq m\omega+(1-\delta)\omega,\\[2mm]
g'(t)=-\sum\limits_{i=1}^{2}\mu_i\int_{g(t)}^{h(t)}\int_{-\infty}^{g(t)}J_i(x-y)u_i(t,x)dydx, &m\omega<t\leq m\omega+(1-\delta)\omega,\\[2mm]
h(t)=h(m\omega+(1-\delta)\omega), &m\omega+(1-\delta)\omega<t\leq (m+1)\omega,\\[2mm]
g(t)=g(m\omega+(1-\delta)\omega), &m\omega+(1-\delta)\omega<t\leq (m+1)\omega,\\[2mm]
h(0)=-g(0)=h_0,\\[2mm]
u_i(0,x)=u_{i, 0}(x),   &-h_0\leq x\leq h_0,\,i=1,2,
\end{array} \right.
\label{a07}\ee}
where $k\geq b_2$,  the initial time ($t=0$) is chosen as the starting time of the warm season of the first year ($m=0$), the parameter $\omega$ is the length of one year, $(1-\delta)\omega$ accounts for the length of the warm season and the initial function $u_{i, 0}(x)$,
$i=1, 2$, satisfies
\be\left\{
\begin{array}{ll}
u_{i, 0}(x)\in C([-h_0, h_0]), &u_{i, 0}(-h_0)=u_{i, 0}(h_0)=0,\\
0<u_{i, 0}(x)\leq e_i, &  x\in(-h_0, h_0), i=1, 2.
\end{array} \right.
\label{a08} \ee
The kernel function $J_i:\mathbb{R}\rightarrow\mathbb{R}\ (i=1, 2)$ is nonnegative and continuous, and satisfies
$$\mathbf{(J)}: J_i\in C(\mathbb{R})\cap L^{\infty}(\mathbb{R}) \, {\rm \, is \, symmetric}, J_i(0)>0, \int_{\mathbb{R}}J_i(x)dx=1, \sup_{\mathbb{R}}J_i<\infty, i=1, 2.$$

Here and in what follows, unless stated otherwise, we always take $m = 0, 1, 2,\cdots$.

In \eqref{a07}, we take a year as the cycle $[m\omega, (m+1)\omega]$ and divide a year into the warm season and the cold season. From spring to autumn, because of the warm climate and abundant food, birds and mosquitoes have more opportunities to capture food and reproduce, we define this period as the warm season
$(m\omega, m\omega+(1-\delta)\omega]$, and assume that the spatiotemporal distribution and spreading of species are respectively controlled by the first two equations and free boundary conditions in \eqref{a07}. Correspondingly, the cold season $(m\omega+(1-\delta)\omega, (m+1)\omega]$ is from autumn to spring of the next year. Due to formidable survival conditions (the cold climate, shortage of resources, or species can not capture enough food to feed up their offspring), we assume that the number of mosquitoes conforms to Malthusian growth rule which shows that it decays at an exponential rate \cite{Malthus}. During this season, the hypothesis $h(t)=h(m\omega+(1-\delta)\omega)$ together with $g(t)=g(m\omega+(1-\delta)\omega)$ implies biologically that mosquitoes winter in the form of eggs or larvae, and do not migrate in the interval $[g(m\omega+(1-\delta)\omega), h(m\omega+(1-\delta)\omega)]$. Meanwhile, the infected birds can only spread in the above interval. The free boundaries conditions in \eqref{a07}
\bess\left\{
\begin{array}{ll}
h'(t)=\sum\limits_{i=1}^{2}\mu_i\int_{g(t)}^{h(t)}\int_{h(t)}^{+\infty}J_i(x-y)u_i(t,x)dydx, \\[2mm]
g'(t)=-\sum\limits_{i=1}^{2}\mu_i\int_{g(t)}^{h(t)}\int_{-\infty}^{g(t)}J_i(x-y)u_i(t,x)dydx
\end{array} \right.
\eess
mean that the expansion rate of the common population range of the two species (mosquitoes and birds) are directly proportional to the outward flux of the two species, where the factors $\mu_1$ and $\mu_2$ measure the spreading abilities of birds and mosquitoes to the new area.

For the convenience of discussions, we introduce some notations. Given $h_0,  \omega>0$, we denote
{\small\bess
\begin{array}{lll}
&\mathbb{H}_{h_0,\omega}:=\{h\in C([0, +\infty))\cap C^1([m\omega, m\omega+(1-\delta)\omega]): h(0)=h_0, h(t) \,{\rm is \, nondecreasing}\}, \\[2mm]
&\mathbb{G}_{h_0,\omega}:=\{g\in C([0, +\infty))\cap C^1([m\omega, m\omega+(1-\delta)\omega]): -g\in \mathbb{H}_{h_0,\omega}\}.
\end{array}
\eess}
For any given $g\in \mathbb{G}_{h_0,\omega}, h\in \mathbb{H}_{h_0,\omega}$ and $U_0=(u_{1,0},\, u_{2,0})$ satisfying \eqref{a08}, set
$$\Omega^{g, h}:=\{(t ,x)\in \mathbb{R}^2: t\in [0, +\infty), g(t)<x<h(t)\},$$
\begin{equation*}
\begin{array}{ll}
\mathbb{X}_\omega=\mathbb{X}_{\omega, U_0}^{g, h}:=\{(\xi_1, \xi_2): &\xi_i\in C(\Omega^{g, h}),\, \xi_i\geq0,\, \xi_i(0, x)=u_{i, 0}(x) \,\, {\rm for }\,\,x\in[-h_0, h_0], \\[2mm]
&{\rm and}\,\,\xi_i(t, g(t))=\xi_i(t, h(t))=0 \,\, {\rm for }\,\,t\in[0, +\infty), i=1, 2\}.
\end{array}
\end{equation*}
Our main results for \eqref{a07} are given in the following.
\begin{thm}
\label{thm1.1}
$($Existence and uniqueness$)$
Suppose $\mathbf{(J)}$ holds, and the initial functions satisfy \eqref{a08}. Then problem \eqref{a07} admits a unique positive solution $(u_1(t, x), u_2(t, x); g(t), h(t))$ satisfying $(g, h)\in\mathbb{G}_{h_0,\omega}\times\mathbb{H}_{h_0,\omega}, (u_1, u_2)\in\mathbb{X}_{\omega}$ for any given $\omega>0$. Furthermore,
\be\left\{
\begin{array}{ll}
0<u_i\leq e_i \, \, &{\rm for} \,\, t>0, x\in[g(t), h(t)],\\[2mm]
g'(t)\leq0, h'(t)\geq0 &{\rm for} \,\, t>0.
\end{array} \right.
\ee
\end{thm}

\begin{thm}
\label{thm1.2}
$($Spreading-vanishing dichotomy$)$
Suppose $\mathbf{(J)}$ holds, and the initial functions satisfy \eqref{a08}. Let $(u_1(t, x), u_2(t, x); g(t), h(t))$ be the solution to problem \eqref{a07} and denote
$$g_{\infty}:=\lim\limits_{t\rightarrow\infty}g(t) \quad {\rm and} \quad h_{\infty}:=\lim\limits_{t\rightarrow\infty}h(t).$$
Then one of the following alternatives must happen for \eqref{a07}:

$(i)$ Spreading:\  $-g_{\infty}=h_{\infty}=\infty$ and
$$\lim\limits_{n\rightarrow\infty}(u_1(t+n\omega,x), u_2(t+n\omega,x))= (U_1^{\vartriangle}, U_2^{\vartriangle})$$
 uniformly for $t\in [0, \omega]$ and locally uniformly for $x\in\mathbb{R}$, where $(U_1^{\vartriangle}, U_2^{\vartriangle})(t)$ is the unique positive periodic solutions of the following problem
{\small
\bess
\left\{
\begin{array}{ll}
U_{1t}=a_1(e_1-U_1)U_2 -b_1 U_1, &t\in (0, (1-\delta)\omega], \\[2mm]
U_{2t}=a_2(e_2-U_2) U_1 -b_2 U_2,  &t\in (0, (1-\delta)\omega],\\[2mm]
U_{1t}= -b_1
U_1, &t\in ((1-\delta)\omega, \omega], \\[2mm]
U_{2t}=-k_2 U_2, &t\in ((1-\delta)\omega, \omega], \\[2mm]
U_i(0)=U_i(\omega)(i=1,2).
\end{array} \right.
\eess
}

$(ii)$ Vanishing: $$\lim\limits_{t\rightarrow\infty}(u_1(t,x), u_2(t,x))=(0, 0)\, \,  {\rm uniformly \, \,  for} \, \,  x\in[g(t), h(t)].$$
\end{thm}

We mention here that $h_\infty-g_\infty<\infty$ when vanishing happens in some special cases, see Lemma \ref{vanish}.

\begin{thm}
\label{thm1.3}
$($Spreading-vanishing criteria 1$)$
Suppose $\mathbf{(J)}$, $J_1(x)=J_2(x)$ and $0\leq\delta<1$ hold, and the initial functions satisfy \eqref{a08}. Fix $a_i, e_i, b_2$ and $k$,
then the following conclusions hold.

$(i)$ If $b_1\geq b_1^*:=a_1a_2e_1e_2/b_2$, then vanishing happens; while for $b_1\in(0, b_1^*)$, there exists $\delta^*\in (0, 1)$ such that vanishing happens for $\delta\in[\delta^*, 1]$.

$(ii)$ If $b_1\in(0, b_1^*)$ and $\delta\in(0, \delta^*)$, there exists $h_0^*>0$ such that spreading happens for $h_0\geq h_0^*$.

$(iii)$ If $b_1\in(0, b_1^*)$, $\delta\in(0, \delta^*)$ and $h_0\in (0, h_0^*)$, then
for any given initial datum $(u_{1, 0}, u_{2, 0})$ satisfying \eqref{a08}, there exists $\mu^*\geq\mu_*>0$ such that vanishing happens for $0<\mu_1+\mu_2\leq \mu_*$ and spreading happens for $\mu_1+\mu_2>\mu^*$.
\end{thm}

Theorem \ref{thm1.3} implies that, spreading-vanishing of virus depends firstly on the recovery rate ($b_1$ or $\gamma_b$) of birds from WNv,
large recovery rate is benefit for the vanishing of virus. For small recovery rate, large length ($\delta\omega$) of the cold season leads to
the vanishing of virus. However, even the recovery rate and the length of the cold season are small, large initial infected interval ($[-h_0, h_0]$) or
the expanding factors ($\mu_1$ and $\mu_2$) can bring about the spreading of WNv.

The following result presents the sufficient conditions for the spreading or vanishing of WNv in the case $J_1(x)\neq J_2(x)$.
\begin{thm}
\label{thm1.4}
$($Spreading-vanishing criteria 2$)$
Suppose $\mathbf{(J)}$ holds, and the initial functions satisfy \eqref{a08}. There exists $0\leq \mu_\vartriangle\leq\mu^\vartriangle\leq +\infty$ such that vanishing happens for $(0, 0)<(\mu_1, \mu_2)< (\mu_\vartriangle, \mu_\vartriangle)$, and spreading happens for $(\mu_1, \mu_2)>(\mu^\vartriangle, \mu^\vartriangle)$.
In particular, if $\delta=1$, or $0\leq\delta<1$ and $\min\{b_1, k\}\delta>c_1(1-\delta)$, then $\mu_\vartriangle=\mu^\vartriangle=+\infty$, which means that vanishing happens for any $\mu_1$ and $\mu_2$, where $$c_1=\frac{-(b_1+b_2)+\sqrt{(b_1-b_2)^2+4a_1a_2e_1e_2}}{2}.$$
While if $0\leq\delta<1$ and $\max\{b_1, k\}\delta<c_1(1-\delta)$, then $\mu_\vartriangle=\mu^\vartriangle=0$, which means that spreading happens for any $\mu_1$ and $\mu_2$, provided that $h_0$ is sufficiently large.
\end{thm}

The rest of the paper is organized as follows: In Section 2, we introduce the associated eigenvalue problem and present the detailed analysis of principal eigenvalue in a bounded domain.
A fixed boundary problem with seasonal succession and the corresponding preliminaries are addressed
in Section 3. These results pave the way for proving the main results
in Section 4. While some of the ideas
are adopted from \cite{Du2020, PZ2013}, considerable variations are needed as our principal problem here does not satisfy some conditions required in \cite{B2017}, which is different from the problem in \cite{Du2020}. Section 4 is devoted to the proofs of Theorems \ref{thm1.2}, \ref{thm1.3} and \ref{thm1.4}.
Some ecological explanations of our theoretical results are given in the last section.

\section{An associated eigenvalue problem}

To understand the asymptotical behavior of the solution to problem \eqref{a07}, we
study the corresponding periodic eigenvalue problem
{\small
\be\left\{
\begin{array}{ll}
\phi_t -d_1\mathcal{L}_1[\phi]=a_1e_1\psi -b_1 \phi+\lambda\phi, &0<t\leq (1-\delta)\omega, -L_1\leq x\leq L_2, \\[2mm]
\psi_t-d_2\mathcal{L}_2[\psi]=a_2e_2\phi -b_2 \psi+\lambda\psi,  &0<t\leq (1-\delta)\omega, -L_1\leq x\leq L_2,\\[2mm]
\phi_t-d_1\mathcal{L}_1[\phi]=-b_1 \phi+\lambda\phi, &(1-\delta)\omega<t\leq \omega, -L_1\leq x\leq L_2, \\[2mm]
\psi_t=-k\psi+\lambda\psi,  &(1-\delta)\omega<t\leq \omega, -L_1\leq x\leq L_2,\\[2mm]
\phi(0,x)=\phi(\omega,x), \,\,\psi(0,x)=\psi(\omega,x),  &-L_1\leq x\leq L_2.
\end{array} \right.
\label{b01}\ee
}
In order to get some properties of periodic eigenvalue problem \eqref{b01}, we turn to investigate its associated spatial-independent version. More generally, we consider the following problem
{\small
\be\left\{
\begin{array}{ll}
\phi_t =a_1e_1\psi -b_1 \phi+\lambda \phi, &0<t\leq (1-\delta)\omega, \\[2mm]
\psi_t=a_2e_2\phi -b_2 \psi+\lambda \psi,  &0<t\leq (1-\delta)\omega,\\[2mm]
\phi_t=-k_1 \phi+\lambda \phi, &(1-\delta)\omega<t\leq \omega, \\[2mm]
\psi_t=-k_2\psi+\lambda \psi,  &(1-\delta)\omega<t\leq \omega,\\[2mm]
\phi(0)=\phi(\omega), \,\,\psi(0)=\psi(\omega).
\end{array} \right.
\label{b13}\ee
}

As in \cite{N2009}, we define the generalized principal eigenvalues $\overline{\lambda}_1^{O}$ and $\underline{\lambda}_{1}^{O}$ of problem \eqref{b13} as
\bess
\begin{array}{ll}
\overline{\lambda}_1^{O}:=\inf\{\lambda\in \mathbb{R}\,|\, \exists\phi, \psi\in C^{1}(\mathbb{R}), \phi, \psi>0 {\rm\  and\ }  \phi, \psi {\ \rm are \ } \omega-{\rm periodic \ so \ as}\\
\qquad \qquad\ \ {\rm \ the\ inequalities\ in\  }\eqref{b19} {\rm \ hold}\},
\end{array}
\eess
where
\be\left\{
\begin{array}{ll}
\phi_t \leq a_1e_1\psi -b_1 \phi+\lambda \phi, &0<t\leq (1-\delta)\omega, \\[2mm]
\psi_t \leq a_2e_2\phi -b_2 \psi+\lambda \psi,  &0<t\leq (1-\delta)\omega,\\[2mm]
\phi_t \leq -k_1 \phi+\lambda\phi, &(1-\delta)\omega<t\leq \omega, \\[2mm]
\psi_t \leq -k_2\psi+\lambda \psi,  &(1-\delta)\omega<t\leq \omega,
\end{array} \right.
\label{b19}\ee
and
\bess
\begin{array}{ll}
\underline{\lambda}_{1}^{O}:=\sup\{\lambda\in \mathbb{R}\,|\, \exists\phi, \psi\in C^{1}(\mathbb{R}), \phi, \psi>0 {\rm\  and\ }  \phi, \psi {\ \rm are \ } \omega-{\rm periodic\ } {\rm so \ as\  } \\
\qquad \qquad\ \ \eqref{b19}\ \textrm{ with  the inequalities revised  holds}\}.
\end{array}
\eess
If there exists a  positive function pair $(\phi, \psi)\in C^{1}(\mathbb{R})\times C^{1}(\mathbb{R})$ such that problem \eqref{b13} with $\lambda=\lambda_1^O$ holds, then
$\lambda_1^O$ is called a principal eigenvalue of problem \eqref{b13}.

We first see the special case $\delta=1$. In this case, the spatial-independent problem \eqref{b13} becomes
\bess\left\{
\begin{array}{ll}
\phi_t=-k_1 \phi+\lambda \phi, &0<t\leq \omega, \\[2mm]
\psi_t=-k_2\psi+\lambda \psi,  &0<t\leq \omega,\\[2mm]
\phi(0)=\phi(\omega), \,\,\psi(0)=\psi(\omega).
\end{array} \right.
\eess
If $k_1=k_2$, then $\overline{\lambda}_1^O=\underline{\lambda}_1^O=\lambda_1^O=k_1,$ and if $k_1\neq k_2$, then
$\overline{\lambda}_1^O=\max\{k_1, k_2\}$, $\underline{\lambda}_1^O=\min\{k_1, k_2\}$ and the two  generalized principal eigenvalues are not equal obviously.

Suppose that $0\leq\delta<1$, we claim that $\overline{\lambda}_1^{O}=\underline{\lambda}_{1}^{O}=\lambda_1^{O}$, where $\lambda_1^{O}$ is the principal eigenvalue of problem \eqref{b13} with a positive eigenfunction pair $(\phi, \psi)\in [C^1([0, \omega])]^2$. In fact,
we can further provide a detailed calculation process of principal eigen-pair $(\lambda_1^{O}, \phi, \psi)$ in the following.

With a view to the periodicity of $\phi$ and  $\psi$, problem \eqref{b13} can be rewritten as
\be\left\{
\begin{array}{ll}
\phi_t =a_1e_1\psi -b_1 \phi+\lambda_1^{O}\phi, &0<t\leq (1-\delta)\omega, \\[2mm]
\psi_t=a_2e_2\phi -b_2 \psi+\lambda_1^{O}\psi,  &0<t\leq (1-\delta)\omega,\\[2mm]
\phi((1-\delta)\omega)=\phi(0)e^{(k_1-\lambda_1^{O})\delta\omega},  \\[2mm]
\psi((1-\delta)\omega)=\phi(0)e^{(k_2-\lambda_1^{O})\delta\omega},
\end{array} \right.
\label{b14}\ee
the first two equations in \eqref{b14} are abbreviated as
\be
\begin{pmatrix}
 \phi_t\\
  \psi_t
 \end{pmatrix}=
\begin{pmatrix}
 -b_1+\lambda_1^{O} & a_1e_1 \\
  a_2e_2 & -b_2+\lambda_1^{O}
 \end{pmatrix}
 \begin{pmatrix}
 \phi\\
  \psi
 \end{pmatrix}:=M\begin{pmatrix}
 \phi\\
  \psi
 \end{pmatrix}.
\label{b15}
\ee
Then, we see from the characteristic equation $|M-\mu E|=0$ that
\be
\mu_{1,2}=\lambda_1^{O}+\frac{-(b_1+b_2)\pm \sqrt{(b_1-b_2)^2+4a_1a_2e_1e_2}}{2}:=\lambda_1^{O}+c_{1,2}.
\label{b26}\ee
Without loss of generality, we can assume that $c_1\geq c_2$. Direct calculation yields that $ b_1+c_1=-(b_2+c_2)>0$.

The linearly independent eigenvectors $(k_{11}, k_{12})$ and $(k_{21}, k_{22})$ associated with eigenvalues $\mu_1$ and $\mu_2$ satisfy
\be
\begin{pmatrix}
 k_{i1}& k_{i2}
 \end{pmatrix}
\begin{pmatrix}
 -b_1+\lambda_1^{O}-\mu_i & a_1e_1 \\
  a_2e_2 & -b_2+\lambda_1^{O}-\mu_i
 \end{pmatrix}=\begin{pmatrix}
 0 & 0
 \end{pmatrix}
\label{e08}\ee
for $i=1,2$, and it is easily seen that
$$(k_{11}, k_{12})=(a_2e_2, b_1-\lambda_1^{O}+\mu_1)=(a_2e_2, b_1+c_1)$$
and
$$(k_{21}, k_{22})=(b_2-\lambda_1^{O}+\mu_2, a_1e_1)=(b_2+c_2, a_1e_1).$$

Subsequently,  we consider the following algebraic equations
\be
\begin{pmatrix}
 a_2e_2&  b_1+c_1 \\
  b_2+c_2 & a_1e_1
 \end{pmatrix}
\begin{pmatrix}
 \phi\\
  \psi
 \end{pmatrix}=
 \begin{pmatrix}
 e^{\mu_1t}\\
  me^{\mu_2t}
 \end{pmatrix},
 \label{e09}\ee
its solution is naturally given by
\bess
(\phi,\, \psi)=\bigg(\frac{a_1e_1e^{\mu_1t}-(b_1+c_1)me^{\mu_2t}}{C_0},\, \frac{-(b_2+c_2)e^{\mu_1t}+a_2e_2me^{\mu_2t}}{C_0}\bigg),
\eess
where
\be
C_0=a_1a_2e_1e_2-(b_1+c_1)(b_2+c_2)=a_1a_2e_1e_2+(b_1+c_1)^2>0.
\label{cc}
\ee
Using the third and forth equations of \eqref{b14}, we have
{\small
\be
\left\{
\begin{array}{ll}
a_1e_1e^{c_1(1-\delta)\omega}e^{\lambda_1^{O}\omega}-(b_1+c_1)e^{c_2(1-\delta)\omega}m e^{\lambda_1^{O}\omega}
+(b_1+c_1)e^{k_1\delta\omega}m=
a_1e_1e^{k_1\delta\omega},\\[2mm]
-(b_2+c_2)e^{c_1(1-\delta)\omega}e^{\lambda_1^{O}\omega}+a_2e_2e^{c_2(1-\delta)\omega}m e^{\lambda_1^{O}\omega}
-a_2e_2e^{k_2\delta\omega}m=
-(b_2+c_2)e^{k_2\delta\omega}.
\end{array} \label{b16}\right.
\ee}
For convenience, we denote $\Lambda=e^{\lambda_1^{O}\omega}$. Next, we show that problem \eqref{b16}  admits a unique solution $(m, \Lambda)$ such that $\phi>0$ and $\psi>0$.
For this purpose, we need consider the following three cases.

\begin{figure}[ht]
\centering
\includegraphics[width=0.45\textwidth]{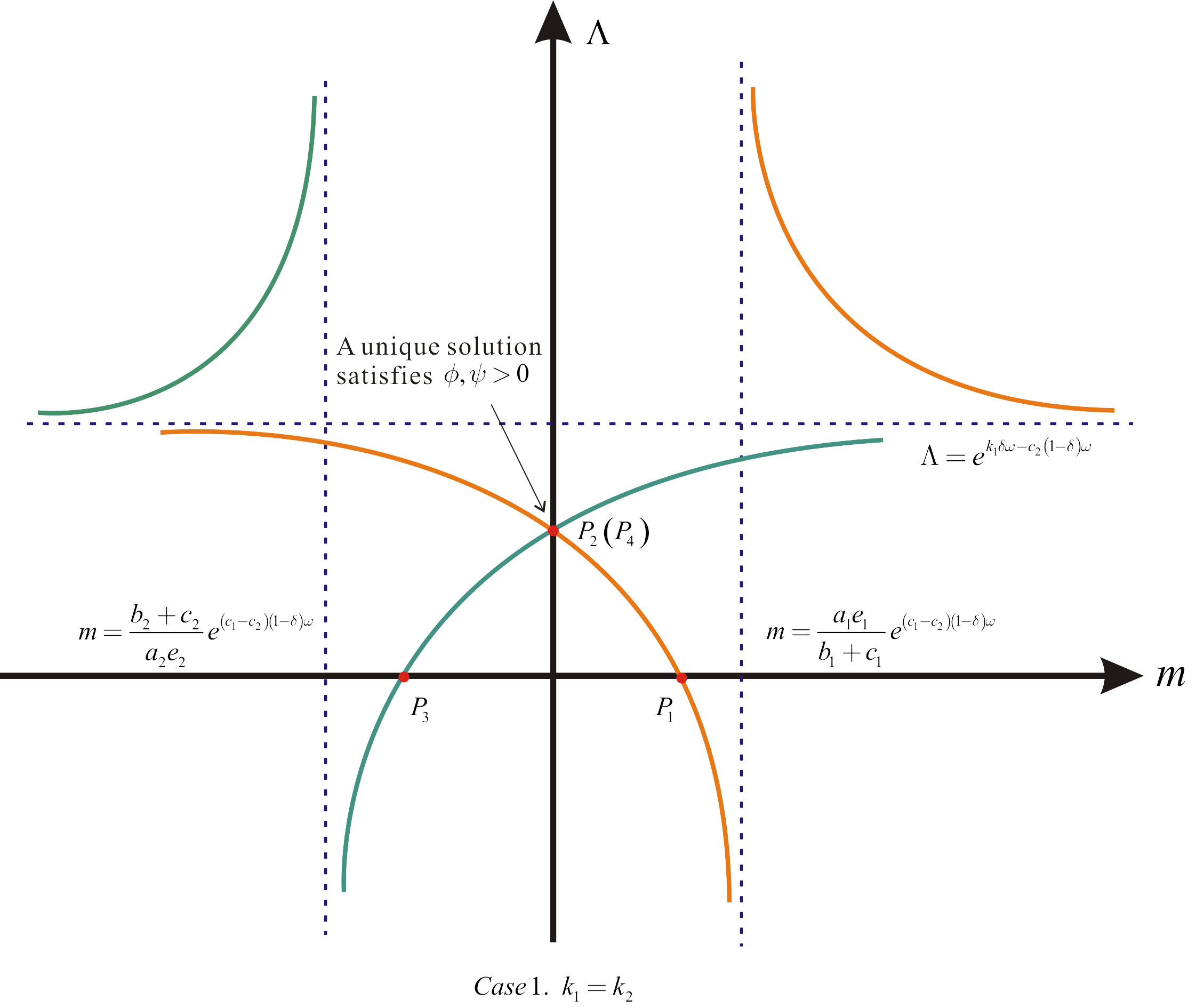}
\caption{\scriptsize Figure shows the distribution of the unique solution in Case 1. }
\label{tu1}
\end{figure}

\textbf{Case 1.} If $k_1=k_2$, multiplying both sides of the first equation of \eqref{b16} by $a_2e_2$,
and the second equation of that by $(b_1+c_1)$, then adding the two equations, yields
$$[a_1a_2e_1e_2-(b_1+c_1)(b_2+c_2)][e^{c_1(1-\delta)\omega}\Lambda
-e^{k_1\delta\omega}]=0,$$
hence problem \eqref{b16} has a unique solution $(0, e^{k_1\delta\omega-c_1(1-\delta)\omega})$ since that $[a_1a_2e_1e_2-(b_1+c_1)(b_2+c_2)]=[a_1a_2e_1e_2+(b_1+c_1)^2]>0$, see also Fig. \ref{tu1}. The principal eigenvalue $\lambda_1^{O}$ can be explicitly expressed, that is, $\lambda_1^{O}=(k_1+c_1)\delta-c_1$. Moreover, the corresponding eigenfunction pair is
\bess
(\phi,\, \psi)=\bigg(\frac{a_1e_1e^{(k_1+c_1)\delta t}}{C_0},\, \frac{-(b_2+c_2)e^{(k_1+c_1)\delta t}}{C_0}\bigg).
\eess

\begin{figure}[ht]
\centering
\subfigure[]{ {
\includegraphics[width=0.45\textwidth]{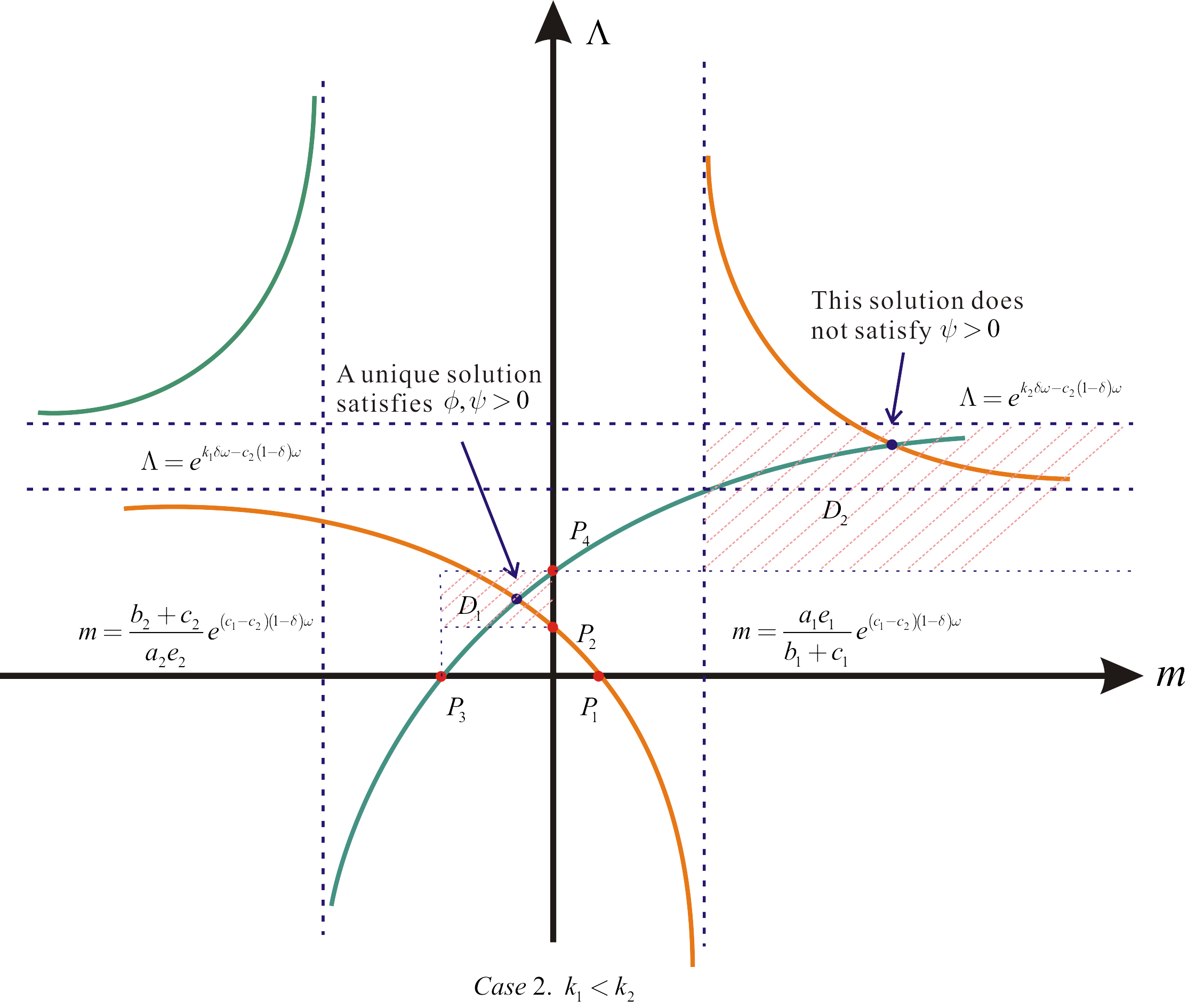}
} }
\subfigure[]{ {
\includegraphics[width=0.45\textwidth]{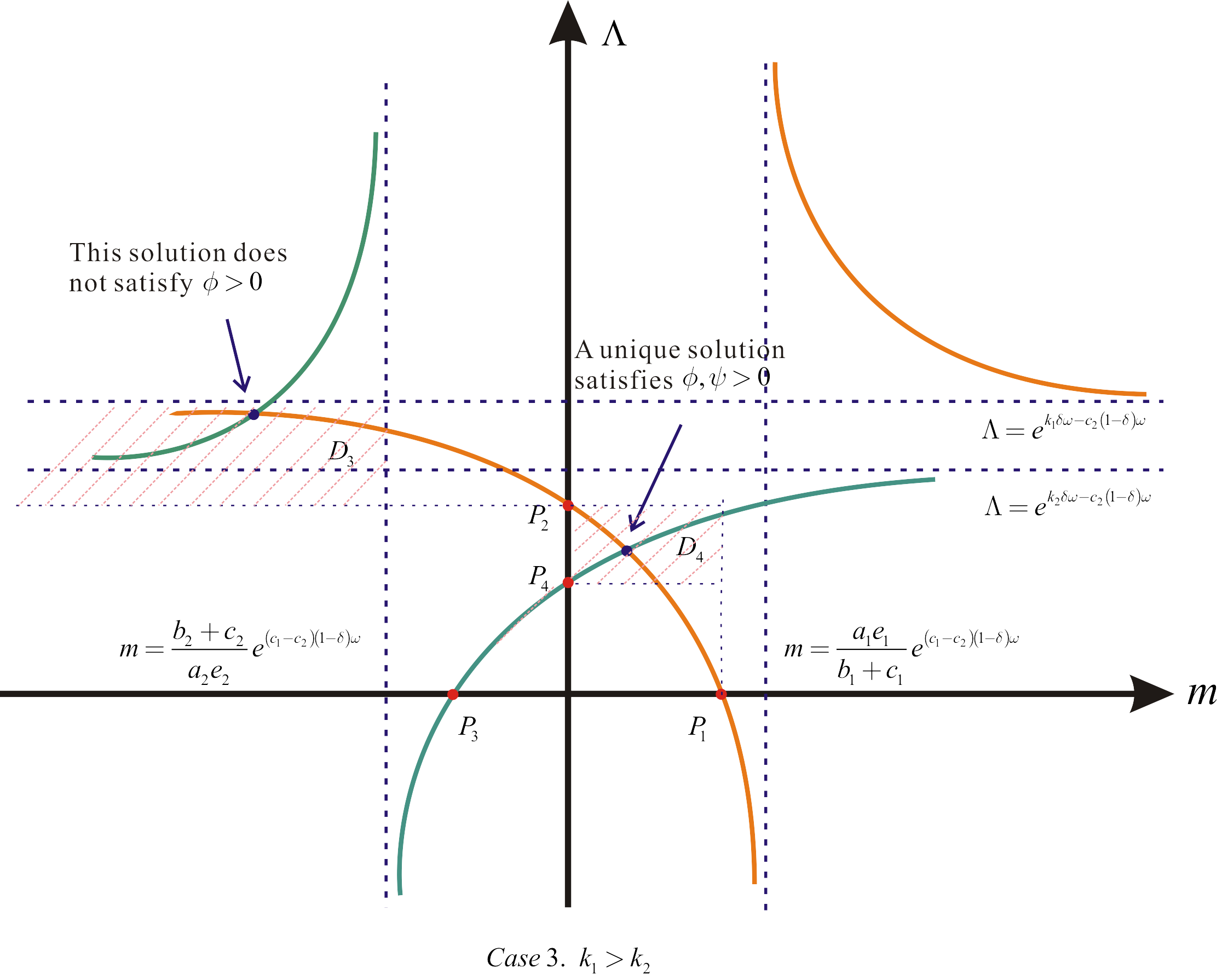}
} }
\caption{\scriptsize Figures (a)-(b) show the distribution of the unique solution in Cases 2 and 3, respectively. }
\label{tu2}
\end{figure}

If $k_1\neq k_2$, for abbreviation, we denote
{\small
\bess
\begin{array}{ll}
A_{11}=a_1e_1e^{c_1(1-\delta)\omega},\, A_{12}=(b_1+c_1)e^{c_2(1-\delta)\omega}, \, A_{13}=(b_1+c_1)e^{k_1\delta\omega},\,A_{14}=a_1e_1e^{k_1\delta\omega}, \\[2mm]
A_{21}=-(b_2+c_2)e^{c_1(1-\delta)\omega},\, A_{22}=a_2e_2e^{c_2(1-\delta)\omega},\, A_{23}=a_2e_2e^{k_2\delta\omega},\,A_{24}=-(b_2+c_2)e^{k_2\delta\omega}.
\end{array}
\eess}
Thus \eqref{b16} becomes
{\small
\be
\left\{
\begin{array}{ll}
A_{11}\Lambda-A_{12}\Lambda m+A_{13}m=A_{14},\\[2mm]
A_{21}\Lambda+A_{22}\Lambda m-A_{23}m=A_{24},
\end{array} \label{b17}\right.
\ee}
where $A_{ij}$ are all positive since that $b_1+c_1=-(b_2+c_2)>0$.
By using the elimination method, we can get the quadratic equation with variable $\Lambda$,
{\small
\be
(A_{11}A_{22}+A_{12}A_{21})\Lambda^2-(A_{13}A_{21}+A_{11}A_{23}-A_{22}A_{14}-A_{12}A_{24})
\Lambda+A_{23}A_{14}+A_{13}A_{24}=0,
\label{b20}\ee}
and that with variable $m$,
{\small
\be
(A_{13}A_{22}-A_{12}A_{23})m^2-(A_{22}A_{14}+A_{12}A_{24}-A_{13}A_{21}-A_{11}A_{23})m+
A_{11}A_{24}-A_{21}A_{14}=0.
\label{b21}\ee
}
It follows from Vieta's theorem that two roots $\Lambda_1, \Lambda_2$ of \eqref{b20} satisfy $\Lambda_1\Lambda_2>0$, and two roots $m_1, m_2$ of \eqref{b21} satisfy $$
m_1m_2=\frac {A_{11}A_{24}-A_{21}A_{14}}{A_{13}A_{22}-A_{12}A_{23}}
=\frac{a_1e_1(b_2+c_2)e^{c_1(1-\delta)\omega}}{a_2e_2(b_1+c_1)e^{c_2(1-\delta)\omega}}
=\frac{-a_1e_1e^{c_1(1-\delta)\omega}}{a_2e_2e^{c_2(1-\delta)\omega}} <-\frac{a_1e_1}{a_2e_2}<0.
$$
It is not easy to derive the explicit solution of problem \eqref{b17}, and we therefore consider the solutions by image method.
In fact, problem \eqref{b17} can be rewritten as
{\small
\be
\left\{
\begin{array}{ll}
\Lambda=(A_{14}- A_{13}m)/(A_{11}-A_{12}m),\\[2mm]
\Lambda=(A_{24}+ A_{23}m)/(A_{21}+A_{22}m).
\end{array} \label{b171}\right.
\ee}
The first equation in \eqref{b171} shows that its curve must go through points $P_1(\frac{a_1e_1}{b_1+c_1}, 0)$ and $P_2(0, e^{k_1\delta\omega-c_1(1-\delta)\omega})$, what's more, $\Lambda$ is strictly decreasing with respect to $m$ in  the interval $(-\infty, \frac{a_1e_1}{b_1+c_1}e^{(c_1-c_2)(1-\delta)\omega})$ $\cup(\frac{a_1e_1}{b_1+c_1}e^{(c_1-c_2)(1-\delta)\omega}, +\infty)$, see Fig. \ref{tu2}. Meanwhile, the second equation indicates that $\Lambda$ is strictly increasing on $m$ in  $(-\infty, \frac{b_2+c_2}{a_2e_2}e^{(c_1-c_2)(1-\delta)\omega})\cup(\frac{b_2+c_2}{a_2e_2}e^{(c_1-c_2) (1-\delta)\omega}, +\infty)$, and its curve pass the points $P_3(\frac{b_2+c_2}{a_2e_2}, 0)$ and $P_4(0, e^{k_2\delta\omega-c_1(1-\delta)\omega})$.

\vspace{1mm}
\textbf{Case 2.} If $k_1<k_2,$ two curves have one intersection in the regions $D_1:= (\frac{b_2+c_2}{a_2e_2}, 0)\times(e^{k_1\delta\omega-c_1(1-\delta)\omega}, $  $e^{k_2\delta\omega-c_1(1-\delta)\omega})$ and $D_2:=(\frac{a_1e_1}{b_1+c_1}e^{(c_1-c_2)(1-\delta)\omega}, +\infty)\times(e^{k_2\delta\omega-c_1(1-\delta)\omega}, e^{k_2\delta\omega-c_2(1-\delta)\omega})$ respectively,
that is to say, problem \eqref{b16} have solutions $(m_1, \Lambda_1)\in D_1$ and $(m_2, \Lambda_2)\in D_2$, see Fig. \ref{tu2} (a).
However, because of $\psi<0$ in  $D_2$, which fails to meet the requirement of the eigenfunction pair, so problem \eqref{b16} has a unique solution, which is in $D_1$, such that $\phi>0$ and $\psi>0$.

\vspace{1mm}
\textbf{Case 3.} If $k_1>k_2,$ the discussion is similar to the  above situation.
We obtain that  problem \eqref{b16} have solutions  $(m_3, \Lambda_3)\in D_3:=(-\infty, \frac{b_2+c_2}{a_2e_2}e^{(c_1-c_2)(1-\delta)\omega})\times( e^{k_1\delta\omega-c_1(1-\delta)\omega}, e^{k_1\delta\omega-c_2(1-\delta)\omega})$ and $(m_4, \Lambda_4)\in D_4:=(0, \frac{a_1e_1}{b_1+c_1})\times(e^{k_2\delta\omega-c_1(1-\delta)\omega}, e^{k_1\delta\omega-c_1(1-\delta)\omega})$, see Fig. \ref{tu2} (b).
Noticing that $\phi<0$ in  $D_3$, so problem \eqref{b16} has a unique solution in $D_4$ such that $\phi>0$ and $\psi>0$.

In summary, if $0\leq\delta<1$, problem \eqref{b16} (i.e. problem \eqref{b13}) has a unique solution with $\phi, \psi>0$.\epf
\vspace{3mm}

The above discussion enables us to obtain the following conclusions.
\begin{thm}
\label{thm2.2}
$(i)$ Assume that
$$\mathbf{(H_1)}:\ (a)\ 0\leq \delta<1\ \textrm{or}\ (b)\ \delta=1, \, k_1=k_2$$
holds. Then $\overline{\lambda}_1^O=\underline{\lambda}_1^O=\lambda_1^O$, the eigenvalue
problem \eqref{b13} admits a principal
eigenvalue $\lambda_1^{O}$ with a positive eigenfunction pair $(\phi, \psi) \in [C^1([0, \omega])]^2$.
Especially, if $k_1=k_2$, for any given $\delta\in[0,1]$, $\lambda_1^O$ can be explicitly expressed as $\lambda_1^O=(k_1+c_1)\delta-c_1.$

$(ii)$ If $\delta=1$ and $k_1\neq k_2$, then
$\overline{\lambda}_1^O=\max\{k_1, k_2\}>\min\{k_1, k_2\}=\underline{\lambda}_1^O.$
So the generalized principal eigenvalues $\overline{\lambda}_1^O$ and $\underline{\lambda}_1^O$ of the eigenvalue problem \eqref{b13} are not equal.
\end{thm}

As in \cite{DHM},
$\lambda_1^O$ is an important threshold of epidemiology models. $\lambda_1^O<0 (>0)$ implies high (low) risk.
It follows from Theorem \ref{thm2.2} that if $\delta=1$, then $\overline{\lambda}_1^O\geq \underline{\lambda}_1^O>0$. We then come to the following conclusion for $0\leq \delta<1$.
\begin{thm}
\label{thm2.3}\
The following statements  are valid:

$(i)$\ If $\delta=0$,  the necessary and sufficient condition for $\lambda_1^{O}=0$   is that $a_1a_2e_1e_2= b_1b_2$.

$(ii)$\ If $0<\delta<1$,
then $\lambda_1^{O}=0$ if and only if
\be
\frac{a_1e_1(e^{c_1(1-\delta)\omega}-e^{k_1\delta\omega})}
{(b_1+c_1)(e^{c_2(1-\delta)\omega}-e^{k_1\delta\omega})}=
\frac{(b_2+c_2)(e^{c_1(1-\delta)\omega}-e^{k_2\delta\omega})}{a_2e_2(e^{c_2(1-\delta)\omega}
-e^{k_2\delta\omega})},
\label{e06}
\ee
where $c_i(i=1, 2)$ is defined in \eqref{b26}.
\end{thm}
\bpf
$(i)$\ We first prove the necessity.
Suppose that $\delta=0$ and $\lambda_1^{O}=0$, then \eqref{b13} becomes
{\small
\be\left\{
\begin{array}{ll}
\phi_t =a_1e_1\psi -b_1 \phi, &0<t\leq \omega, \\[2mm]
\psi_t=a_2e_2\phi -b_2 \psi,  &0<t\leq \omega,\\[2mm]
\phi(0)=\phi(\omega), \,\,\psi(0)=\psi(\omega).
\end{array} \right.
\label{e07}\ee}
Repeat the calculation process from \eqref{b15} to \eqref{e09} with $\lambda_1^O=0$, then the solution of problem \eqref{e07} is
\bess
(\phi,\, \psi)=\bigg(\frac{a_1e_1e^{c_1t}-(b_1+c_1)me^{c_2t}}{C_0},\, \frac{-(b_2+c_2)e^{c_1t}+a_2e_2me^{c_2t}}{C_0}\bigg),
\eess
where $C_0$ is defined in \eqref{cc}.

Using the periodic conditions $\phi(0)=\phi(\omega), \psi(0)=\psi(\omega)$ yields
\bess
\left\{
\begin{array}{ll}
a_1e_1e^{c_1\omega}-(b_1+c_1)e^{c_2\omega}m +(b_1+c_1)m=a_1e_1,\\[2mm]
-(b_2+c_2)e^{c_1\omega}+a_2e_2e^{c_2\omega}m-a_2e_2m=-(b_2+c_2).
\end{array} \right.
\eess
We further have
$$[a_1a_2e_1e_2-(b_1+c_1)(b_2+c_2)][e^{c_1\omega}-1]=0,$$
therefore, $e^{c_1\omega}=1$, that is, $c_1=0$, which means that $a_1a_2e_1e_2=b_1b_2$ by using the expression of $c_1$ in \eqref{b26}.

We now prove the sufficiency.
Suppose that $\delta=0$, \eqref{b13} becomes
{\small
\be\left\{
\begin{array}{ll}
\phi_t =a_1e_1\psi -b_1 \phi+\lambda_1^O\phi, &0<t\leq \omega, \\[2mm]
\psi_t=a_2e_2\phi -b_2 \psi+\lambda_1^O\psi,  &0<t\leq \omega,\\[2mm]
\phi(0)=\phi(\omega), \,\,\psi(0)=\psi(\omega).
\end{array} \right.
\label{e10}\ee
}
Owing to $a_1a_2e_1e_2=b_1b_2$, we have
$$c_1=\frac{-(b_1+b_2)+\sqrt{(b_1-b_2)^2+4a_1a_2e_1e_2}}{2}=0,$$
$$c_2=\frac{-(b_1+b_2)-\sqrt{(b_1-b_2)^2+4a_1a_2e_1e_2}}{2}=-(b_1+b_2),$$
 and the general solution of problem \eqref{e10} is
\bess
(\phi,\, \psi)=\bigg(\frac{a_1e_1e^{\lambda_1^Ot}-b_1me^{[\lambda_1^O-(b_1+b_2)]t}}{C_0},\, \frac{-b_1e^{\lambda_1^O t}+a_2e_2me^{[\lambda_1^O-(b_1+b_2)]t}}{C_0}\bigg).
\eess
Using the periodicity of $\phi$ and $\psi$ yields $\lambda_1^O=0$ and $m=0$.

$(ii)$\ If $0<\delta<1$, it can be calculated directly from \eqref{b16} that
\be
\frac{a_1e_1(e^{c_1(1-\delta)\omega}e^{\lambda_1^O\omega}-e^{k_1\delta\omega})}
{(b_1+c_1)(e^{c_2(1-\delta)\omega}e^{\lambda_1^O\omega}-e^{k_1\delta\omega})}=
\frac{(b_2+c_2)(e^{c_1(1-\delta)\omega}e^{\lambda_1^O\omega}-e^{k_2\delta\omega})}
{a_2e_2(e^{c_2(1-\delta)\omega}e^{\lambda_1^O\omega}-e^{k_2\delta\omega})}.
\label{e11}\ee
therefore, $\lambda_1^O=0$ if and only if \eqref{e06} holds.
\epf

\vspace{3mm}
Now, we give conditions for $\lambda_1^{O}>(<)0$, which means that the risk is low(high).
\begin{cor}
\label{cor1.3}\
The following statements are valid:

$(i)$\ If $\delta=0$,  then $\lambda_1^{O}>(<)0$ if and only if $a_1a_2e_1e_2<(>) b_1b_2$.

$(ii)$\ If $0< \delta<1$,
then the sufficient condition for
$\lambda_1^{O}>0$ is  $$\min\{k_1, k_2\}\delta>c_1(1-\delta),$$
and the necessary condition for $\lambda_1^{O}>0$ is $$\max\{k_1, k_2\}\delta> c_1(1-\delta),$$
as well as the sufficient condition for
$\lambda_1^{O}<0$ is  $$\max\{k_1, k_2\}\delta<c_1(1-\delta),$$
and the necessary condition of $\lambda_1^{O}<0$ is $$\min\{k_1, k_2\}\delta< c_1(1-\delta),$$
where $$c_1=\frac{-(b_1+b_2)+\sqrt{(b_1-b_2)^2+4a_1a_2e_1e_2}}{2}$$ defined in \eqref{b26}.

Specially, if $k_1=k_2$, then $\lambda_1^{O}>(<)0$ if and only if $$k_1\delta>(<) c_1(1-\delta).$$

$(iii)$\ If $\delta=1$, then $\overline{\lambda}_1^O\geq \underline{\lambda}_1^O>0$.

\end{cor}
\bpf
$(i)$ The result for  $\lambda_1^{O}>0$ (or $\lambda_1^{O}<0$) can be verified by the similar manner as in Theorem \ref{thm2.3} (i), we omit it here.

$(ii)$\ It can be easily seen from the discussion of Cases 1-3 in Theorem \ref{thm2.2} (i) that the following results for three cases are valid:

\qquad $(a)$\ if $k_1=k_2$, $\lambda_1^O=k_1\delta-c_1(1-\delta)$;

\qquad $(b)$\ if $k_1<k_2$, $k_1\delta-c_1(1-\delta)<\lambda_1^O<k_2\delta-c_1(1-\delta)$;

\qquad $(c)$\ if $k_1>k_2$, $k_2\delta-c_1(1-\delta)<\lambda_1^O<k_1\delta-c_1(1-\delta)$.

As a result, $\lambda_1^O>0$ provided that
\bess
\min\{k_1, k_2\}{\delta}> c_1({1-\delta}),
\eess
and if $\lambda_1^O>0$, we have
\bess
\max\{k_1, k_2\}{\delta}>c_1({1-\delta}).
\eess
$(iii)$\ The result for the case $\delta=1$ is directly from Theorem \ref{thm2.2} $(ii)$.
\epf

\vspace{3mm}

Now we study the monotonicity of $\lambda_1^O$ with respect to parameters of problem \eqref{b13}. Especially, we focus on its monotonicity  with respect to $\delta$, which provides a way to study the impact of the length $\delta \omega$ of cold season on prevention and control of virus.

To emphasize the dependence of the principal eigenvalue $\lambda_1^{O}$ on the parameters of  problem  \eqref{b13}, we now denote it by $\lambda_1^{O}(\Gamma, b_1, b_2, \delta)$, where $\Gamma$ is
collection of $a_i, e_i, \omega, k_1$ and $k_2$.

\begin{lem}
\label{lem2.5}
The following statements  are valid:

$(i)$\ $\lambda_1^O(b_1, b_2)$ is strictly increasing in $b_1$ and $b_2$ for $b_1, b_2\geq 0$;

$(ii)$\  $\lambda_1^O(\delta)$ is strictly increasing in $\delta$ for $\delta\in[0, 1)$ provided that $k_i\geq b_i (i=1,2)$.
\end{lem}
\bpf
$(i)$\ It can easily seen from problem \eqref{b14} that
$\lambda_1^O$ is strictly increasing in $b_1$ and $b_2$.

$(ii)$\ Let $0\leq\delta_1<\delta_2<1$,  and $(\lambda_1^O(\delta_i), \phi_i(t,x), \psi_i(t,x))$  is the principal eigen-pair of eigenvalue problem \eqref{b13} with $\delta=\delta_i(i=1, 2)$.

Let $\phi_{2}^*(t)=\phi_2((1-\delta_2)\omega-t)$ and $\psi_{2}^*(t)=\psi_2((1-\delta_2)\omega-t)$, it follows from \eqref{b13} with $\delta=\delta_2$ that
{\small
\be\left\{
\begin{array}{ll}
-\phi_{2t}^* = a_1e_1\psi_2^* -b_1 \phi_2^*+\lambda_{1}^{O}(\delta_2)\phi_2^*, &0<t\leq (1-\delta_2)\omega, \\[2mm]
-\psi_{2t}^*= a_2e_2\phi_2^* -b_2 \psi_2^*+\lambda_{1}^{O}(\delta_2)\psi_2^*,  &0<t\leq (1-\delta_2)\omega, \\[2mm]
-\phi_{2t}^*= -k_1 \phi_2^*+\lambda_{1}^{O}(\delta_2)\phi_2^*, &(1-\delta_2)\omega<t\leq \omega,   \\[2mm]
-\psi_{2t}^*= -k_2\psi_2^*+\lambda_{1}^{O}(\delta_2)\psi_2^*,  &(1-\delta_2)\omega<t\leq \omega,  \\[2mm]
\phi_2^*(0)=\phi_2^*(\omega), \,\,\psi_2^*(0)=\psi_2^*(\omega).
\end{array} \right.
\label{e01}\ee
}
Multiplying both sides of the first and third equation of \eqref{e01}  by $\phi_1$,
 then integrating over $(0, (1-\delta_2)\omega]$ and $((1-\delta_2)\omega, \omega]$, respectively, and
 recalling that $(\lambda_1^O(\delta_1), \phi_1(t,x), \psi_1(t,x))$  is the principal eigen-pair of eigenvalue problem \eqref{b13} with $\delta=\delta_1$, we obtain
{\small
\be
\begin{array}{lll}
&&\int_{0}^{(1-\delta_2)\omega}[a_1e_1 \psi_2^* \phi_1
-b_1 \phi_2^* \phi_1 +\lambda_{1}^{O}(\delta_2) \phi_2^* \phi_1]dt\\[2mm]
&=&- \phi_{2}^* \phi_1|_0^{(1-\delta_2)\omega}+\int_{0}^{(1-\delta_2)\omega}\phi_{2}^*
\phi_{1t} dt \\[2mm]
&=&- \phi_{2}^* \phi_1|_0^{(1-\delta_2)\omega}+\int_{0}^{(1-\delta_2)\omega}\phi_{2}^*[a_1e_1\psi_1-b_1\phi_1
+\lambda_{1}^{O}(\delta_1) \phi_1]dt
\end{array}
\label{e02}
\ee}
and
{\small
\be
\begin{array}{lll}
&&\int_{(1-\delta_2)\omega}^{\omega}[-k_1 \phi_2^* \phi_1 +\lambda_{1}^{O}(\delta_2) \phi_2^* \phi_1]dt\\[2mm]
&=&- \phi_{2}^* \phi_1|_{(1-\delta_2)\omega}^{\omega}+\int_{(1-\delta_2)\omega}^{\omega}\phi_{2}^* \phi_{1t} dt \\[2mm]
&=&- \phi_{2}^* \phi_1|_{(1-\delta_2)\omega}^{\omega}+\int_{(1-\delta_2)\omega}^{(1-\delta_1)\omega}
\phi_{2}^*[a_1e_1\psi_1-b_1\phi_1
+\lambda_{1}^{O}(\delta_1) \phi_1]dt\\[2mm]
&&+\int_{(1-\delta_1)\omega}^{\omega}\phi_{2}^*[-k_1\phi_1+\lambda_{1}^{O}(\delta_1) \phi_1]dt\\[2mm]
&=&-\phi_{2}^* \phi_1|_{(1-\delta_2)\omega}^{\omega}+\int_{(1-\delta_2)\omega}^{(1-\delta_1)\omega}
a_1e_1\phi_{2}^*\psi_1dt-\int_{(1-\delta_2)\omega}^{(1-\delta_1)\omega}b_1\phi_1\phi_{2}^*dt\\[2mm]
&&-\int_{(1-\delta_1)\omega}^{\omega}k_1\phi_1\phi_{2}^*dt
+\int_{(1-\delta_2)\omega}^{\omega}\lambda_{1}^{O}(\delta_1) \phi_1\phi_{2}^*dt.
\end{array}
\label{e03}
\ee}
Adding \eqref{e02} and \eqref{e03} yields
{\small
\be
\begin{array}{lll}
&&(\lambda_{1}^{O}(\delta_2)-\lambda_{1}^{O}(\delta_1))\int_{0}^{\omega}\phi_1\phi_2^*dt\\[2mm]
&=&a_1e_1[\int_0^{(1-\delta_1)\omega}\psi_1\phi_2^*dt-\int_0^{(1-\delta_2)\omega}\phi_1\psi_2^*dt]
-\int_{(1-\delta_2)\omega}^{(1-\delta_1)\omega}(k_1-b_1)\phi_1\phi_2^*dt.
\end{array}\label{e04}
\ee}
Similarly, multiplying both sides of the second and forth equation of \eqref{e01} by $\psi_1$,
  integrating over $(0, (1-\delta_2)\omega]$ and $((1-\delta_2)\omega, \omega]$, respectively,
and then adding them give
 {\small
\be
\begin{array}{lll}
&&(\lambda_{1}^{O}(\delta_2)-\lambda_{1}^{O}(\delta_1))\int_{0}^{\omega}\psi_1\psi_2^*dt\\[2mm]
&=&a_2e_2[\int_0^{(1-\delta_1)\omega}\phi_1\psi_2^*dt-\int_0^{(1-\delta_2)\omega}\psi_1\phi_2^*dt]
-\int_{(1-\delta_2)\omega}^{(1-\delta_1)\omega}(k_2-b_2)\psi_1\psi_2^*dt.
\end{array}\label{e05}
\ee}
It follows from \eqref{e04} and \eqref{e05} that
{\small
\bess
\begin{array}{lll}
&&(\lambda_{1}^{O}(\delta_2)-\lambda_{1}^{O}(\delta_1))[\frac{1}{a_1e_1}\int_{0}^{\omega}\phi_1\phi_2^*dt
+\frac{1}{a_2e_2}\int_{0}^{\omega}\psi_1\psi_2^*dt]\\[2mm]
&\geq&\int_0^{(1-\delta_1)\omega}\psi_1\phi_2^*dt-\int_0^{(1-\delta_2)\omega}\phi_1\psi_2^*dt
+[\int_0^{(1-\delta_1)\omega}\phi_1\psi_2^*dt-\int_0^{(1-\delta_2)\omega}\psi_1\phi_2^*dt]\\[2mm]
&>& \int_0^{(1-\delta_2)\omega}[\psi_1\phi_2^*-\phi_1\psi_2^*+\phi_1\psi_2^*-\psi_1\phi_2^*]dt=0
\end{array}
\eess}
provided that $k_i\geq b_i (i=1,2)$, which gives $\lambda_{1}^{O}(\delta_2)> \lambda_{1}^{O}(\delta_1)$.
\epf
\vspace{3mm}

\begin{rmk}
\label{rmk1.4}
The monotonicity of ${\lambda}_1^O$ in Lemma \ref{lem2.5} (ii) holds for $\delta\in [0, 1)$.
If $\delta=1$ and $b_1\neq k$, we have $\overline{\lambda}_1^O(\delta)>\underline{\lambda}_1^O(\delta)$, the principal eigenvalue does not exist. However, we can similarly obtain that $\overline{\lambda}_1^O(\delta)$ and $\underline{\lambda}_1^O(\delta)$ are strictly increasing in $\delta\in [0, 1]$.
\end{rmk}

\begin{figure}[ht]
\centering
\subfigure[$k_1\neq k_2$]{ {
\includegraphics[width=0.33\textwidth]{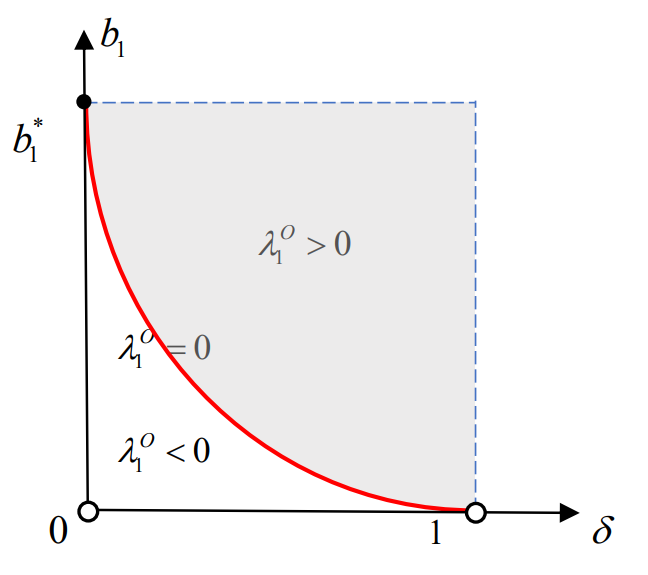}
} }
\subfigure[$k_1=k_2$]{ {
\includegraphics[width=0.33\textwidth]{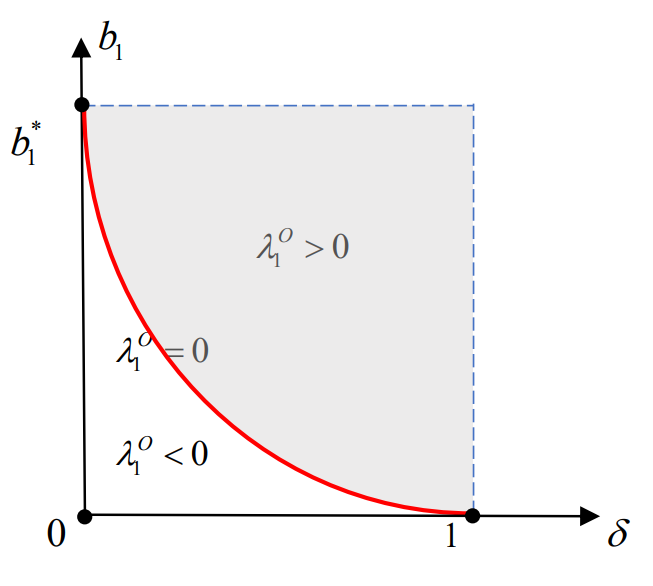}
} }
\caption{\scriptsize Figures (a) and (b) show that the contour lines (red colour) of $\lambda_1^O=0$ when $k_1\neq k_2$ and when $k_1=k_2$, respectively. }
\label{tu3}
\end{figure}

It follows from the above discussions on the monotonicity of $\lambda_1^O$ with respect to $b_1$ and $\delta$, we can get the contour lines of $\lambda_1^O=0$, see Fig. \ref{tu3}, where $b^*_1=a_1a_2e_1e_2/b_2$.

\begin{rmk}
As in our model, let $k_1=b_1$ and $k_2=k$, the above conclusions still hold.
\end{rmk}

Come back to problem \eqref{b01} involving nonlocal diffusion, we define the corresponding generalized
principal eigenvalues $\overline{\lambda}_{1}^{P}$ and $\underline{\lambda}_{1}^{P}$ of \eqref{b01} as
\bess
\begin{array}{ll}
\overline{\lambda}_{1}^{P}:=\inf\{\lambda\in \mathbb{R}\,|\, \exists \, \phi, \psi\in  C^{1, 0}(\mathbb{R}\times [-L_1, L_2]), \phi, \psi>0 {\rm\  and\ }  \phi, \psi {\ \rm are \ } \omega-{\rm periodic \ so \ as}\\
\qquad \qquad\ \ {\rm \ the\ inequalities\ in\  }\eqref{b27} {\rm \ hold}\},
\end{array}
\eess
where
{\small
\be\left\{
\begin{array}{ll}
\phi_t -d_1\mathcal{L}_1[\phi]\leq a_1e_1\psi -b_1 \phi+\lambda\phi, &0<t\leq (1-\delta)\omega, -L_1\leq x\leq L_2, \\[2mm]
\psi_t-d_2\mathcal{L}_2[\psi]\leq a_2e_2\phi -b_2 \psi+\lambda\psi,  &0<t\leq (1-\delta)\omega, -L_1\leq x\leq L_2,\\[2mm]
\phi_t-d_1\mathcal{L}_1[\phi]\leq -b_1 \phi+\lambda\phi, &(1-\delta)\omega<t\leq \omega, -L_1\leq x\leq L_2, \\[2mm]
\psi_t\leq -k\psi+\lambda\psi,  &(1-\delta)\omega<t\leq \omega, -L_1\leq x\leq L_2
\end{array} \right.
\label{b27}\ee}
and
\bess
\begin{array}{ll}
\underline{\lambda}_{1}^{P}:=\sup\{\lambda\in \mathbb{R}\,|\, \exists \, \phi, \psi\in  C^{1, 0}(\mathbb{R}\times [-L_1, L_2]), \phi, \psi>0 {\rm\  and\ }  \phi, \psi {\ \rm are \ } \omega-{\rm periodic\ } {\rm so \ as\  } \\
\qquad \qquad\ \ \eqref{b27}\textrm{\ with\ the\ inequalities\ reversed\ holds}\}.
\end{array}
\eess
If there exists a  positive function pair $(\phi, \psi)\in [C^{1, 0}(\mathbb{R}\times [-L_1, L_2])]^2$ such that problem \eqref{b01} with $\lambda=\lambda_1^P$ holds, then $\lambda_1^P$ is called a principal
eigenvalue of problem \eqref{b01}.

Based the above definitions, we have the following properties.

\begin{lem}
\label{lem2.7} $(i)$ $\underline{\lambda}_{1}^{P}([-L_1, L_2])=\underline{\lambda}_{1}^{P}([0, L_1+L_2])$ and $\overline{\lambda}_{1}^{P}([-L_1, L_2])=\overline{\lambda}_{1}^{P}([0, L_1+L_2])$.

 $(ii)$ $\underline{\lambda}_{1}^{P}([-L_1, L_2])$ and $\overline{\lambda}_{1}^{P}([-L_1, L_2])$ are decreasing with respect to $L_1+L_2$, in the sense, $\underline{\lambda}_{1}^{P}([-L_1, L_2])\leq \underline{\lambda}_{1}^{P}([-L_3, L_4])$ and  $\overline{\lambda}_{1}^{P}([-L_1, L_2])\leq \overline{\lambda}_{1}^{P}([-L_3, L_4])$ if $L_1+L_2>L_3+L_4$.
\end{lem}
\bpf
$(i)$ It follows from the transformation $(t, x)\to (t, z)$ where $z=x+L_1$ and $x\in [-L_1, L_2]$.

$(ii)$ Without loss of generality, we only need to prove that $\underline\lambda_{1}^{P}([-L_1, L_2])$ are decreasing with respect to $L_1+L_2$, that is, $\underline\lambda_{1}^{P}([-L_1, L_2])\leq \underline\lambda_{1}^{P}([-L_3, L_4])$ if $L_1+L_2>L_3+L_4$.
It suffices to prove $\underline\lambda_{1}^{P}([0, L_1+L_2])\leq\underline\lambda_{1}^{P}([0, L_3+L_4])$ if $L_1+L_2>L_3+L_4$ by (i).

Let $(\underline\lambda_{1}^{P}, \phi(t,x), \psi(t,x))$ is the generalized principal eigen-pair of
\eqref{b01} on $[0,\omega)\times[0, L_1+L_2]$, and it satisfies
{\small
\be\left\{
\begin{array}{ll}
\phi_t -d_1\mathcal{L}_1[\phi]\geq a_1e_1\psi -b_1 \phi+\underline\lambda_{1}^{P}\phi, &0<t\leq (1-\delta)\omega,\ 0\leq x\leq L_1+L_2, \\[2mm]
\psi_t-d_2\mathcal{L}_2[\psi]\geq a_2e_2\phi -b_2 \psi+\underline\lambda_{1}^{P}\psi,  &0<t\leq (1-\delta)\omega,\ 0\leq x\leq L_1+L_2,\\[2mm]
\phi_t-d_1\mathcal{L}_1[\phi]\geq -b_1 \phi+\underline\lambda_{1}^{P}\phi, &(1-\delta)\omega<t\leq \omega,\ 0\leq x\leq L_1+L_2, \\[2mm]
\psi_t\geq -k\psi+\underline\lambda_{1}^{P}\psi,  &(1-\delta)\omega<t\leq \omega,\ 0\leq x\leq L_1+L_2,\\[2mm]
\phi(0, x)=\phi(\omega, x), \psi(0, x)=\psi(\omega, x), &0\leq x\leq L_1+L_2.
\end{array} \right.
\label{b29}\ee}
Because of
\bess
\begin{array}{lll}
\,\,\,\,\,\, -d_1[\int_{0}^{L_1+L_2}J_1(x-y)\phi(t, y)dy-\phi(t, x)]\\[2mm]
=-d_1[\int_{0}^{L_3+L_4}J_1(x-y)\phi(t, y)dy-\phi(t, x)]
-d_1\int_{L_3+L_4}^{L_1+L_2}J_1(x-y)\phi(t, y)dy\\[2mm]
\leq, {\not\equiv} -d_1[\int_{0}^{L_3+L_4}J_1(x-y)\phi(t, y)dy-\phi(t, x)]
\end{array}
\eess
for $t\in(0, (1-\delta)\omega]$ and $x\in[0, L_3+L_4]$, we have
\bess
\phi_t-d_1[\int_{0}^{L_3+L_4}J_1(x-y)\phi(t, y)dy-\phi(t,x)]
\geq, {\not\equiv} a_1e_1\psi-b_1\phi+\underline\lambda_1^P([0, L_1+L_2])\phi
\eess
for $t\in(0, (1-\delta)\omega]$ and $x\in[0, L_3+L_4]$, and
\bess
\phi_t-d_1[\int_{0}^{L_3+L_4}J_1(x-y)\phi(t, y)dy-\phi(t,x)]
\geq, {\not\equiv} -b_1\phi+\underline\lambda_1^P([0, L_1+L_2])\phi
\eess
for $t\in((1-\delta)\omega, \omega]$ and $x\in[0, L_3+L_4]$.
Furthermore, we have
\bess
\psi_t-d_2[\int_{0}^{L_3+L_4}J_2(x-y)\psi(t, y)dy-\psi(t,x)]
\geq, {\not\equiv} a_2e_2\phi-b_2\psi+\underline\lambda_1^P([0, L_1+L_2])\psi
\eess
for $t\in(0, (1-\delta)\omega]$ and $x\in[0, L_3+L_4]$, and
\bess
\psi_t\geq -k\psi+\underline\lambda_1^P([0, L_1+L_2])\psi
\eess
for $t\in((1-\delta)\omega, \omega]$ and $x\in[0, L_3+L_4]$ in the similar manner.
Since $\phi(0,x)=\phi(\omega,x), \, \psi(0,x)=\phi(\omega,x)$ for $x\in[0, L_1+L_2]$, we then have
\bess
\phi(0,x)=\phi(\omega,x), \, \psi(0,x)=\phi(\omega,x),\,\, x\in[0, L_3+L_4].
\eess
It follows from the definition of $\underline\lambda_1^P([0, L_3+L_4])$ that
$$\underline\lambda_1^P([0, L_3+L_4])\geq\underline\lambda_1^P([0, L_1+L_2]).$$
\epf
\begin{thm}
\label{thm2.1}
Assume $\mathbf{(J)}$ and the kernel function $J_1(x)=J_2(x)$ for $x\in\mathbb{R}$ hold.

$(i)$ If the following condition holds
$$\mathbf{(H_2)}:\ (a)\ 0\leq \delta<1\ \textrm{or} \ (b)\ \delta=1, \, b_1-d_1\lambda_1^*=k.$$
Then $\overline{\lambda}_1^P=\underline{\lambda}_1^P=\lambda_1^P$, where $\lambda_1^{P}([-L_1, L_2])$ is the principal eigenvalue of eigenvalue problem \eqref{b01} with a positive eigenfunction pair $(\phi, \psi)  \in [C^{1, 0}([0, \omega]\times[-L_1, L_2])]^2$. Especially, if $b_1-d_1\lambda_1^*=k$, for any given $\delta\in[0,1]$, $\lambda_1^P$ can be explicitly written as $$\lambda_1^P=(b_1-d_1\lambda_1^*+s_1)\delta-s_1,$$
where $\lambda_1^*$ is defined in \eqref{ed} and
\be
s_{1}=\frac{-(b_1-d_1\lambda_1+b_2-d_2\lambda_1)+ \sqrt{[(b_1-d_1\lambda_1^*)-(b_2-d_2\lambda_1^*)]^2+4a_1a_2e_1e_2}}{2}.
\label{b28}
\ee

$(ii)$ If $\delta=1$ and $ b_1-d_1\lambda_1^*\neq k$, then
$$\overline{\lambda}_1^P=\max\{b_1-d_1\lambda_1^*, k\}> \min\{ b_1-d_1\lambda_1^*, k\}=\underline{\lambda}_1^P$$
 and the eigenvalue problem \eqref{b01} has unequal generalized principal eigenvalues $\overline{\lambda}_1^P$ and $\underline{\lambda}_1^P$.
\end{thm}
\bpf
We notice that the kernel functions $J_1$ and $J_2$ satisfy $J_1(x)=J_2(x)(:= J(x))$ in $\mathbb{R}$, then $\mathcal{L}_1=\mathcal{L}_2(:=\mathcal{L})$. Consider the following eigenvalue problem
\be
\mathcal{L}[u(x)]:=\int_{-L_1}^{L_2} J(x-y)u(y)dy-u(x)= \lambda u(x),\ -L_1\leq x\leq L_2.
\label{ed}
\ee
It follows from [\cite{Cao2019} Proposition 3.4] that problem \eqref{ed} admits a principal eigen-pair
$(\lambda^*_1, g(x))$ with $\lambda^*_1([-L_1, L_2])<0$, $g\in C[-L_1, L_2]$ and $g(x)>0$ in $[-L_1, L_2]$.

Suppose that $0\leq\delta<1$, we consider the spatial-independent eigenvalue problem
{\small
\bess\left\{
\begin{array}{lll}
f_{1t}=a_1e_1f_2-(b_1- d_1\lambda_{1}^*)f_1+\lambda_1^{O} f_1, &0<t\leq (1-\delta)\omega,\\[2mm]
f_{2t}=a_2e_2f_1-(b_2 -d_2\lambda_{1}^*)f_2+\lambda_1^{O} f_2, &0<t\leq (1-\delta)\omega,\\[2mm]
f_{1t}=-(b_1- d_1\lambda_1^*)f_1+\lambda_1^{O} f_1, &(1-\delta)\omega<t\leq \omega,\\[2mm]
f_{2t}=-k f_2+\lambda_1^{O} f_2, &(1-\delta)\omega<t\leq \omega,\\[2mm]
f_1(0)=f_1(\omega), \ f_2(0)=f_2(\omega).
\end{array}\right.
\eess
}
Theorem \ref{thm2.2}, where $b_1$ and $b_2$ are replaced by $b_1-d_1\lambda_1^*$ and $b_2-d_2\lambda_1^*$, implies that the above problem has a principal eigenvalue $\lambda_1^{O}(\Gamma, b_1-d_1\lambda_1^*, b_2-d_2\lambda_1^*)$ with the eigenfunction pair $(f_1(t), f_2(t))\in C^1[0,\omega]\times C^1[0,\omega]$.

\ Set
\bess
\phi(t,x)=f_1(t)g(x), \,\,\psi(t,x)=f_2(t)g(x),
\eess
for $t\in[0,\omega]$,\, $x\in[-L_1, L_2]$.
It is easy to check that
{\small
\be\left\{
\begin{array}{lll}
\phi_{t}-d_1\mathcal{L}[\phi]=a_1e_1\psi-b_1\phi+\lambda_1^{O} \phi, &0<t\leq (1-\delta)\omega,\\[2mm]
\psi_{t}-d_1\mathcal{L}[\psi]=a_2e_2\phi-b_2\psi+\lambda_1^{O} \psi, &0<t\leq (1-\delta)\omega,\\[2mm]
\phi_t-d_1\mathcal{L}[\phi]= -b_1 \phi+\lambda_1^{O}\phi, &(1-\delta)\omega<t\leq \omega, \\[2mm]
\psi_t= -k\psi+\lambda_1^{O}\psi,  &(1-\delta)\omega<t\leq \omega,\\[2mm]
\phi(0)=\phi(\omega), \ \psi(0)=\psi(\omega).
\end{array}
\label{b18}\right.
\ee
}
Therefore $\overline \lambda_1^{P}=\underline \lambda_1^{P}=\lambda_1^{P}=\lambda_1^{O}(\Gamma, b_1-d_1\lambda_1^*, b_2-d_2\lambda_1^*)$.
 The remaining proof of Theorem \ref{thm2.1} is similar to that of Theorem \ref{thm2.2}, we omit it here.

Specially, if  $\delta=1$, we have to consider two cases.
When $b_1-d_1\lambda_1^*=k$, then $\overline{\lambda}_1^P=\underline{\lambda}_1^P=b_1-d_1\lambda_1^*$. However, when $b_1-d_1\lambda_1^*\neq k$, then
$\overline{\lambda}_1^P=\max\{b_1-d_1\lambda_1^*, k\}$ and $\underline{\lambda}_1^P=\min\{b_1-d_1\lambda_1^*, k\}$, which means that $\overline{\lambda}_1^P\neq\underline{\lambda}_1^P.$\epf

\begin{cor}
\label{rmk1.3}
Similar to the discussions of Theorem \ref{thm2.2}, it follows from Theorem \ref{thm2.1} that
when $\delta\in[0, 1)$, the following hold true:

\qquad $(i)$\ if $k=b_1-d_1\lambda_1^*$,
$\lambda_1^{P}=(b_1-d_1\lambda_1^*+s_1)\delta-s_1;$

\qquad $(ii)$\ if $k>b_1-d_1\lambda_1^*$, $(b_1-d_1\lambda_1^*)\delta-s_1(1-\delta)<\lambda_1^P<k\delta-s_1(1-\delta)$;

\qquad $(iii)$\ if $k<b_1-d_1\lambda_1^*$, $k\delta-s_1(1-\delta)<\lambda_1^P<(b_1-d_1\lambda_1^*)\delta-s_1(1-\delta),$
where $\lambda_1^*$ and $s_{1}$ are  defined in \eqref{ed} and \eqref{b28}, respectively.
\end{cor}
From Theorem \ref{thm2.1} and Corollary \ref{cor1.3}, we have the following result.
\begin{cor}
\label{rmk1.1}
Assume $\mathbf{(J)}$ and the kernel function $J_1(x)=J_2(x)$ for $x\in\mathbb{R}$ hold.
The following statements  are valid:

$(i)$\ If $\delta=0$,  then

\ \quad$(a)$\ $\lambda_1^{P}= 0$  if and only if
$a_1a_2e_1e_2= (b_1-d_1\lambda_1^*)(b_2-d_2\lambda_1^*);$

\ \quad$(b)$\ $\lambda_1^{P}>(<) 0$  if and only if
$a_1a_2e_1e_2< (>)(b_1-d_1\lambda_1^*)(b_2-d_2\lambda_1^*);$

$(ii)$\ If $0<\delta<1$,
then

\ \quad$(a)$\ $\lambda_1^{P}= 0$  if and only if
\bess
\frac{a_1e_1(e^{s_1(1-\delta)\omega}-e^{(b_1-d_1\lambda_1^*)\delta\omega})}
{(b_1-d_1\lambda_1^*+s_1)(e^{s_2(1-\delta)\omega}-e^{(b_1-d_1\lambda_1^*)\delta\omega})}=
\frac{(b_2-d_2\lambda_1^*+s_2)(e^{s_1(1-\delta)\omega}-e^{k\delta\omega})}
{a_2e_2(e^{s_2(1-\delta)\omega}-e^{k\delta\omega})},
\eess
\ where $s_1$ is defined by \eqref{b28} and
\be
s_{2}=\frac{-(b_1-d_1\lambda_1^*+b_2-d_2\lambda_1^*)- \sqrt{[(b_1-d_1\lambda_1^*)-(b_2-d_2\lambda_1^*)]^2+4a_1a_2e_1e_2}}{2};
\label{e15}
\ee

\ \quad$(b)$\ the sufficient condition for $\lambda_1^{P}> 0$  is
\bess
\min\{b_1-d_1\lambda_1^*, \,k\}\delta> s_1(1-\delta),
\eess
\ and the necessary condition for $\lambda_1^{P}> 0$  is
\bess
\max\{b_1-d_1\lambda_1^*, \,k\}\delta> s_1(1-\delta),
\eess
\ specially, if $b_1-d_1\lambda_1^*=k$, then $\lambda_1^{P}> 0$  if and only if
$
(b_1-d_1\lambda_1^*)\delta> s_1(1-\delta);
$

\ \quad$(c)$\ the sufficient condition for $\lambda_1^{P}<0$  is
\bess
\max\{b_1-d_1\lambda_1^*,\,k\}\delta< s_1({1-\delta}),
\eess
\ and the necessary condition for $\lambda_1^{P}<0$  is
\bess
\min\{b_1-d_1\lambda_1^*, \,k\}\delta< s_1({1-\delta}),
\eess
\ specially, if $b_1-d_1\lambda_1^*=k$, then $\lambda_1^{P}< 0$  if and only if
$
(b_1-d_1\lambda_1^*)\delta< s_1(1-\delta).
$

$(iii)$\ If $\delta=1$, then $\overline\lambda_1^{P}\geq \underline\lambda_1^{P}>0$ always holds.
\end{cor}
\vspace{5mm}

 To stress the dependence of $\lambda_1^{P}$ on the parameters of problem \eqref{b01}, we denote it by $$\lambda_1^{P}(\Gamma, [-L_1, L_2], b_i, d_i, \delta),$$
where $\Gamma=\{a_i, e_i, k, \omega\}$. With the above definition, we have the monotonicity of $\lambda_1^P$ by using Lemma \ref{lem2.5} and Remark \ref{rmk1.4}.

\begin{cor}
\label{rmk1.2}
Assume $\mathbf{(J)}$ and the kernel function $J_1(x)=J_2(x)$ for $x\in\mathbb{R}$ hold.
The following statements are valid:

$(i)$\ $\overline{\lambda}_1^P$ and $\underline{\lambda}_1^P$ are strictly increasing with respect to $b_1\in (0, +\infty)$.

$(ii)$\ $\lambda_1^P(\delta)$ is strictly increasing in $\delta\in [0, 1)$. Specially, if $\delta=1$, $\lambda_1^P(\delta)$ is strictly increasing in $\delta$ when $b_1-d_1\lambda_1^*=k$, and when $b_1-d_1\lambda_1^*\neq k$, then $\overline{\lambda}_1^P(\delta)$ and $\underline{\lambda}_1^P(\delta)$ are strictly increasing in $\delta$.
\end{cor}

In particular, if $\delta=0$, the periodic eigenvalue problem \eqref{b01} is without seasonal succession, the following result is known in \cite{Du2020}.
\begin{lem}

$($\cite{Du2020}  Proposition 2.1\,$)$
Assume $(J)$ holds. Then the eigenvalue problem \eqref{b01} with $\delta=0$ has a principal
eigenvalue $\lambda_1^P$ with a positive eigenfunction pair $(\phi, \psi)(x)  \in [C[-L_1, L_2])]^2$. Moreover, $\lambda_1^P$ is an algebraically simple eigenvalue.
\end{lem}

Next we present a comparison principle about the generalized principal eigenvalue, which will be used in the sequel.
\begin{lem}
\label{lem2.1}
$(i)$ Suppose $\mathbf{(J)}$ holds and let $\overline{\lambda}_{1}^{P}$ be the generalized principal eigenvalue of \eqref{b01}. If there exist two functions $\tilde{\phi}, \tilde{\psi}\in C^{1,0}([0, \omega]\times[-L_1, L_2])$ with $\tilde{\phi}(t, x), \tilde{\psi}(t,x)\geq, {\not\equiv}0$ in $[0, \omega]\times[-L_1, L_2]$ such that
{\small
\be\left\{
\begin{array}{ll}
\tilde{\phi}_t- d_1\mathcal{L}_1[\tilde{\phi}]\geq a_1e_1\tilde{\psi} -b_1 \tilde{\phi}+\lambda_1\tilde{\phi}, &0<t\leq (1-\delta)\omega, -L_1\leq x\leq L_2, \\[2mm]
\tilde{\psi}_t-d_2\mathcal{L}_2[\tilde{\psi}]\geq a_2e_2\tilde{\phi} -b_2 \tilde{\psi}+\lambda_1\tilde{\psi},  &0<t\leq (1-\delta)\omega, -L_1\leq x\leq L_2, \\[2mm]
\tilde{\phi}_t-d_1\mathcal{L}_1[\tilde{\phi}]\geq -b_1 \tilde{\phi}+\lambda_1\tilde{\phi}, &(1-\delta)\omega<t\leq \omega, -L_1\leq x\leq L_2, \\[2mm]
\tilde{\psi}_t\geq-k\tilde{\psi}+\lambda_1\tilde{\psi},  &(1-\delta)\omega<t\leq \omega, -L_1\leq x\leq L_2, \\[2mm]
\tilde{\phi}(0,x)\geq \tilde{\phi}(\omega,x), \,\,\tilde{\psi}(0,x)\geq\tilde{\psi}(\omega,x),  &-L_1\leq x\leq L_2
\end{array} \right.
\label{b02}\ee}
for a constant $\lambda_1$, then $\overline{\lambda}_{1}^{P}\geq \lambda_1$ and the equality holds only when $(\overline{\lambda}_{1}^{P}, \tilde{\phi}, \tilde{\psi})$ is the principal eigen-pair
of eigenvalue problem \eqref{b01}.

$(ii)$ Suppose $\mathbf{(J)}$ holds and let $\underline{\lambda}_{1}^{P}$ be the generalized principal eigenvalue of \eqref{b01}. If there exist two functions $\hat{\phi}, \hat{\psi}\in C^{1, 0}([0, \omega]\times[-L_1, L_2])$ with $\hat{\phi}(t, x), \hat{\psi}(t, x)\geq, {\not\equiv} 0$ in $[0, \omega]\times[-L_1, L_2]$ such that
{\small
\be\left\{
\begin{array}{ll}
\hat{\phi}_t-d_1\mathcal{L}_1[\hat{\phi}]\leq a_1e_1\hat{\psi} -b_1 \hat{\phi}+\lambda_2\hat{\phi}, &0<t\leq (1-\delta)\omega, -L_1\leq x\leq L_2, \\[2mm]
\hat{\psi}_t-d_2\mathcal{L}_2[\hat{\psi}]\leq a_2e_2\hat{\phi} -b_2 \hat{\psi}+\lambda_2\hat{\psi},  &0<t\leq (1-\delta)\omega, -L_1\leq x\leq L_2,\\[2mm]
\hat{\phi}_t-d_1\mathcal{L}_1[\hat{\phi}]\leq -b_1 \hat{\phi}+\lambda_2\hat{\phi}, &(1-\delta)\omega<t\leq \omega, -L_1\leq x\leq L_2, \\[2mm]
\hat{\psi}_t\leq-k\hat{\psi}+\lambda_2\hat{\psi},  &(1-\delta)\omega<t\leq \omega, -L_1\leq x\leq L_2, \\[2mm]
\hat{\phi}(0,x)\leq \hat{\phi}(\omega,x), \,\,\hat{\psi}(0,x)\leq\hat{\psi}(\omega,x),  &-L_1\leq x\leq L_2
\end{array} \right.
\label{b03}\ee}
for a constant $\lambda_2$, then $\underline{\lambda}_{1}^{P}\leq \lambda_2$ and the equality holds only when $(\underline{\lambda}_{1}^{P}, \hat{\phi}, \hat{\psi})$ is the principal eigen-pair
of eigenvalue problem \eqref{b01}.
\end{lem}
\bpf
We only verify the assertion (i), because (ii) can be proved in the same manner.

Define
\be
\langle u,\, v\rangle:=\int_{-L_1}^{L_2}u(t,x)v(t,x)dx.
\label{f05}
\ee
Let $(\phi(t,x), \psi(t,x))$ be positive eigenfunctions corresponding to the generalized principal eigenvalue $\overline{\lambda}_{1}^{P}$ in \eqref{b01}.
Set $(\phi_1(t,x)=\phi((1-\delta)\omega-t, x),\, \psi_1(t,x)=\psi((1-\delta)\omega-t, x)$, it is easy to see that
{\small
\be\left\{
\begin{array}{ll}
-\phi_{1t} -d_1\mathcal{L}_1[\phi_1]\leq a_1e_1\psi_1 -b_1 \phi_1+\overline{\lambda}_{1}^{P}\phi_1, &0<t\leq (1-\delta)\omega, -L_1\leq x\leq L_2, \\[2mm]
-\psi_{1t}-d_2\mathcal{L}_2[\psi_1]\leq a_2e_2\phi_1 -b_2 \psi_1+\overline{\lambda}_{1}^{P}\psi_1,  &0<t\leq (1-\delta)\omega, -L_1\leq x\leq L_2,\\[2mm]
-\phi_{1t}-d_1\mathcal{L}_1[\phi_1]\leq -b_1 \phi_1+\overline{\lambda}_{1}^{P}\phi_1, &(1-\delta)\omega<t\leq \omega, -L_1\leq x\leq L_2, \\[2mm]
-\psi_{1t}\leq -k\psi_1+\overline{\lambda}_{1}^{P}\psi_1,  &(1-\delta)\omega<t\leq \omega,
-L_1\leq x\leq L_2,\\[2mm]
\phi_1(0,x)=\phi_1(\omega,x), \,\,\psi_1(0,x)=\psi_1(\omega,x),  &-L_1\leq x\leq L_2.
\end{array} \right.
\label{b08}\ee
}
We first restrict $t\in(0,(1-\delta)\omega]$,
multiplying both sides of the first  equation of \eqref{b08}  by $\tilde{\phi}$,
and the second  equation of that by $\tilde{\psi}$, then integrating over $[-L_1, L_2]$, yields
{\small
\bess
\begin{array}{lll}
-\langle\phi_{1t}, \tilde{\phi}\rangle-d_1\langle\mathcal{L}_1[\phi_1], \tilde{\phi}\rangle\leq a_1e_1\langle\psi_1, \tilde{\phi}\rangle
-b_1\langle\phi_1, \tilde{\phi}\rangle+\overline{\lambda}_{1}^{P}\langle\phi_1, \tilde{\phi}\rangle,\\[2mm]
-\langle\psi_{1t}, \tilde{\psi}\rangle-d_2\langle\mathcal{L}_2[\psi_1], \tilde{\psi}\rangle\leq a_2e_2\langle\phi_1, \tilde{\psi}\rangle
-b_2\langle\psi_1, \tilde{\psi}\rangle+\overline{\lambda}_{1}^{P}\langle\psi_1, \tilde{\psi}\rangle.
\end{array}
\eess}
Since $\langle\mathcal{L}_1[\phi_1], \tilde{\phi}\rangle=\langle\phi_1, \mathcal{L}_1[\tilde{\phi}]\rangle$ and $\langle\mathcal{L}_2[\psi_1], \tilde{\psi}\rangle=\langle\psi_1, \mathcal{L}_2[\tilde{\psi}]\rangle$, it follows from \eqref{b02} that
{\small
\bess
\langle\phi_1, \tilde{\phi}_t-a_1e_1\tilde{\psi}-\lambda_1\tilde{\phi}\rangle\geq\langle\phi_1,
d_1\mathcal{L}_1[\tilde{\phi}]-b_1\tilde{\phi}\rangle\geq -\langle\phi_{1t}, \tilde{\phi}\rangle-a_1e_1\langle\psi_1, \tilde{\phi}\rangle
-\overline{\lambda}_{1}^{P}\langle\phi_1, \tilde{\phi}\rangle,\\[2mm]
\langle\psi_1, \tilde{\psi}_t-a_2e_2\tilde{\phi}-\lambda_1\tilde{\psi}\rangle\geq\langle\psi_1,
d_2\mathcal{L}_2[\tilde{\psi}]-b_2\tilde{\psi}\rangle\geq -\langle\psi_{1t}, \tilde{\psi}\rangle-a_2e_2\langle\phi_1, \tilde{\psi}\rangle
-\overline{\lambda}_{1}^{P}\langle\psi_1, \tilde{\psi}\rangle,
\eess}
equivalently,
{\small
\be
(\lambda_1-\overline{\lambda}_{1}^{P})\langle\phi_1, \tilde{\phi}\rangle\leq a_1e_1(\langle\psi_1, \tilde{\phi}\rangle-\langle\phi_1, \tilde{\psi}\rangle)+\langle\phi_{1t}, \tilde{\phi}\rangle+\langle\phi_1, \tilde{\phi}_t\rangle,
\label{b09}
\ee
\be
(\lambda_1-\overline{\lambda}_{1}^{P})\langle\psi_1, \tilde{\psi}\rangle\leq a_2e_2(\langle\phi_1, \tilde{\psi}\rangle-\langle\psi_1, \tilde{\phi}\rangle)+\langle\psi_{1t}, \tilde{\psi}\rangle+\langle\psi_1, \tilde{\psi}_t\rangle
\label{b22}
\ee}
for $t\in(0, (1-\delta)\omega]$.
For the third equation of \eqref{b08}, repeating the above process,  we can easily obtain
{\small\be
(\lambda_1-\overline{\lambda}_{1}^{P})\langle\phi_1, \tilde{\phi}\rangle\leq \langle\phi_{1t}, \tilde{\phi}\rangle+\langle\phi_1, \tilde{\phi}_t\rangle\,\,
\label{b10}
\ee
for $t\in((1-\delta)\omega, \omega]$.
For \eqref{b09} and  \eqref{b10}, integrating over $(0, (1-\delta)\omega)$ and $((1-\delta)\omega, \omega)$ respectively, and adding them together, we obtain
{\small
\be
\begin{array}{lll}
&(\lambda_1-\overline{\lambda}_{1}^{P})\int_{0}^{\omega}\langle\phi_1,\, \tilde{\phi}\rangle dt\\[2mm]
\leq& \int_{0}^{(1-\delta)\omega}a_1e_1(\langle\psi_1, \tilde{\phi}\rangle-\langle\phi_1, \tilde{\psi}\rangle)dt
+\int_{0}^{\omega}\langle\phi_1, \tilde{\phi}_t\rangle+\langle\phi_{1t}, \tilde{\phi}\rangle dt\\[2mm]
=&\int_{0}^{(1-\delta)\omega}a_1e_1(\langle\psi_1, \tilde{\phi}\rangle-\langle\phi_1, \tilde{\psi}\rangle)dt+\int_{-L_1}^{L_2}
[\phi_1(\omega,x)\tilde{\phi}(\omega,x)-\phi_1(0,x)\tilde{\phi}(0,x)]dx\\[2mm]
\leq& \int_{0}^{(1-\delta)\omega}a_1e_1(\langle\psi_1, \tilde{\phi}\rangle-\langle\phi_1, \tilde{\psi}\rangle)dt,\label{b24}
\end{array}
\ee}
since $\phi_1(\omega,x)=\phi_1(0,x)$ in \eqref{b08} and $\tilde{\phi}(\omega,x)\leq\tilde{\phi}(0,x)$ by \eqref{b02}.
Similarly, multiplying both sides of the forth equation of \eqref{b08} by $\tilde{\psi}$ and integrating over $[-L_1, L_2]$, yields
\be
-\langle\psi_{1t}, \tilde{\psi}\rangle\leq -k\langle\psi_1, \tilde{\psi}\rangle+\overline{\lambda}_{1}^{P}\langle\psi_1, \tilde{\psi}\rangle.
\label{b23}
\ee}
For \eqref{b22} and  \eqref{b23}, integrating over $(0, (1-\delta)\omega)$ and $((1-\delta)\omega, \omega)$ respectively, we obtain
{\small
\bess
\begin{array}{llll}
&&(\lambda_1-\overline{\lambda}_{1}^{P})\int_{0}^{(1-\delta)\omega}\langle\psi_1, \tilde{\psi}\rangle dt\\[2mm]
&\leq&\int_{0}^{(1-\delta)\omega}a_2e_2(\langle\phi_1, \tilde{\psi}\rangle-\langle\psi_1, \tilde{\phi}\rangle)dt
+\int_{0}^{(1-\delta)\omega}\langle\psi_1, \tilde{\psi}_t\rangle+\langle\psi_{1t}, \tilde{\psi}\rangle dt\\[2mm]
&\leq&\int_{0}^{(1-\delta)\omega}a_2e_2(\langle\phi_1, \tilde{\psi}\rangle-\langle\psi_1, \tilde{\phi}\rangle)dt\\[2mm]
&&+\int_{-L_1}^{L_2}[\tilde{\psi}((1-\delta)\omega, x)\psi_1((1-\delta)\omega, x)-\tilde{\psi}(0, x)\psi_1(0, x)]dx,
\end{array}
\eess}
and
{\small
\bess
(\lambda_1-\overline{\lambda}_{1}^{P})\int_{(1-\delta)\omega}^{\omega}\langle\psi_1, \tilde{\psi}\rangle dt
\leq\int_{-L_1}^{L_2}[\tilde{\psi}(\omega, x)\psi_1(\omega, x)-\tilde{\psi}((1-\delta)\omega, x)\psi_1((1-\delta)\omega, x)]dx,
\eess
}
and then adding them together, yields
{\small
\be
(\lambda_1-\overline{\lambda}_{1}^{P})\int_{0}^{\omega}\langle\psi_1, \tilde{\psi}\rangle dt
\leq \int_{0}^{(1-\delta)\omega}a_2e_2(\langle\phi_1, \tilde{\psi}\rangle-\langle\psi_1, \tilde{\phi}\rangle)dt.
\label{b25}
\ee}
It follows from \eqref{b24} and \eqref{b25} that
{\small
\bess
\begin{array}{lll}
&&(\lambda_1-\overline{\lambda}_{1}^{P})\big[\frac{1}{a_1e_1}\int_{0}^{\omega}\langle\phi_1, \tilde{\phi}\rangle dt+\frac{1}{a_2e_2}\int_{0}^{\omega}\langle\psi_1, \tilde{\psi}\rangle dt\big]\\[2mm]
&\leq&\int_{0}^{(1-\delta)\omega}(\langle\psi_1, \tilde{\phi}\rangle-\langle\phi_1, \tilde{\psi}\rangle+\langle\phi_1, \tilde{\psi}\rangle-\langle\psi_1, \tilde{\phi}\rangle)dt = 0.
\end{array}
\eess}
Because of $\langle\phi_1,\, \tilde{\phi}\rangle>0$ and $\langle\psi_1,\, \tilde{\psi}\rangle>0$, we admit $\overline{\lambda}_{1}^{P}\geq\lambda_1$.
\epf

\begin{lem}
\label{lem2.2}
Assume $\mathbf{(J)}$, $\mathbf{(H_2)}$ and $J_1(x)=J_2(x)$ in $x\in \mathbb{R}$ hold, and let $\lambda_1^P([-L_1, L_2])$  be the principal eigenvalue of \eqref{b01}. Then
 $\lambda_1^P([-L_1, L_2])$ is strictly decreasing and continuous in $L:=L_1+L_2$.
\end{lem}
\bpf
Notice that $\lambda_1^P([-L_1, L_2])$ is well-defined under the assumptions.
 It follows from Lemma \ref{lem2.7} (ii) that $\lambda_1^P([-L_1, L_2])$ is decreasing in $L_1+L_2$.
 Moreover, by the definition of the principal eigenvalue, we can see that $\lambda_1^P([-L_1, L_2])$ is strict decreasing from the proof of Lemma \ref{lem2.7} (ii).

Next, we prove the continuity of  $\lambda_1^P([-L_1, L_2])$. We only need to prove that $\lambda_1^P([0, L])$ is continuous in $L$ by Lemma \ref{lem2.7} (i), that is, for any $\epsilon>0$,
there is $\gamma^*>0$ such that as long as $|L_1^*-L|<\gamma^*$,
\bess
|\lambda_1^P([0, L_1^*])-\lambda_1^P([0, L])|<\epsilon.
\eess

The following proof is divided into two cases, we mainly follow the approach of \cite{Du2020}, but some modifications are required for our model, which characterizes the behaviors in warm and cold seasons.

\textbf{Case 1.}  If $L_1^*\leq L$. Clearly, $\lambda_1^P([0, L_1^*])>\lambda_1^P([0, L])>\lambda_1^P([0, L])-\epsilon$ holds by the conclusion $(i)$.  We next conclude that for $L_1^*\in (L-\gamma^*, L]$ such that $\lambda_1^P([0, L_1^*])<\lambda_1^P([0, L])+\epsilon$.
In fact, let $(\phi(t,x), \psi(t,x))$ be the positive eigenfunction pair corresponding to $\lambda_1^P([0, L])$, and
\bess
m:=\min\Big\{\min\limits_{(t,x)\in[0,\omega]\times[-L,L]}\phi(t,x),\ \min\limits_{(t,x)\in[0,\omega]\times[-L,L]}\psi(t,x)\Big\},\\
M:=\max\Big\{\max\limits_{(t,x)\in[0,\omega]\times[-L,L]}\phi(t,x),\ \max\limits_{(t,x)\in[0,\omega]\times[-L,L]}\psi(t,x)\Big\}.
\eess
For $(t, x)\in(0, (1-\delta)\omega]\times[0, L_1^*]$
\be
\begin{array}{lll}
\,\,\,\,\, -d_1[\int_{0}^{L}J_1(x-y)\phi(t,y)dy-\phi(t,x)]\\[2mm]
=-d_1[\int_{0}^{L_1^*}J_1(x-y)\phi(t,y)dy-\phi(t,x)]-d_1\int_{L_1^*}^{L}J_1(x-y)\phi(t,y)dy\\[2mm]
\geq-d_1[\int_{0}^{L_1^*}J_1(x-y)\phi(t,y)dy-\phi(t,x)]-d_1M\|J_1\|_{\infty}(L-L_1^*)\\[2mm]
\geq-d_1[\int_{0}^{L_1^*}J_1(x-y)\phi(t,y)dy-\phi(t,x)]-\frac{d_1M\|J_1\|_{\infty}(L-L_1^*)}{m}\phi.
\end{array}
\label{b04}
\ee
It follows that
\bess
\begin{array}{lll}
&&\phi_t-d_1[\int_{0}^{L_1^*}J_1(x-y)\phi(t,y)dy-\phi(t,x)]\\
&\leq& a_1e_1\psi(t,x)-b_1\phi(t,x)
+[\lambda_1^P([0, L])+\frac{d_1M\|J_1\|_{\infty}(L-L_1^*)}{m}]\phi(t,x)
\end{array}
\eess
for $(t, x)\in(0, (1-\delta)\omega]\times[0, L_1^*]$ and
\bess
\begin{array}{lll}
&&\phi_t-d_1[\int_{0}^{L_1^*}J_1(x-y)\phi(t,y)dy-\phi(t,x)]\\[2mm]
 &\leq& -b_1\phi
+[\lambda_1^P([0, L])+\frac{d_1M\|J_1\|_{\infty}(L-L_1^*)}{m}]\phi
\end{array}
\eess
for $(t, x)\in((1-\delta)\omega, \omega]\times[0, L_1^*].$

Analogously, we have
{\small
\bess
\begin{array}{lll}
&&\psi_t-d_2[\int_{0}^{L_1^*}J_2(x-y)\psi(t,y)dy-\psi(t,x)]\\[2mm]
&\leq& a_2e_2\phi-b_2\psi
+[\lambda_1^P([0, L])+\frac{d_2M\|J_2\|_{\infty}(L-L_1^*)}{m}]\psi
\end{array}
\eess
}
for $(t, x)\in(0, (1-\delta)\omega]\times[0, L_1^*]$ and
{\small
\bess
\begin{array}{lll}
\psi_t\leq -k\psi
+[\lambda_1^P([0, L])+\frac{d_2M\|J_2\|_{\infty}(L-L_1^*)}{m}]\psi
\end{array}
\eess
}
for $(t, x)\in((1-\delta)\omega, \omega]\times[0, L_1^*].$
Denote
\bess
K:=\max\bigg\{\frac{d_1M\|J_1\|_{\infty}}{m}, \ \frac{d_2M\|J_2\|_{\infty}}{m}\bigg\},
\eess
we apply Lemma \ref{lem2.1} to obtain that
\bess
\lambda_1^P([0, L_1^*])\leq \lambda_1^P([0, L])+K(L-L_1^*)<\lambda_1^P([0, L])+\epsilon.
\eess
holds if $L-L_1^*<\gamma^*:=\frac{\epsilon}{K}$.

\textbf{Case 2.}  If $L_1^*\geq L$. By the conclusion $(i)$, $\lambda_1^P([0, L_1^*])<\lambda_1^P([0, L])<\lambda_1^P([0, L])+\epsilon$ holds,  we only need to show
that $\lambda_1^P([0, L_1^*])>\lambda_1^P([0, L])-\epsilon$ when $L_1^*\in [L, L+\gamma^*)$.
Let $(\phi(t,x), \psi(t,x))$ is the eigenfunction pair of \eqref{b01} on $[0,\omega)\times[0,L]$.
We extend $\phi$ and $\psi$ by $\phi(t,x)=\phi(t,L), \psi(t,x)=\psi(t,L)$ for $x\geq L$. In such a case, $m$ and $M$ are still the minimum and maximum of the extend eigenfunctions $\phi$ and $\psi$.

Firstly, in the warm season $t\in(0, (1-\delta)\omega]$. Similarly to \eqref{b04}, it is calculated that
\be\begin{array}{lll}
\,\,\,\,\, -d_1[\int_{0}^{L_1^*}J_1(x-y)\phi(t,y)dy-\phi(t,x)]\\[2mm]
\geq-d_1[\int_{0}^{L}J_1(x-y)\phi(t,y)dy-\phi(t,x)]-\frac{d_1M\|J_1\|_{\infty}(L_1^*-L)}{m}\phi(t,x)
\end{array}
\label{b05}
\ee
for $x\in [0, +\infty)$.
We denote, for $x\in[0, +\infty)$,
\bess
A(t,x)=\phi_t-d_1[\int_{0}^{L}J_1(x-y)\phi(t,y)dy-\phi(t,x)]-a_1e_1\psi(t,x)+b_1\phi(t,x).
\eess
Note that $A(t,x)$ is continuous and
\be
\label{b06}
A(t,x)=\lambda_1^P([0, L])\phi(t,x)\qquad {\rm for} \,\,x\in [0, L],
\ee
then there exists $\gamma_1>0$ such that
\be
|A(t,x)-\lambda_1^P([0, L])\phi(t,x)|\leq\frac{m\epsilon}{2} \qquad{\rm for}\,\,x\in[L,L+\gamma_1].
\label{b07}
\ee
Thanks to \eqref{b05}, \eqref{b06} and \eqref{b07}, we know, for $x\in [0, L]$,
\bess
\begin{array}{lll}
&\phi_t-d_1[\int_{0}^{L_1^*}J_1(x-y)\phi(t,y)dy-\phi(t,x)]-a_1e_1\psi(t,x)+b_1\phi(t,x)\\[2mm]
\geq& A(t,x)-\frac{d_1M\|J_1\|_{\infty}(L_1^*-L)}{m}\phi(t,x)\\[2mm]
=&\lambda_1^P([0, L])\phi(t,x)-\frac{d_1M\|J_1\|_{\infty}(L_1^*-L)}{m}\phi(t,x)\\[2mm]
\geq& [\lambda_1^P([0, L])-K(L_1^*-L)]\phi(t,x),
\end{array}
\eess
in addition, for  $x\in[L,L+\gamma_1]$,
\bess
\begin{array}{ll}
 &\phi_t-d_1[\int_{0}^{L_1^*}J_1(x-y)\phi(t,y)dy-\phi(t,x)]-a_1e_1\psi(t,x)+b_1\phi(t,x)\\[2mm]
\geq& A(t,x)-\frac{d_1M\|J_1\|_{\infty}(L_1^*-L)}{m}\phi(t,x)\\[2mm]
\geq&\lambda_1^P([0, L])\phi(t,x)-\frac{m\epsilon}{2}-\frac{d_1M\|J_1\|_{\infty}(L_1^*-L)}{m}\phi(t,x)\\[2mm]
\geq& [\lambda_1^P([0, L])-\frac{\epsilon}{2}-K(L_1^*-L)]\phi(t,x).
\end{array}
\eess
Moreover, when $t\in((1-\delta)\omega, \omega]$,
\bess
\begin{array}{lll}
&\phi_t-d_1[\int_{0}^{L_1^*}J_1(x-y)\phi(t,y)dy-\phi(t,x)]+b_1\phi(t,x)\\[2mm]
\geq &[\lambda_1^P([0, L])-K(L_1^*-L)]\phi(t,x), \quad\quad  x\in [0, L],\\[2mm]
&\phi_t-d_1[\int_{-L_1^*}^{L_1^*}J_1(x-y)\phi(t,y)dy-\phi(t,x)]+b_1\phi(t,x)\\[2mm]
\geq& [\lambda_1^P([0, L])-\frac{\epsilon}{2}-K(L_1^*-L)]\phi(t,x),\quad x\in[L,L+\gamma_1].\\[2mm]
\end{array}
\eess
Then it follows from analogous calculations that
\bess
\begin{array}{lll}
&\psi_t-d_2[\int_{0}^{L_1^*}J_2(x-y)\psi(t,y)dy-\psi(t,x)]-a_2e_2\phi(t,x)+b_2\psi(t,x)\\[2mm]
\geq &[\lambda_1^P([0, L])-K(L_1^*-L)]\psi(t,x), \qquad\qquad  x\in [0, L],\\[2mm]
&\psi_t-d_2[\int_{0}^{L_1^*}J_2(x-y)\psi(t,y)dy-\psi(t,x)]-a_2e_2\phi(t,x)+b_2\psi(t,x)\\[2mm]
\geq& [\lambda_1^P([0, L])-\frac{\epsilon}{2}-K(L_1^*-L)]\psi(t,x),\qquad x\in[L, L+\gamma_1],\\[2mm]
&\psi_t+k\psi(t,x)
\geq [\lambda_1^P([0, L])-K(L_1^*-L)]\psi(t,x), \quad  x\in [0, L],\\[2mm]
&\psi_t+k\psi(t,x)
\geq [\lambda_1^P([0, L])-\frac{\epsilon}{2}-K(L_1^*-L)]\psi(t,x),\quad x\in[L, L+\gamma_1].
\end{array}
\eess
If $L_1^*-L<\gamma^*:=\min\big\{\gamma_1, \frac{\epsilon}{2K}\big\}$, by Lemma \ref{lem2.1} we obtain
\bess
\lambda_1^P([0, L_1^*])>\lambda_1^P([0, L])-\epsilon.
\eess
Combining the aforementioned  two situations yields directly that $\lambda_1^P([0, L])$ is continuous.
\epf

\vspace{3mm}
Now, we consider the nonlocal free boundary problem \eqref{a07}, and denote
\bess
\begin{array}{rl}
\lambda_1^{F}(t):=&\lambda_1^P(\Gamma, [g(t), h(t)], b_i, d_i, \delta)\\[2mm]
=&\lambda_1^O(\Gamma, b_1-d_1\lambda_1^*([g(t), h(t)]), b_2-d_2\lambda_1^*([g(t), h(t)]), \delta),
\end{array}
\eess
where $\lambda_1^*$ is defined in \eqref{ed} with $[-L_1, L_2]$ replaced by $[g(t), h(t)]$.
\begin{lem}
\label{lem2.15}
Assume $\mathbf{(J)}$, $\mathbf{(H_2)}$ and $J_1(x)=J_2(x)$ in $x\in \mathbb{R}$ hold.
Then $\lambda_1^{F}(t)$ is nonincreasing for $t>0$.  Moreover, if $0\leq \delta <1$, then
$\lambda_1^{F}(t)>\lambda_1^{F}(t+\omega)$ for any $t\geq0$.
\end{lem}
\bpf Recalling that $h(t)-g(t)$ is nondecreasing with respect to $t$ (see Theorem \ref{thm1.1}), we conclude that $\lambda_1^{F}(t_1)\geq \lambda_1^{F}(t_2)$ if $t_1<t_2$ by Lemma \ref{lem2.2}.

Notice that the free boundaries are strictly monotone in the warm season. If $0\leq \delta <1$,
$g(t)>g(t+\omega)$ and $h(t)<h(t+\omega)$ since $[t, t+\omega]$ contains some warm days ($\delta=1$, all days are cold), so $\lambda_1^{P}([g(t), h(t)])>\lambda_1^{P}([g(t+\omega), h(t+\omega)])$, that is, $\lambda_1^{F}(t)>\lambda_1^{F}(t+\omega)$.
\epf

\begin{lem}
\label{lem2.16}
Assume $\mathbf{(J)}$ and $J_1(x)=J_2(x)$ in $x\in \mathbb{R}$ hold. The generalized principal eigenvalue of \eqref{b01} satisfied $\overline \lambda_1^P([-L_1, L_2])\geq 0$ for any $L_1$ and $L_2$
provided that $\underline \lambda_1^O\geq 0$.
\end{lem}
\bpf
 Let $(\underline\lambda_1^O, \phi(t), \psi(t))$ be the generalized principal eigen-pair of eigenvalue problem \eqref{b13} with $k_1=b_1$ and $k_2=k$, that is
{\small
\be\left\{
\begin{array}{ll}
\phi_t \geq a_1e_1\psi -b_1 \phi+\underline\lambda_1^O\phi, &0<t\leq (1-\delta)\omega,\\[2mm]
\psi_t\geq a_2e_2\phi -b_2 \psi+\underline\lambda_1^O\psi,  &0<t\leq (1-\delta)\omega, \\[2mm]
\phi_t\geq -b_1 \phi+\underline\lambda_1^O\phi, &(1-\delta)\omega<t\leq \omega,  \\[2mm]
\psi_t\geq -k\psi+\underline\lambda_1^O\psi,  &(1-\delta)\omega<t\leq \omega,\\[2mm]
\phi(0)=\phi(\omega), \psi(0)=\psi(\omega).
\end{array} \right.
\label{e17}\ee}

We notice that the kernel functions $J_1$ and $J_2$ satisfy $J_1(x)=J_2(x)(:= J(x))$ in $\mathbb{R}$,  problem \eqref{ed} admits a principal eigen-pair
$(\lambda^*_1, g(x))$ with $\lambda^*_1([-L_1, L_2])<0$, $g\in C[-L_1, L_2]$ and $g(x)>0$ in $[-L_1, L_2]$.
Let $w(t,x)=\phi(t)g(x), v(t,x)=\psi(t)g(x)$, then
{\small
\be\left\{
\begin{array}{ll}
w_t -d_1\mathcal{L}_1[w]\geq a_1e_1v -b_1 w+(-d_1\lambda^*_1+\underline\lambda_1^O)w, &0<t\leq (1-\delta)\omega, -L_1\leq x\leq L_2, \\[2mm]
v_t-d_2\mathcal{L}_2[v]\geq a_2e_2w -b_2 v+(-d_2\lambda^*_1+\underline\lambda_1^O)v,  &0<t\leq (1-\delta)\omega, -L_1\leq x\leq L_2,\\[2mm]
w_t-d_1\mathcal{L}_1[w]\geq -b_1 w+(-d_1\lambda^*_1+\underline\lambda_1^O)w, &(1-\delta)\omega<t\leq \omega, -L_1\leq x\leq L_2, \\[2mm]
v_t\geq -kv+\underline\lambda_1^Ov,  &(1-\delta)\omega<t\leq \omega, -L_1\leq x\leq L_2,\\[2mm]
w(0, x)=w(\omega, x), v(0, x)=v(\omega, x), &-L_1\leq x\leq L_2.
\end{array} \right.
\label{e181}\ee}
It follows from Lemma \ref{lem2.1} (i) that $\overline \lambda_1^P\geq \underline\lambda_1^O$ since $\lambda^*_1<0$.
\epf

\begin{lem}
\label{lem2.3}
Assume $\mathbf{(J)}$ and $J_1(x)=J_2(x)$ in $x\in \mathbb{R}$ hold. If $0\leq \delta <1$, then the principal eigenvalue $\lambda_1^P([-L_1, L_2])(-L_1<L_2)$ of \eqref{b01} has the following properties:

$(i)$ If $\lambda_1^O\geq 0$, then we have  $ \lambda_1^P([-L_1, L_2])>0$ for any $L_1$ and $L_2$.

$(ii)$ If $\lambda_1^O< 0$ and $\lambda_1^{P}([-h_0, h_0])\leq 0$,
then $\lambda_1^P([-L_1, L_2])<0$ for any $L_1$ and $L_2$.

$(iii)$ If $\lambda_1^O< 0$ and $\lambda_1^{P}([-h_0, h_0])> 0,$ then there exists $t^*>0$ such that $\lambda_1^{P}([g(t^*), h(t^*)])= 0$
and $(t-t^*)\lambda_1^{P}([g(t), h(t)])<0$ for $t\in(0, t^*-\omega)\bigcup(t^*+\omega, \infty)$.
\end{lem}
\bpf
$(i)$ If $0\leq\delta <1$, we have  $\overline{\lambda}_1^P=\underline{\lambda}_1^P={\lambda}_1^P$, and
$\overline\lambda_1^O=\lambda_1^O$.
Under the assumption of $J_1(x)=J_2(x)(:= J(x))$ in $\mathbb{R}$,  problem \eqref{ed} admits a principal eigen-pair
$(\lambda^*_1, g(x))$ with $\lambda^*_1([-L_1, L_2])<0$, $g\in C[-L_1, L_2]$ and $g(x)>0$ in $[-L_1, L_2]$.
Let $w(t,x)=\phi(t)g(x), v(t,x)=\psi(t)g(x)$, where $(\lambda_1^O, \phi(t), \psi(t))$  is the principal eigen-pair of eigenvalue problem \eqref{b13} with $k_1=b_1$ and $k_2=k$, then
{\footnotesize
\be\left\{
\begin{array}{llll}
w_t -d_1\mathcal{L}_1[w]= a_1e_1v -b_1 w+(-d_1\lambda^*_1+\lambda_1^O)w\\[2mm]
\qquad \qquad \qquad \quad >a_1e_1v -b_1 w+\lambda_1^Ow, &0<t\leq (1-\delta)\omega, -L_1\leq x\leq L_2, \\[2mm]
v_t-d_2\mathcal{L}_2[v]=a_2e_2w -b_2 v+(-d_2\lambda^*_1+\lambda_1^O)v\\[2mm]
\qquad \qquad \qquad \ \ >a_2e_2w -b_2 v+\lambda_1^Ov,  &0<t\leq (1-\delta)\omega, -L_1\leq x\leq L_2,\\[2mm]
w_t-d_1\mathcal{L}_1[w]= -b_1 w+(-d_1\lambda^*_1+\lambda_1^O)w> -b_1 w+\lambda_1^Ow, &(1-\delta)\omega<t\leq \omega, -L_1\leq x\leq L_2, \\[2mm]
v_t=-kv+\lambda_1^Ov,  &(1-\delta)\omega<t\leq \omega, -L_1\leq x\leq L_2,\\[2mm]
w(0, x)=w(\omega, x), v(0, x)=v(\omega, x), &-L_1\leq x\leq L_2.
\end{array} \right.
\label{e1811}\ee}
It follows from Lemma \ref{lem2.1} (i) that $\lambda_1^P>\lambda_1^O$ since $\lambda^*_1<0$.

$(ii)$ This result is obvious from  Lemma \ref{lem2.15}.

$(iii)$  It follows from Lemma \ref{lem2.15} that $\lambda_1^{F}(t)$ is nonincreasing for $t>0$, therefore there exists $t^*>0$ such that $\lambda_1^{P}([g(t^*), h(t^*)])= 0$.
 Recalling that $\lambda_1^{P}([g(t), h(t)])>\lambda_1^{P}([g(t+\omega), h(t+\omega)])$ by Lemma \ref{lem2.15}, we have
\bess
\lambda_1^{P}([g(t), h(t)])\left\{
\begin{array}{ll}
>0,\  &t\in(0, t^*-\omega),\\[2mm]
<0, &t\in(t^*+\omega, \infty).
\end{array} \right.
\eess
\epf

\section{A fixed boundary problem with seasonal succession}

In this section, we first present some known results, such as maximum principle (Lemma \ref{lem3.3}) and comparison principle (Lemma \ref{lem3.10}), as well as give some conclusions for the corresponding fixed boundary problem with seasonal succession (Lemmas \ref{lem3.1} and \ref{lem3.2}).

\begin{lem}
\label{lem3.3}
$($Maximum principle [\cite{Du2020} Lemma 3.1]\,$)$
Assume that $\mathbf{(J)}$ holds, and $g\in\mathbb{G}_{h_0,\omega}, h\in\mathbb{H}_{h_0,\omega}$ for some $\omega, h_0>0$. Suppose that $u_i, u_{it}\in C(\Omega^{g, h}), d_i, D_i, c_{ij}, C_{ij}\in L^{\infty}(\Omega^{g, h}), d_i\geq 0$, $D_i\geq 0$ ($i,j=1,2$) and
{\small
\bess
\left\{
\begin{array}{ll}
u_{it}\geq d_i\mathcal{L}_i[u_i]+\sum\limits_{j=1}^{2}c_{ij}u_j
 &m\omega<t\leq m\omega+(1-\delta)\omega,\, g(t)<x<h(t), \\[2mm]
 u_{it}\geq D_i\mathcal{L}_i[u_i]+\sum\limits_{j=1}^{2}C_{ij}u_j
 &m\omega+(1-\delta)\omega<t\leq (m+1)\omega, \,  g(t)<x<h(t),  \\[2mm]
u_i(t, g(t))\geq 0, u_i(t, h(t))\geq 0,  &t>0,\\[2mm]
u_i(0,x)\geq 0, i=1, 2,   &x\in[-h_0,h_0],
\end{array} \right.
\eess}
where $m=0, 1, \cdots$. Then the following conclusions hold:

$(i)$ If $c_{ij}, C_{ij}\geq 0(i\neq j)$ for $i,j={1,2}$, then $u_i\geq 0$ in $\Omega^{g, h}$.

$(ii)$ If $d_{i_0}>0$, $0\leq \delta<1$ and $u_{i_0}(0,x)\not\equiv 0$ in $[-h_0, h_0]$, then $u_{i_0}>0$ in $\Omega^{g, h}$.
\end{lem}

Applying maximum principle, we obtain comparison principal of the problem with seasonal succession, which will be used in the sequel.

\begin{lem}
\label{lem3.10}
Assume $\mathbf{(J)}$ holds, and $g\in\mathbb{G}_{h_0,\omega}, h\in\mathbb{H}_{h_0,\omega}$ for some $\omega, h_0>0$. Suppose that $u_i, \tilde{u}_i\in C(\Omega^{g, h})$ satisfy

$(i)$ $u_{it}, \tilde{u}_{it}\in C(\Omega^{g, h})$.

$(ii)$  $0<u_i\leq e_i, 0<\tilde{u}_i\leq e_i$.

$(iii)$ $(\tilde{u}_1, \tilde{u}_2)$ satisfies
{\small
\be
\left\{
\begin{array}{ll}
\tilde{u}_{1t}\geq d_1\mathcal{L}_1[\tilde{u}_1]+a_1(e_1-\tilde{u}_1)\tilde{u}_2-b_1\tilde{u}_1,
 &(t,x)\in (m\omega, m\omega+(1-\delta)\omega]\times(g(t), h(t)), \\[2mm]
\tilde{u}_{2t}\geq d_2\mathcal{L}_2[\tilde{u}_2]+a_2(e_2-\tilde{u}_2)\tilde{u}_1-b_2\tilde{u}_2,
 &(t,x)\in (m\omega, m\omega+(1-\delta)\omega]\times(g(t), h(t)), \\[2mm]
\tilde{u}_{1t}\geq d_1\mathcal{L}_1[\tilde{u}_1]-b_1\tilde{u}_1,
 &(t,x)\in (m\omega+(1-\delta)\omega, (m+1)\omega]\times(g(t), h(t)), \\[2mm]
\tilde{u}_{2t}\geq -k\tilde{u}_2,
 &(t,x)\in (m\omega+(1-\delta)\omega, (m+1)\omega]\times(g(t), h(t)).
\end{array} \right.
\label{c10}
\ee
}

$(iv)$ $(u_1, u_2)$ satisfies \eqref{c10} but with the inequalities reversed.

$(v)$ At the boundary,
\bess
u_i(t, x)\leq \tilde{u}_i(t, x)\quad {\rm for}\,\,  t>0, x=\{g(t), h(t)\}.
\eess

$(vi)$ At the initial time,
\bess
u_i(0,x)\leq \tilde{u}_i(0,x)\quad {\rm for}\,\, x\in[-h_0,h_0].
\eess

Then
\bess
u_i(t,x)\leq \tilde{u}_i(t,x)\quad {\rm for}\  (t,x)\in \Omega^{g, h}, i=1,2.
\eess
\end{lem}
\bpf
Define
\bess
H:=\tilde{u}_1-u_1, V:=\tilde{u}_2-u_2,
\eess
and
\bess
c_{11}:=-(b_1+a_1u_2),\  c_{12}:=a_1(e_1-\tilde{u}_1)>0, \ c_{21}:=a_2(e_2-\tilde{u}_2)>0,\  c_{22}:=-(b_2+a_2u_1)
\eess for $t\in(m\omega, m\omega+(1-\delta)\omega]$,

\bess
C_{11}:=-b_1,\  C_{12}:=0,\  C_{21}:=0,\  C_{22}:=-k, D_1:=d_1,\ D_2:=0
\eess for $t\in(m\omega+(1-\delta)\omega, (m+1)\omega]$.
It follows from Lemma \ref{lem3.3} that $H\geq 0$ and $V\geq 0$ in $\Omega^{g, h}$.
\epf

\vspace{3mm}
For our later purpose, we consider the corresponding fixed boundary problem with seasonal succession of \eqref{a07}:
{\small
\be\left\{
\begin{array}{ll}
u_{1t}=d_1\int_{-L_1}^{L_2}J_1(x-y)u_1(t,y)dy-d_1u_1
+a_1(e_1-u_1)u_2 -b_1 u_1, &(t,x)\in Q^m_{\textrm w}, \\[2mm]
u_{2t}=d_2\int_{-L_1}^{L_2}J_2(x-y)u_2(t,y)dy-d_2u_2
+a_2(e_2-u_2) u_1 -b_2 u_2,   &(t,x)\in Q^m_{\textrm w},\\[2mm]
u_{1t}=d_1\int_{-L_1}^{L_2}J_1(x-y)u_1(t,y)dy-d_1u_1-b_1 u_1, &(t,x)\in Q^m_{\textrm c}, \\[2mm]
u_{2t}=-k u_2, &(t,x)\in Q^m_{\textrm c}, \\[2mm]
u_i(0,x)=u_{i,0}(x),\, i=1,2,   &-L_1\leq x\leq L_2,
\end{array} \right.
\label{c01}
\ee
}
where $Q^m_{\textrm w}:=(m\omega, m\omega+(1-\delta)\omega]\times(-L_1, L_2)$, $Q^m_{\textrm c}:=(m\omega+(1-\delta)\omega, (m+1)\omega]\times(-L_1, L_2)$, $m=0, 1, \cdots$, $u_{i,0}(x)\in C([-L_1, L_2])\backslash \{0\}$ and $ 0\leq u_{i,0}(x)\leq e_i (i=1,2)$. It is well-known that \eqref{c01}
has a unique positive solution which is defined for all $t>0$.
\begin{defi}
A pair of functions $(\tilde{u}_1, \tilde{u}_2)(t, x)\in [C^{1,0}((0, \infty)\times[-L_1, L_2])]^2$ is called an upper  solution to problem \eqref{c01} if
{\small
\bess
\left\{
\begin{array}{ll}
\tilde{u}_{1t}\geq d_1\int_{-L_1}^{L_2}J_1(x-y)\tilde{u}_1(t,y)dy-d_1\tilde{u}_1
+a_1(e_1-\tilde{u}_1)\tilde{u}_2 -b_1 \tilde{u}_1, &(t,x)\in Q^m_{\textrm w}, \\[2mm]
\tilde{u}_{2t}\geq d_2\int_{-L_1}^{L_2}J_2(x-y)\tilde{u}_2(t,y)dy-d_2\tilde{u}_2
+a_2(e_2-\tilde{u}_2) \tilde{u}_1 -b_2 \tilde{u}_2,   &(t,x)\in Q^m_{\textrm w},\\[2mm]
\tilde{u}_{1t}\geq d_1\int_{-L_1}^{L_2}J_1(x-y)\tilde{u}_1(t,y)dy-d_1\tilde{u}_1-b_1 \tilde{u}_1, &(t,x)\in Q^m_{\textrm c}, \\[2mm]
\tilde{u}_{2t}\geq -k \tilde{u}_2, &(t,x)\in Q^m_{\textrm c}, \\[2mm]
\tilde{u}_i(0,x)\geq u_{i,0}(x),\,i=1,2,   &-L_1\leq x\leq L_2.
\end{array} \right.
\eess
}
\end{defi}
Similarly, we call $(\hat{u}_1, \hat{u}_2)(t, x)$ a lower solution of \eqref{c01} if all reversed inequalities in the above places are satisfied.

To understand the asymptotic behavior of solution to the initial value problem \eqref{c01},
we now consider the corresponding periodic problem of \eqref{c01}
{\small
\be\left\{
\begin{array}{ll}
U_{1t}=d_1\int_{-L_1}^{L_2}J_1(x-y)U_1(t,y)dy-d_1U_1
+a_1(e_1-U_1)U_2 -b_1 U_1, &(t,x)\in Q^0_{\textrm w}, \\[2mm]
U_{2t}=d_2\int_{-L_1}^{L_2}J_2(x-y)U_2(t,y)dy-d_2U_2
+a_2(e_2-U_2) U_1 -b_2 U_2,  &(t,x)\in Q^0_{\textrm w},\\[2mm]
U_{1t}=d_1\int_{-L_1}^{L_2}J_1(x-y)U_1(t,y)dy-d_1U_1-b_1 U_1, &(t,x)\in Q^0_{\textrm c}, \\[2mm]
U_{2t}=-k U_2, &(t,x)\in Q^0_{\textrm c}, \\[2mm]
U_i(0,x)=U_i(\omega,x), \,i=1,2,   &-L_1\leq x\leq L_2.
\end{array} \right.
\label{c02} \ee
}
\begin{lem}
\label{lem3.1}
Suppose $\mathbf{(J)}$ holds, let $(u_1(t, x), u_2(t, x))$ be the unique positive solution of \eqref{c01}, $\underline{\lambda}_1^P([-L_1, L_2])$ and $ \overline{\lambda}_1^P([-L_1, L_2])$ are the generalized principal eigenvalue of \eqref{b01}. The following conclusions hold.

$(i)$ If $\underline{\lambda}_1^P([-L_1, L_2])\geq 0$, then problem \eqref{c02} has the only nonnegative solution $(0, 0)$, and
$\lim\limits_{t\rightarrow\infty}(u_1(t, x), u_2(t, x))=(0,0)$ uniformly for $x\in[-L_1, L_2]$.

$(ii)$ If $\overline{\lambda}_1^P([-L_1, L_2])<0$, then problem \eqref{c02} admits a unique positive periodic solution $(U^*_1(t,x), U^*_2(t,x))$, and $0<U^*_{i}\leq e_i (i=1,2)$. Moreover,
$\lim\limits_{n\rightarrow\infty}(u_1(t+n\omega, x), u_2(t+n\omega, x))=(U^*_1(t,x), U^*_2(t,x))$ uniformly for $x\in[-L_1, L_2]$.
\end{lem}
\bpf
Suppose $\overline{\lambda}_1^P([-L_1, L_2])<0$, we first consider the existence of positive periodic solution $(U^*_1(t,x), U^*_2(t,x))$ to problem \eqref{c02}. 
Let $(\phi, \psi)$ be a positive eigenfunction pair with respect to $\overline{\lambda}_1^P([-L_1, L_2])$, we are in a position to check that for a small $\varepsilon>0$, $(\varepsilon\phi, \varepsilon\psi)$ and $(e_1, e_2)$ are the lower and upper solutions of \eqref{c02}.

Let us denote by $(\overline{U}_1, \overline{U}_2)(t, x)$  the unique positive solution of \eqref{c01} with the initial function pair $(e_1, e_2)$, our next goal is to prove that $(\overline{U}_1, \overline{U}_2)$ is nonincreasing in $t$. In fact, considering $(e_1, e_2)$ as an upper solution of \eqref{c01}, we apply Lemma \ref{lem3.10}, where $(g(t), h(t))\equiv(-L_1, L_2)$, to assert that
$$\overline{U}_i(t,x)\leq e_i, \ \  {\rm for}\ (t,x)\in Q^0_{\textrm w}\cup Q^0_{\textrm c}, \ i=1,2.$$
 It follows from the definition of $(\overline{U}_1, \overline{U}_2)$ and Lemma \ref{lem3.10}  that, for any $t_1>0$,
$$\overline{U}_i(t+t_1, x)\leq \overline{U}_i(t,x)\ \  {\rm in} \ Q^0_{\textrm w}\cup Q^0_{\textrm c},$$
which gives the monotonicity in $t$.
Meanwhile, by regarding $(\varepsilon\phi, \varepsilon\psi)$ as a lower solution of \eqref{c01}, we can deduce that
$$\overline{U}_1(t,x)\geq \varepsilon\phi\ \  {\rm and}\ \ \overline{U}_2(t,x)\geq \varepsilon\psi\ \ {\rm in}\ Q^0_{\textrm w}\cup Q^0_{\textrm c}.$$
As a consequence, the limitation of $(\overline{U}_1, \overline{U}_2)(t, x)$ exists, denoted by $(U^*_1, U^*_2)(t, x)$, which is a positive periodic solution of \eqref{c02} by the dominated convergence theorem.
Furthermore, for any positive periodic solution $(U_1, U_2)(t, x)$ to \eqref{c02} satisfying $U_i(t, x)\leq e_i$,  Lemma \ref{lem3.10} gives that $\overline U_i(t, x)\geq U_i(t, x)$ for $(t,x)\in Q^0_{\textrm w}\cup Q^0_{\textrm c}$, that is to say, $(U^*_1, U^*_2)(t, x)$ is the maximal positive solution of \eqref{c02}.

Next, we proceed to show the uniqueness of the positive periodic solution of \eqref{c02}. Assume $(U_{11}, U_{21})$ and $(U_{12}, U_{22})$ are two positive periodic solutions, there is no loss of generality in assuming $(U_{12}, U_{22})=(U^*_1, U^*_2)$, we have
\bess
(0,0)<(U_{11}, U_{12})\leq, \not\equiv (U_{12}, U_{22}),
\eess
and it follows from \eqref{c02} that
{\small
\bess
\begin{array}{rll}
\langle U_{11t}-d_1\mathcal{L}_1[U_{11}],\, U_{12}\rangle&=&\langle a_1(e_1-U_{11})U_{21} -b_1 U_{11},\,  U_{12}\rangle,\\[2mm]
\langle U_{21t}-d_2\mathcal{L}_2[U_{21}],\,  U_{22}\rangle&=&\langle a_2(e_2-U_{21})U_{11} -b_2 U_{21},\,  U_{22}\rangle,\\[2mm]
\langle U_{12t}-d_1\mathcal{L}_1[U_{12}],\,  U_{11}\rangle&=&\langle a_1(e_1-U_{12})U_{22} -b_1 U_{12},\,  U_{11}\rangle,\\[2mm]
\langle U_{22t}-d_2\mathcal{L}_2[U_{22}],\,  U_{21}\rangle&=&\langle a_2(e_2-U_{22})U_{12} -b_2 U_{22},\,  U_{21}\rangle
\end{array}
\eess
}
hold in warm season ($t\in (0,(1-\delta)\omega]$), where $\langle \cdot,\, \cdot\rangle$ is defined in \eqref{f05} and
{\small
\bess
\begin{array}{rll}
\langle U_{11t}-d_1\mathcal{L}_1[U_{11}],\,  U_{12}\rangle&=&\langle -b_1 U_{11},\,  U_{12}\rangle,\\[2mm]
\langle U_{21t},\,  U_{22}\rangle&=&\langle -k U_{21},\,  U_{22}\rangle,\\[2mm]
\langle U_{12t}-d_1\mathcal{L}_1[U_{12}],\,  U_{11}\rangle&=&\langle -b_1 U_{12},\,  U_{11}\rangle,\\[2mm]
\langle U_{22t},\,  U_{21}\rangle&=&\langle -k U_{22},\,  U_{21}\rangle
\end{array}
\eess
}
for $t\in ((1-\delta)\omega, \omega]$.
Using $\langle\mathcal{L}_1[U_{i1}], U_{i2}\rangle=\langle\mathcal{L}_1[U_{i2}], U_{i1}\rangle (i=1,2)$ and
transferring items yield that
{\small
\be
\begin{array}{rll}
\langle U_{11t},\, U_{12}\rangle+a_1\langle(e_1-U_{12})U_{22},\,  U_{11}\rangle&=&\langle U_{12t},\,  U_{11}\rangle+a_1\langle(e_1-U_{11})U_{21},\,  U_{12}\rangle,\\[2mm]
\langle U_{21t},\,  U_{22}\rangle+a_2\langle(e_2-U_{22})U_{12},\,  U_{21}\rangle&=&\langle U_{22t},\,  U_{21}\rangle+a_2\langle(e_2-U_{21})U_{11},\,  U_{22}\rangle\\[2mm]
\end{array}
\label{c03}
\ee
}
for $t\in (0,(1-\delta)\omega]$, and
{\small
\be
\begin{array}{rll}
\langle U_{11t},\,  U_{12}\rangle&=&\langle U_{12t},\,  U_{11}\rangle,\\[2mm]
\langle U_{21t},\,  U_{22}\rangle&=&\langle U_{22t},\,  U_{21}\rangle
\end{array}
\label{c0311}
\ee
}
for $t\in ((1-\delta)\omega, \omega]$.
Integrating the first equation of \eqref{c03} over $(0,(1-\delta)\omega)$, the fist equation of \eqref{c0311} over $((1-\delta)\omega, \omega)$, respectively, then adding them together give
{\small
\bess
\begin{array}{rll}
&&a_1\int_{0}^{(1-\delta)\omega}e_1[\langle U_{21},\,  U_{12}\rangle-\langle U_{22},\,  U_{11}\rangle]-\langle U_{11}U_{21},\,  U_{12}\rangle+\langle U_{12}U_{22},\,  U_{11}\rangle dt\\[3mm]
&=&\int_{0}^{\omega}\langle U_{11t},\,  U_{12}\rangle dt-\int_{0}^{\omega}\langle U_{12t},\,  U_{11}\rangle dt=0
\end{array}
\eess
}
since $U_{11}(\omega,x)=U_{11}(0,x)$ and $U_{12}(\omega,x)=U_{12}(0,x)$ in \eqref{c02}. Therefore,
\be
\int_{0}^{(1-\delta)\omega} [\langle U_{21},\,  U_{12}\rangle-\langle U_{22},\,  U_{11}\rangle]dt=\int_{0}^{(1-\delta)\omega}\frac{1}{e_1} \langle U_{11}U_{12},\,  U_{21}-U_{22}\rangle dt.
\label{e18}\ee
The same procedure for the second equation of \eqref{c03} and that of \eqref{c0311} may be easily adapted to obtain
\be
\int_{0}^{(1-\delta)\omega} [\langle U_{11},\,  U_{22}\rangle-\langle U_{12},\,  U_{21}\rangle]dt=\int_{0}^{(1-\delta)\omega} \frac{1}{e_2} \langle U_{21}U_{22},\,  U_{11}-U_{12}\rangle dt.
\label{e19}\ee
Adding \eqref{e18} and \eqref{e19} yields
\bess
0=\int_{0}^{(1-\delta)\omega}[\frac{1}{e_1}\langle U_{11}U_{12},\, U_{21}-U_{22}\rangle+\frac{1}{e_2}\langle U_{21}U_{22},\, U_{11}-U_{12}\rangle] dt<0
\eess
since that $(0, 0)<(U_{11}, U_{12})\leq, \not\equiv (U_{12}, U_{22})$, we then lead to
 a contradiction. This finishes the proof of the uniqueness of the positive periodic solution.

We notice that $\int_{-L_1}^{L_2}J_i(x-y)U^*_i(t,y)dy:=G_i (i=1,2)$ are continuous in $[0,\omega]\times[-L_1, L_2]$, it follows from \eqref{c02} that $U^*_1$ and $U^*_2$ can be expressed by a quadratic formula involving $G_i$ and the constant parameters. Therefore both $U^*_1$ and $U^*_2$ are continuous functions in $[0,\omega]\times[-L_1, L_2]$.

Finally, we claim that \eqref{c02} has no positive periodic solution if $\underline{\lambda}_1^P([-L_1, L_2])\geq 0$. Suppose, contrary to our claim, that $(\tilde{U}_1, \tilde{U}_2)$ is a positive periodic solution to problem \eqref{c02} and satisfies
\be
\left\{
\begin{array}{ll}
\tilde{U}_{1t}-d_1\mathcal{L}_1[\tilde{U}_1]=
a_1(e_1-\tilde{U}_1)\tilde{U}_2 -b_1 \tilde{U}_1<a_1e_1\tilde{U}_2 -b_1 \tilde{U}_1, &(t,x)\in Q^0_{\textrm w}, \\[2mm]
\tilde{U}_{2t}-d_2\mathcal{L}_2[\tilde{U}_2]=
a_2(e_2-\tilde{U}_2) \tilde{U}_1 -b_2 \tilde{U}_2<a_2e_2\tilde{U}_1 -b_2 \tilde{U}_2,  &(t,x)\in Q^0_{\textrm w},\\[2mm]
\tilde{U}_{1t}-d_1\mathcal{L}_1[\tilde{U}_1]=-b_1 \tilde{U}_1, &(t,x)\in Q^0_{\textrm c}, \\[2mm]
\tilde{U}_{2t}=-k \tilde{U}_2, &(t,x)\in Q^0_{\textrm c}, \\[2mm]
\tilde{U}_i(0,x)=\tilde{U}_i(\omega,x), i=1,2,   &-L_1\leq x\leq L_2.
\end{array} \right.
\label{e20}
\ee
It follows from the definition of the generalized principal eigenvalue that $\overline{\lambda}_1^P([-L_1, L_2])\leq 0$. If $\overline{\lambda}_1^P([-L_1, L_2])\neq \underline{\lambda}_1^P([-L_1, L_2])$, we then have $\underline{\lambda}_1^P([-L_1, L_2])<0$, which is in contradiction with the assumption $\underline{\lambda}_1^P([-L_1, L_2])\geq 0$. If $\overline{\lambda}_1^P([-L_1, L_2])= \underline{\lambda}_1^P([-L_1, L_2])$, then ${\lambda}_1^P([-L_1, L_2])$ is well-defined, and
${\lambda}_1^P([-L_1, L_2])<0$ by applying Lemma \ref{lem2.1} to the system \eqref{e20}, we also lead to a contradiction.

Through the above discussion, we conclude that problem \eqref{c02} has a unique positive periodic solution if $\overline{\lambda}_1^P([-L_1, L_2])<0$, while when $\underline{\lambda}_1^P([-L_1, L_2])\geq 0$, $(0,0)$ is the only nonnegative periodic solution. Next we prove the stabilities in this theorem.

For $(i)$, since $\underline{\lambda}_1^P([-L_1, L_2])\geq 0$, we have $(U^*_1, U^*_2)=(0,0)$. And then
$\lim\limits_{t\rightarrow\infty}(\overline{U}_1, \overline{U}_2)=(0,0)$ uniformly in $[-L_1, L_2]$, which is due to Dini's theorem. We conclude from Lemma \ref{lem3.10} that $0<u_i(t,x)\leq \overline{U}_i(t,x)$ for $(t,x)\in Q^0_{\textrm w}\cup Q^0_{\textrm c}$, hence that $(u_1, u_2)\rightarrow(0,0)$ as $t\rightarrow\infty$.

For $(ii)$, if $\overline{\lambda}_1^P([-L_1, L_2])<0$, problem \eqref{c02} has a unique positive periodic solution $(U^*_1, U^*_2)$.
According to the fact that $(\overline{U}_1, \overline{U}_2)$ is nonincreasing in $t$ and $\overline{U}_i$ is bounded for $(t,x)\in Q^0_{\textrm w}\cup Q^0_{\textrm c}$. Hence $\lim\limits_{n\rightarrow\infty}(\overline{U}_1, \overline{U}_2)(t+n\omega, x)=(U^*_1, U^*_2)(t, x)$
uniformly for $x\in[-L_1, L_2]$ by Dini's theorem.

To find a positive lower solution, now we reset a initial time so that the initial value is positive. Note that $U_i(1,x)>0$ for $x\in[-L,L]$, for the small enough $\varepsilon>0$, $U_1(1,x)\geq\varepsilon\phi, \,U_2(1,x)\geq\varepsilon\psi\,\,\,\textrm{in}\,\,[-L_1, L_2].$
Let $(\underline{U}_1, \underline{U}_2)$ be the solution of \eqref{c01}
with initial function pair $(\varepsilon\phi, \varepsilon\psi)$. Then $\underline{U}_i$ are nondecreasing in $t$ and by Lemma \ref{lem3.10}, we admit
$\underline{U}_i(t,x)\leq U_i(1+t,x)\leq \overline{U}_i(1+t,x)$ in $Q^0_{\textrm w}\cup Q^0_{\textrm c}$. Apply similar above arguments to this case, we obtain that $\lim\limits_{n\rightarrow\infty}(\underline{U}_1, \underline{U}_2)(t+n\omega, x)=({U}^*_1, {U}^*_2)(t, x)$ uniformly for $x\in[-L_1, L_2]$. This completes the proof.
\epf

\begin{rmk} It follows from Corollary \ref{rmk1.2} that $\overline{\lambda}_1^P(\delta)$ and $\underline{\lambda}_1^P(\delta)$ are strictly increasing in $\delta$.
That is, if $\delta_1<\delta_2$, then $\underline{\lambda}_1^P(\delta_1)<\underline{\lambda}_1^P(\delta_2)$ and $\overline{\lambda}_1^P(\delta_1)<\overline{\lambda}_1^P(\delta_2)$.
Recalling that $\overline{\lambda}_1^P(1)\geq \underline{\lambda}_1^P(1)>0$
by Theorem \ref{thm2.1} (ii), so
 there exist a  $\delta_2^*$ such that $\underline{\lambda}_1^P(\delta_2)\geq 0$ for $\delta\geq \delta^*_2$,  which gives  $$\lim\limits_{t\rightarrow\infty}(u_1(t, x), u_2(t, x))=(0,0)$$
by Lemma \ref{lem3.1}.
That is to say, the bigger $\delta$ is, the more likely $(u_1, u_2)$ is to become extinct.
On the other hand,  for small enough $\delta_1$, we have $\overline{\lambda}_1^P(\delta_1)<0$,  it follows from  Lemma \ref{lem3.1} that
\bess
\lim\limits_{n\rightarrow\infty}(u_1(t+n\omega, x), u_2(t+n\omega, x))=(U^*_1(t,x), U^*_2(t,x)),
\eess
which implies that the smaller $\delta$ is, the more likely $(u_1, u_2)$ converges to a positive periodic solution $(U^*_1, U^*_2)$.
\label{rmk111}
\end{rmk}

\vspace{3mm}
For problem \eqref{c02}, we consider the corresponding spatial-independent problem
\be
\left\{
\begin{array}{ll}
U_{1t}=a_1(e_1-U_1)U_2 -b_1 U_1, &t\in (0, (1-\delta)\omega], \\[2mm]
U_{2t}=a_2(e_2-U_2) U_1 -b_2 U_2,  &t\in (0, (1-\delta)\omega],\\[2mm]
U_{1t}= -b_1 U_1, &t\in ((1-\delta)\omega, \omega], \\[2mm]
U_{2t}=-k U_2, &t\in ((1-\delta)\omega, \omega], \\[2mm]
U_i(0)=U_i(\omega)(i=1,2).
\end{array} \right.
\label{c11}
\ee

 Similarly to the proof of Lemma \ref{lem3.1}, we can easily obtain the existence and uniqueness of the positive periodic solution to problem \eqref{c11}.

\begin{lem}
\label{lem3.2}
Suppose $\mathbf{(J)}$ and $J_1(x)=J_2(x)$ hold.  If $0\leq \delta<1$ and $\lambda_1^O<0$, then
\bess
\lim\limits_{L_1, L_2\rightarrow+\infty}(U^*_{1, [-L_1, L_2]}, U^*_{2, [-L_1, L_2]})(t, x) =(U_1^{\vartriangle}, U_2^{\vartriangle})(t)
\eess
uniformly for $t\in [0, \omega]$ and locally uniformly for $x\in\mathbb{R}$, where  $(U^*_{1, [-L_1, L_2]}, U^*_{2, [-L_1, L_2]})(t, x)$ and $(U_1^{\vartriangle}, U_2^{\vartriangle})(t)$ are the unique positive periodic solution of problem \eqref{c02} and \eqref{c11}, respectively.
\end{lem}
\bpf
If $0\leq \delta<1$ and $\lambda_1^O<0$, it follows from Lemma \ref{lem2.3} that there exists $T_1$ large enough such that $\lambda_1^P([g(T_1), h(T_1)])<0$.
Let $L_0=\max\{-g(T_1), h(T_1)\}$, using the monotonicity in Lemma \ref{lem2.2} and the fact $\lambda_1^P([g(T_1), h(T_1)])<0$ gives that $\lambda_1^P([-L_0, L_0])<0$, which means that problem \eqref{c02} has a unique positive periodic solution $(U^*_{1, [-L, L]}, U^*_{2, [-L, L]})$ for all $L\geq L_0$.
To derive the desired result, the proof falls naturally into three parts.

\textbf{Part 1.}  Let us first prove that  for $L_0\leq L$, $U_1^*$ and  $U_2^*$ are nondecreasing in $L$.

Suppose
$L_0\leq L_1^*\leq L_2^*$, let $(U_{1[-L_j^*, L_j^*]}, U_{2, [-L_j^*, L_j^*]})$ be the positive solution of \eqref{c01} with $L_1=L_2=L_j^*$ and initial function pair $(U_{1, [-L_j^*, L_j^*]}(0,x), U_{2, [-L_j^*, L_j^*]}(0,x)) (j=1,2)$, and  $0<U_{i, [-L_1^*, L_1^*]}(0, x) \leq U_{i, [-L_2^*, L_2^*]}(0, x) \leq e_i (i=1,2).$
Thanks to
\bess
\int_{-L_2^*}^{L_2^*}J_1(x-y)U_{i, [-L_2^*, L_2^*]}(t, y)dy\geq \int_{-L_1^*}^{L_1*}J_1(x-y)U_{i, [-L_2^*, L_2^*]}(t, y)dy,
\eess
then $(U_{1, [-L_2^*, L_2^*]}, U_{2, [-L_2^*, L_2^*]})$ is an upper solution over the restriction of $[0, \omega]\times[-L_1^*, L_1^*]$,
as a result, $U_{i, [-L_2^*, L_2^*]}\geq U_{i, [-L_1^*, L_1^*]}$ in $[0, \omega]\times[-L_1^*, L_1^*]$ by applying comparison principle.
Because of $\lambda_1^P<0$, it follows from Lemma \ref{lem3.1} (ii) that
\bess
U^*_{1, [-L_1^*, L_1^*]}\leq U^*_{1, [-L_2^*, L_2^*]},\ U^*_{2, [-L_1^*, L_1^*]}\leq U^*_{2, [-L_2^*, L_2^*]}\,\,{\rm for}\,\,(t,x)\in(0,\infty)\times[-L_1^*, L_1^*].
\eess
The monotone property implies that the pointwise limit exists and
\bess
\lim\limits_{L\rightarrow+\infty}(U^*_{1, [-L, L]}, U^*_{2, [-L, L]})(t, x)=(U_1^*, U_2^*)(t, x).
\eess

Notice that $[-L^m, L^m]\subseteq[-L_1, L_2]\subseteq[-L^M, L^M]$ with $L^m=\min\{L_1, L_2\}$ and $L^M=\max\{L_1, L_2\}$, we have
\bess
\lim\limits_{L_1, L_2\rightarrow+\infty}(U^*_{1, [-L_1, L_2]}, U^*_{2, [-L_1, L_2]})(t,x)=(U_1^*, U_2^*)(t,x)
\eess
by using the monotonicity of positive periodic solution with respect to the interval. It is easily seen that $0<U_i^*(t,x)\leq e_i$ and $(U_1^*, U_2^*)$ is a positive solution of \eqref{c02} with $(-L_1,L_2)$ replaced by $(-\infty, +\infty)$ by the dominated convergence theorem.

\textbf{Part 2.}   We next prove that $U_i^*(t,x)$ is independent of $x$.
It suffices to show that
\be
(U_1^*, U_2^*)(t,x_1)=(U_1^*, U_2^*)(t,0)\,\,{\rm for\,\, any \,\,given}\,\, x_1\in \mathbb{R}.
\label{c04}
\ee
Denote $L^*:=|x_1|$. For $L\geq L_0+2L^*$,
since $$[-L+2L^*, L-2L^*]\subset[-L+L^*-x_1, L-L^*-x_1]\subset[-L,L],$$
it follows from Part 1 that
$$U^*_{i, [-(L-2L^*), L-2L^*]}(t,x)\leq U^*_{i,[-(L-L^*), L-L^*]}(t, x+x_1)\leq U^*_{i, [-L, L]}(t,x)\,(i=1,2).$$
Letting $L\rightarrow\infty$, by the definition of $(U_1^*, U_2^*)$, yields
$$(U_1^*, U_2^*)(t,x)=(U_1^*, U_2^*)(t,x+x_1).$$
Take $x=0$, \eqref{c04} is proved.

\textbf{Part 3.}
Finally, since that  $(U_1^{\vartriangle}, U_2^{\vartriangle})$ is the unique positive periodic solution to problem \eqref{c11}, the conclusion of  Part 2 implies $(U_1^*, U_2^*)=(U_1^{\vartriangle}, U_2^{\vartriangle})$, i.e.
\bess
\lim\limits_{L_1, L_2\rightarrow+\infty}(U^*_{1, [-L_1, L_2]}, U^*_{2, [-L_1, L_2]})(t, x)=(U_1^{\vartriangle}, U_2^{\vartriangle})(t).
\eess
The convergence is locally uniform in $\mathbb{R}$ by Dini's theorem.
\epf
\section{Main results}
We prove Theorems \ref{thm1.1}-\ref{thm1.4} in this section. We always assume (J) holds, the initial
function pair $(u_{1,0}(x), u_{2,0}(x))$ satisfies \eqref{a08}, and
$(u_1(t, x), u_2(t, x); g(t), h(t))$ is the unique positive solution of \eqref{a07}.
We also denote
\be
g_{\infty}=\lim\limits_{t\rightarrow\infty}g(t)\ \textrm{and}\ h_{\infty}=\lim\limits_{t\rightarrow\infty}h(t).
\label{d02}
\ee

\begin{lem}
\label{lem3.4}
(Comparison principle).
Suppose $\mathbf{(J)}$ holds and $(u_1(t, x), u_2(t, x); g(t), h(t))$ is the solution of \eqref{a07} with the initial functions satisfy \eqref{a08}. Let $\overline{g}, \overline{h}\in C([0, \infty))\cap C^1([m\omega, m\omega+(1-\delta)\omega])$. If $(\overline{u}_{1}, \overline{u}_{2})\in [C(\Omega^{\overline{g}, \overline{h}})]^2$ satisfies
{\footnotesize
\be\left\{
\begin{array}{ll}
\overline{u}_{1t}\geq d_1\int_{\overline{g}(t)}^{\overline{h}(t)}J_1(x-y)\overline{u}_1(t,y)dy
-d_1\overline{u}_1(t, x)\\[2mm]
\qquad \ +a_1(e_1-\overline{u}_1)\overline{u}_2 -b_1 \overline{u}_1, &t\in(m\omega, m\omega+(1-\delta)\omega], x\in(\overline{g}(t), \overline{h}(t)),\\[2mm]
\overline{u}_{2t}\geq d_2\int_{\overline{g}(t)}^{\overline{h}(t)}J_2(x-y)\overline{u}_2(t,y)dy
-d_2\overline{u}_2(t, x)\\[2mm]
\qquad \ +a_2(e_2-\overline{u}_2)\overline{u}_1 -b_2 \overline{u}_2, &t\in(m\omega, m\omega+(1-\delta)\omega], x\in(\overline{g}(t), \overline{h}(t)),\\[2mm]
\overline{u}_{1t}\geq d_1\int_{\overline{g}(t)}^{\overline{h}(t)}J_1(x-y)\overline{u}_1(t,y)dy
-d_1\overline{u}_1(t, x)-b_1\overline{u}_1, &t\in(m\omega+(1-\delta)\omega, (m+1)\omega], x\in(\overline{g}(t), \overline{h}(t)),\\[2mm]
\overline{u}_{2t}\geq -k\overline{u}_2, &t\in(m\omega+(1-\delta)\omega, (m+1)\omega], x\in(\overline{g}(t), \overline{h}(t)),\\[2mm]
\overline{u}_i(t, x)\geq 0,&t>0, x\in\{\overline{g}(t), \overline{h}(t)\},\\[2mm]
\overline{h}'(t)\geq\sum\limits_{i=1}^{2}\mu_i\int_{\overline{g}(t)}^{ \overline{h}(t)}\int_{\overline{h}(t)}^{+\infty}J_i(x-y)\overline{u}_i(t,x)dydx, &t\in(m\omega, m\omega+(1-\delta)\omega],\\[2mm]
\overline{g}'(t)\leq-\sum\limits_{i=1}^{2}\mu_i\int_{\overline{g}(t)}^{ \overline{h}(t)}\int_{-\infty}^{\overline{g}(t)}J_i(x-y)\overline{u}_i(t,x)dydx, &t\in(m\omega, m\omega+(1-\delta)\omega],\\[2mm]
\overline{h}(t)\geq \overline{h}(m\omega+(1-\delta)\omega),\,\, \overline{g}(t)\leq\overline{g}(m\omega+(1-\delta)\omega), &t\in(m\omega+(1-\delta)\omega, (m+1)\omega],\\[2mm]
\overline{g}(0)\leq-h_0, \overline{h}(0)\geq h_0, \overline{u}_i(0,x)\geq 0, &x\in[\overline{g}(0), \overline{h}(0)],\\[2mm]
\overline{u}_i(0,x)\geq \overline{u}_{i,0}(x),   &x\in[-h_0, h_0],
\end{array} \right.
\label{d01}\ee
}
then $[g(t), h(t)]\subset[\overline{g}(t), \overline{h}(t)]$ and
\bess
u_i(t, x)\leq \overline{u}_i(t, x), \quad (t, x)\in[0, +\infty)\times[g(t), h(t)].
\eess
If the inequalities in \eqref{d01} are revered, and $(\overline u_1, \overline u_2; \overline g, \overline h)$ is replaced by $(\underline u_1, \underline u_2; \underline g, \underline h)$, then $[g(t),$ $ h(t)]$ $\supset [\underline{g}(t), \underline{h}(t)]$ and
\bess
u_i(t, x)\geq \underline{u}_i(t, x), \quad (t, x)\in[0, +\infty)\times[\underline{g}(t), \underline{h}(t)].
\eess
\end{lem}
\bpf
The proof mainly follows the approach of [\cite{Cao2019} Theorem 3.1], we omit the details with minor modifications in the cold season.
\epf

\vspace{2mm}

By using the comparison principle, we present a corollary to show the monotonicity of the unique solution $(u_1, u_2; g, h)$ to \eqref{a07} with respect to the parameters, and denote the solution by $(u_1^{\mu}, u_2^{\mu}; g^{\mu}, h^{\mu})$, where $\mu:=(\mu_1, \mu_2)$.

\begin{cor}
\label{cor4.2}
Suppose $\mathbf{(J)}$ holds and $(u_1^{\mu}, u_2^{\mu}; g^{\mu}, h^{\mu})$ is the solution of \eqref{a07} with the initial functions satisfying \eqref{a08}. If $\mu_1^*\leq \mu_1^{**}$ and $\mu_2^*\leq \mu_2^{**}$, we have $h^{\mu^*}\leq h^{\mu^{**}}$, $g^{\mu^{*}}\geq g^{\mu^{**}}$ for $t>0$ and $g^{\mu^{*}}<x<h^{\mu^{*}}$.
\end{cor}

\begin{lem}
\label{lem3.5}
 If $h_{\infty}-g_{\infty}<\infty$, then
\be
\lim\limits_{t\rightarrow\infty}\|u_1\|_{C([g(t),h(t)])}=
\lim\limits_{t\rightarrow\infty}\|u_2\|_{C([g(t),h(t)])}=0
\label{c08}
\ee
and $\lambda_1^P([g_{\infty}, h_{\infty}])\geq 0$.
\end{lem}
\bpf
If $\delta=1$, problem \eqref{a07} becomes
\be\left\{
\begin{array}{ll}
u_{1t}=d_1[\int_{g(t)}^{h(t)}J_1(x-y)u_1(t, y)dy-u_1]-b_1 u_1, &t>0,\ -h_0<x< h_0, \\[2mm]
u_{2t}=-k u_2, &t>0,\ -h_0<x< h_0, \\[2mm]
u_i(t,x)=0,  &t>0,\ x\in\{-h_0, h_0\},\\[2mm]
u_i(0,x)=u_{i, 0}(x),   &-h_0\leq x\leq h_0,\,i=1,2,
\end{array} \right.
\label{d18} \ee
which is a nonlocal problem in a fixed interval $[-h_0, h_0]$. It is obviously seen that $h_{\infty}-g_{\infty}=2h_0<\infty$. Let $(\overline u_1,\overline u_2)=(M_1 e^{-b_1 t}, M_2 e^{-k t})$ with $M_i=||u_{i, 0}||_{L^\infty[-h_0, h_0]}$ as an upper solution of \eqref{d18}, we further have \eqref{c08} holds.

If $0\leq \delta<1$, the principal eigenvalue $\lambda_1^P([g_{\infty}, h_{\infty}])$ of \eqref{b01} with $[-L_1, L_2]$ replaced with $[g_{\infty}, h_{\infty}]$ is well-defined. We firstly prove that $\lambda_1^P([g_{\infty}, h_{\infty}])\geq 0$.  If not, $\lambda_1^P([g_{\infty}, h_{\infty}])< 0$. From \eqref{d02}, we admits that  for a sufficiently small $\epsilon\in(0, h_0)$, there exists a large positive integer $n$, such that  $|g(n\omega)-g_{\infty}|<\epsilon$ and $|h(n\omega)-h_{\infty}|<\epsilon$. Using the continuity (Lemma \ref{lem2.2}) of $\lambda_1^P$ gives that $\lambda_1^P([g(n\omega), h(n\omega)])<0$. In view of $J_1(0)>0$, with the above $\epsilon$ small enough such that $J_1(x)>0$ for $x\in[-4\epsilon, 4\epsilon].$ Let $(\underline{u}_1(t, x), \underline{u}_2(t, x))$ be the
solution of \eqref{c01} with $Q^m_{\textrm w}$ and $Q^m_{\textrm c}$ replaced by $(m\omega, m\omega+(1-\delta)\omega]\times(g(n\omega), h(n\omega))$ and $(m\omega+(1-\delta)\omega, (m+1)\omega]\times(g(n\omega), h(n\omega))$ respectively, and initial functions pair $(\underline{u}_{1}(0,x), \underline{u}_{2}(0,x))=(u_1(n\omega, x), u_2(n\omega, x)).$ By comparison principle, we admits
\bess
(\underline{u}_1, \underline{u}_2)(t, x)\leq (u_1, u_2)(t+n\omega, x) \quad {\rm for}\,\,\,\, t\geq 0,\, g(n\omega)\leq x\leq h(n\omega).
\eess
Then from Lemma \ref{lem3.1} (ii), we deduce that
\bess
(0, 0)<(U_1, U_2)(t, x)=\lim\limits_{t\rightarrow\infty}(\underline{u}_1, \underline{u}_2)(t, x)\leq \liminf\limits_{t\rightarrow\infty}(u_1, u_2)(t, x)
\eess
with the convergence uniform for $x\in[g(n\omega), h(n\omega)]$, where $(U_1(t, x), U_2(t, x))$ is the solution of \eqref{c02} with $[-L_1, L_2]$ replaced by $[g(n\omega), h(n\omega)]$. Therefore, there exists $T_1\geq n\omega$ such that
\bess
0<\frac{1}{2}U_i(t, x)<u_i(t, x) \quad {\rm for}\,\,t\geq T_1,\, x\in[g(n\omega), h(n\omega)], i=1, 2.
\eess
Denote
\bess
c_i:=\min_{\substack{t\geq0,\, g(n\omega)\leq x\leq h(n\omega)}} U_i(t, x)>0 \,(i=1,2)\ {\rm and}\ c_3:=\min\limits_{-4\epsilon\leq x\leq 4\epsilon}J_1(x)>0.
\eess
Noting that the fact $[h(t)-2\epsilon, h(t)-\epsilon]\subset[g(n\omega), h(n\omega)]$ for $t\geq n\omega$, together with \eqref{a07} and the above estimations of $u_1$ and $u_2$, yields
\bess
\begin{array}{lll}
h'(t)&=&\sum\limits_{i=1}^{2}\mu_i\int_{g(t)}^{h(t)}\int_{h(t)}^{+\infty}J_i(x-y)u_i(t,x)dydx\\[2mm]
&\geq&\sum\limits_{i=1}^{2}\mu_i\int_{h(t)-2\epsilon}^{h(t)}\int_{h(t)}^{h(t)+2\epsilon}
J_i(x-y)u_i(t,x)dydx\\[2mm]
&\geq& 2 c_3 \epsilon\sum\limits_{i=1}^{2}\mu_i\int_{h(t)-2\epsilon}^{h(t)}u_i(t,x)dx\\[2mm]
&\geq&2 c_3\epsilon\sum\limits_{i=1}^{2}\mu_i\int_{h(t)-2\epsilon}^{h(t)-\epsilon}u_i(t,x)dx\\[2mm]
&\geq&2 c_3\epsilon\sum\limits_{i=1}^{2}\mu_i\int_{h(t)-2\epsilon}^{h(t)-\epsilon}\frac{1}{2}U_i(t, x)dx\\[2mm]
&\geq&c_3\epsilon^2\sum\limits_{i=1}^{2}\mu_ic_i>0
\end{array}
\eess
for $t\geq T_1$, which contradicts with the fact $h_{\infty}<\infty$. Therefore, we have $\lambda_1^P([g_{\infty}, h_{\infty}])\geq 0$.

Let $(\overline{u}_1(t, x), \overline{u}_2(t, x))$ be the
solution of \eqref{c01} with $Q^m_{\textrm w}$ and $Q^m_{\textrm c}$ replaced by $(m\omega, m\omega+(1-\delta)\omega]\times(g_{\infty}, h_{\infty})$ and $(m\omega+(1-\delta)\omega, (m+1)\omega]\times(g_{\infty}, h_{\infty})$ respectively, and initial functions $(\overline{u}_{1,0}, \overline{u}_{2,0})=(e_1, e_2).$ By comparison principle, we have
\bess
(0, 0)\leq (u_1, u_2)(t, x)\leq (\overline{u}_1, \overline{u}_2)(t, x) \quad {\rm for}\,\,t\geq 0,\, g(t)\leq x\leq h(t).
\eess
Since $\lambda_1^P(g_{\infty}, h_{\infty})\geq 0$, by Lemma \ref{lem3.1}, we have $\lim\limits_{t\rightarrow\infty}(\overline{u}_1, \overline{u}_2)=(0, 0)$ uniformly for $x\in[g_{\infty}, h_{\infty}]$. Thus \eqref{c08} is obtained.
\epf

\vspace{3mm}

\begin{rmk}
From the above proof, we know that,
if $\delta=1$, then $$\lim\limits_{t\rightarrow\infty}\|u_1\|_{C[g(t),h(t)]}=\lim\limits_{t\rightarrow\infty}
\|u_2\|_{C[g(t),h(t)]}=0
$$ always holds, and vanishing happens.
\label{vani}
\end{rmk}

Next, we only need to consider the case $0\leq \delta<1$. Let $\lambda_1^P([-h_0, h_0])$ be the principal eigenvalue of \eqref{b01} with $L_1=L_2 = h_0$.
\begin{lem}
\label{lem3.6}
Assume $0\leq \delta<1$ and $\lambda_1^P([-h_0, h_0])>0$. If $\|u_{1,0}\|_{C[-h_0,h_0]}+\|u_{2,0}\|_{C[-h_0,h_0]}$ is sufficiently small, then
\bess
\lim\limits_{t\rightarrow\infty}\|u_1\|_{C([g(t),h(t)])}=\lim\limits_{t\rightarrow\infty}
\|u_2\|_{C([g(t),h(t)])}=0.
\eess
\end{lem}
\bpf
From Lemma \ref{lem3.5}, we only need to prove $h_{\infty}-g_{\infty}<\infty$.
In fact,
$\lambda_1^P([-h_0, h_0])>0$ together with Lemma \ref{lem2.2} implies that
$\lambda_1^P([-h_1, h_1])>0$ for some $h_1=h_0+\epsilon_0$ with small enough
$\epsilon_0>0$. Let $(\phi, \psi)$ be a positive eigenfunction pair corresponding
to $\lambda_1^P([-h_1, h_1])>0$, and
$\phi(t, x)+\psi(t,x )\leq 1$ for $0\leq t\leq\omega$ and $-h_1\leq x\leq h_1$.

Denote
\bess
\gamma:=\lambda_1^P([-h_1, h_1])/2, \quad M:=\gamma \epsilon_0\bigg(\max\limits_{0\leq t\leq\omega}\int_{-h_1}^{h_1}(\mu_1\phi+\mu_2\psi)dx\bigg)^{-1}.
\eess
Take
\bess
\begin{array}{lll}
\overline{h}(t):=h_0+\epsilon_0[1-e^{-\gamma t}],\,\, & \overline{g}(t):=-\overline{h}(t).\\
\overline{u}_1(t, x):=Me^{-\gamma t}\phi(t, x),& \overline{u}_2(t, x):=Me^{-\gamma t}\psi(t, x)
\end{array}
\eess
for $t\geq 0, x\in[-h_1, h_1].$
From \eqref{b01}, for $t\in(m\omega, m\omega+(1-\delta)\omega]$, $x\in(\overline{g}(t),
\overline{h}(t)),$ we easily obtain
\bess
\begin{array}{lll}
&&\overline{u}_{1t}- d_1\int_{\overline{g}(t)}^{\overline{h}(t)}J_1(x-y)\overline{u}_1(t,y)dy
+d_1\overline{u}_1-a_1(e_1-\overline{u}_1)\overline{u}_2 +b_1 \overline{u}_1\\[2mm]
&\geq&\overline{u}_{1t}- d_1\int_{\overline{g}(t)}^{\overline{h}(t)}J_1(x-y)\overline{u}_1(t,y)dy
+d_1\overline{u}_1-a_1e_1\overline{u}_2 +b_1 \overline{u}_1\\[2mm]
&\textcolor{red}{\geq}&Me^{-\gamma t}\phi_t-\gamma\overline{u}_1-Me^{-\gamma t}\phi_t+ a_1e_1\overline{u}_2-b_1\overline{u}_1+\lambda_1^P([-h_1, h_1])\overline{u}_1-
a_1e_1\overline{u}_2 +b_1 \overline{u}_1\\[2mm]
&=&-\gamma\overline{u}_1+\lambda_1^P([-h_1, h_1])\overline{u}_1=\frac{\lambda_1^P([-h_1, h_1])}{2}\overline{u}_1\geq 0.
\end{array}
\eess
Similarly,
\bess
\overline{u}_{2t}-d_2\int_{\overline{g}(t)}^{\overline{h}(t)}J_2(x-y)\overline{u}_2(t,y)dy
+d_2\overline{u}_2(t, x)-a_2(e_2-\overline{u}_2)\overline{u}_1+b_2 \overline{u}_2\geq 0.
\eess
Next, for $t\in(m\omega+(1-\delta)\omega, (m+1)\omega]$, $x\in(\overline{g}(t), \overline{h}(t)),$
we also obtain
\bess
\begin{array}{lll}
&&\overline{u}_{1t}- d_1\int_{\overline{g}(t)}^{\overline{h}(t)}J_1(x-y)\overline{u}_1(t,y)dy
+d_1\overline{u}_1 +b_1 \overline{u}_1\\[2mm]
&\geq &-\gamma\overline{u}_1+\lambda_1^P([-h_1, h_1])\overline{u}_1=\frac{\lambda_1^P([-h_1, h_1])}{2}\overline{u}_1\geq 0,
\end{array}
\eess
clearly,
\bess
\begin{array}{lll}
\overline{u}_{2t}+k\overline{u}_2&=&-\gamma \overline{u}_2+Me^{-\gamma t}\psi_t+k\overline{u}_2\\
&=&-\gamma\overline{u}_2+\lambda_1^P([-h_1, h_1])\overline{u}_2=\frac{\lambda_1^P([-h_1, h_1])}{2}\overline{u}_2\geq 0.
\end{array}\eess
Since $\overline{h}(t)\in [h_0, h_1)$ for $t \geq 0$, then $[\overline{g}(t), \overline{h}(t)]
\subset(-h_1, h_1)$, we have
\bess
\begin{array}{lll}
&&\sum\limits_{i=1}^{2}\mu_i\int_{\overline{g}(t)}^{ \overline{h}(t)}\int_{\overline{h}(t)}^{+\infty}J_i(x-y)\overline{u}_i(t,x)dydx\\[2mm]
&\leq& \sum\limits_{i=1}^{2}\mu_i\int_{\overline{g}(t)}^{ \overline{h}(t)}\overline{u}_i(t,x)dx\\[2mm]
&=&Me^{-\gamma t}\int_{\overline{g}(t)}^{\overline{h}(t)}\mu_1\phi(t,x)dx+Me^{-\gamma t}\int_{\overline{g}(t)}^{\overline{h}(t)}\mu_2\psi(t,x)dx\\[2mm]
&\leq&Me^{-\gamma t}\int_{-h_1}^{h_1}[\mu_1\phi+\mu_2\psi]dx\\[2mm]
&=&\epsilon_0\gamma e^{-\gamma t}=\overline{h}'(t)
\end{array}
\eess
for $t\in(m\omega, m\omega+(1-\delta)\omega],$ and similarly,
\bess
-\sum\limits_{i=1}^{2}\mu_i\int_{\overline{g}(t)}^{ \overline{h}(t)}\int_{-\infty}^{\overline{g}(t)}J_i(x-y)\overline{u}_i(t,x)dydx\geq\overline{g}'(t).
\eess
Owing to $-\overline{g}'(t)=\overline{h}'(t)\geq 0$ for all $t>0$, then
\bess
\overline{h}(t)\geq \overline{h}(m\omega+(1-\delta)\omega),\  \overline{g}(t)\leq\overline{g}(m\omega+(1-\delta)\omega)
\eess
for $t\in(m\omega+(1-\delta)\omega, (m+1)\omega].$
Clearly we also have
\bess
\overline{u}_i(t, x)\geq 0,&t>0, \,x\in\{\overline{g}(t), \overline{h}(t)\}.
\eess
If we denote
\bess
\sigma:=\min\bigg\{\min\limits_{x\in[-h_0, h_0]}\phi(0, x),\  \min\limits_{x\in[-h_0, h_0]}\psi(0, x)\bigg\},
\eess
then
\be
u_{1,0}(x)\leq M\phi(0, x)\leq\overline{u}_1(0,x), \,\,u_{2,0}(x)\leq M\psi(0, x)\leq\overline{u}_2(0,x)  \quad{\rm for}\,\,x\in[-h_0, h_0],
\label{d08}
\ee
provided that
\bess
\|u_{1,0}\|_{C([-h_0,h_0])}+\|u_{2,0}\|_{C([-h_0,h_0])}\leq \sigma M.
\eess
Therefore, using Lemma \ref{lem3.4}, we conclude that $(\overline{u}_1, \overline{u}_2; \overline{g}, \overline{h})$ is an upper solution of \eqref{a07} and  $[g(t), h(t)]\subset[\overline{g}, \overline{h}]$, hence $h_{\infty}-g_{\infty}\leq \overline{h}_{\infty}-\overline{g}_{\infty}=2h_1$, and the proof is completed.
\epf

\begin{lem}
\label{lem3.7}
If $0\leq \delta<1$ and $\lambda_1^O\geq 0$,  then
\be
\lim\limits_{t\rightarrow\infty}(u_1, u_2)(t, x)=(0, 0)
\label{f04}
\ee
uniformly for $x\in[g(t), h(t)]$ and vanishing occurs.
\end{lem}
\bpf
The result can be proved by using the comparison principle. Let $(\overline u_1(t), \overline u_2(t))$ be the solution to
\be\left\{
\begin{array}{lll}
\overline u_{1t}=a_1(e_1-\overline u_1)\overline u_2-b_1\overline u_1, \qquad \qquad &m\omega<t\leq m\omega+(1-\delta)\omega,\\[2mm]
\overline u_{2t}=a_2(e_2-\overline u_2)\overline u_1-b_2\overline u_2,  &m\omega<t\leq m\omega+(1-\delta)\omega,\\[2mm]
\overline u_{1t}=-b_1\overline u_1, &m\omega+(1-\delta)\omega<t\leq (m+1)\omega,\\[2mm]
\overline u_{2t}=-k \overline u_2,  &m\omega+(1-\delta)\omega<t\leq (m+1)\omega,\\[2mm]
\overline u_i(0)=e_{i}\  (i=1,2).
\end{array}
\label{f11}
\right.
\ee
It is easily seen that
\bess
d_i\mathcal{L}_i[\overline u_i](t)=d_i\int_{g(t)}^{h(t)}J_i(x-y)\overline{u}_i dy-d_i\overline{u}_i\leq 0
\eess
and $u_{i}(0, x)\leq e_i=\overline{u}_i(0)$, it follows from Lemma \ref{lem3.10} that
\be
(u_1, u_2)(t, x)\leq (\overline u_1, \overline u_2)(t)
\label{f03}
\ee
for $x\in[g(t),
h(t)]$ and $t\geq 0$.

Set \[K_1=a_1e_2+b_1, ~~~K_2=a_2e_1+b_2+k. \]
Using $(\overline{u}_1^{(0)}, \overline{u}_2^{(0)})=(e_1, e_2)$, as the initial values, iteration sequence $\{(\overline{u}_1^{(n)}, \overline{u}_2^{(n)})\}$ is obtained through the following iterative process
{\small
\bess
\left\{
\begin{array}{lll}
\overline{u}_{1t}^{(n)}+K_1\overline{u}_1^{(n)}
=K_1\overline{u}_1^{(n-1)}+a_1(e_1-\overline{u}^{(n-1)}_1)
\overline{u}^{(n-1)}_2 -b_1 \overline{u}^{(n-1)}_1, &t\in (0, (1-\delta)\omega],\\[2mm]
\overline{u}_{2t}^{(n)}+K_2\overline{u}_2^{(n)}
=K_2\overline{u}_2^{(n-1)}+a_2(e_2-\overline{u}^{(n-1)}_2)
\overline{u}^{(n-1)}_1 -b_2 \overline{u}^{(n-1)}_2, &t\in (0, (1-\delta)\omega],\\[2mm]
\overline{u}_{1t}^{(n)}+K_1\overline{u}_1^{(n)}
=K_1\overline{u}_1^{(n-1)} -b_1 \overline{u}^{(n-1)}_1, &t\in ((1-\delta)\omega, \omega],\\[2mm]
\overline{u}_{2t}^{(n)}+K_2\overline{u}_2^{(n)}
=K_2\overline{u}_2^{(n-1)} -k \overline{u}^{(n-1)}_2, &t\in ((1-\delta)\omega, \omega],\\[2mm]
\overline{u}^{(n)}_i(0)=\overline{u}^{(n-1)}_i(\omega),
\end{array}
\right.
\eess}
where $n=1,2,\cdots,$ $\{(\overline{u}^{(n)}_1, \overline{u}^{(n)}_2)\}$ is called the maximal sequence. From the comparison principle, it easily follows that the above sequence $\{(\overline{u}^{(n)}_1, \overline{u}^{(n)}_2)\}$ admits the monotone property
\bess
\left.
\begin{array}{lll}
(\overline{u}_1^{(n)}, \overline{u}_2^{(n)})\leq(\overline{u}_1^{(n-1)}, \overline{u}_2^{(n-1)})\leq\cdots\leq(\overline{u}_1^{(1)}, \overline{u}_2^{(1)})\leq(\overline{u}_1^{(0)}, \overline{u}_2^{(0)})=(e_1, e_2)
\end{array}
\right.
\eess
for $t\in[0, \omega]$, and therefore the limits exist, denoted by
\[\lim_{n\rightarrow\infty}(\overline{u}_1^{(n)}, \overline{u}_2^{(n)})(t+n\omega)=({U}_1, {U}_2)(t), \]
which yields
\begin{equation*}
\left.
\begin{array}{lll}
({U}_1, {U}_2)
\leq(\overline{u}_1^{(n)}, \overline{u}_2^{(n)})\leq(\overline{u}_1^{(n-1)}, \overline{u}_2^{(n-1)})\leq(e_1, e_2).
\end{array}
\right.
\end{equation*}
Furthermore,  Sobolev imbedding theorem assert that $({U}_1, {U}_2)$  meets with problem \eqref{c11}. It is
clearly seen that $({U}_1$, ${U}_2)$ is the solutions of problem \eqref{c11}.
Since that problem \eqref{c11} has a unique solution $(U_1^{\vartriangle}, U_2^{\vartriangle})$, then $({U}_1, {U}_2)=(U_1^{\vartriangle}, U_2^{\vartriangle})$, that is,
\bess
\lim_{n\rightarrow\infty}(\overline{u}_1^{(n)}, \overline{u}_2^{(n)})=(U_1^{\vartriangle}, U_2^{\vartriangle}).
\eess
Recalling that $(\overline{u}_1, \overline{u}_2)(t)
\leq(e_1, e_2)=(\overline{u}_1^{(0)}, \overline{u}_2^{(0)})(t) $ for $t\in[0, \omega]$, then
$(\overline{u}_1, \overline{u}_2)(\omega)
\leq(e_1, e_2)=(\overline{u}_1^{(0)}, \overline{u}_2^{(0)})(\omega)= (\overline{u}_1^{(1)}, \overline{u}_2^{(1)})(0)$.
Using comparison principle that $(\overline{u}_1, \overline{u}_2)(t+\omega)
\leq(\overline{u}_1^{(1)}, \overline{u}_2^{(1)})(t)$.
We eventually conclude that for $t\in [0, \omega]$ and $n=0, 1, 2\ldots$,
\bess
(\overline{u}_1, \overline{u}_2)(t+n\omega)
\leq(\overline{u}_1^{(n)}, \overline{u}_2^{(n)})(t)
\eess
holds by iterating $\overline{u}^{(n)}_i(0)=\overline{u}^{(n-1)}_i(\omega)$ with respect to $n$.
Therefore
\be
\limsup_{n\rightarrow\infty}(\overline{u}_1, \overline{u}_2)(t+n\omega)=(U_1^{\vartriangle}, U_2^{\vartriangle})(t).
\label{f16}
\ee
Similarly to the proof of Lemmas \ref{lem2.1} and \ref{lem3.1}, we show that problem \eqref{c11} has the only nonnegative solution $(0,0)$ if $\lambda_1^O\geq0$, that is, $(U_1^{\vartriangle}, U_2^{\vartriangle})=(0, 0)$.
As a result,
\bess
\lim_{t\rightarrow\infty}(u_1, u_2)(t, x)=(0, 0)
\eess
uniformly for $x\in[g(t), h(t)]$ because of \eqref{f03} and \eqref{f16}, which means that vanishing occurs.
\epf

\begin{lem}
\label{vanish}
If $a_1a_2e_1e_2\leq b_1b_2$, then $h_\infty-g_\infty<\infty$ and vanishing happens.
\end{lem}
\bpf
 In fact, using Lemma \ref{lem3.5}, it suffices to prove that $h_\infty-g_\infty<\infty$. The result can also be obtained by using the energy-type estimate, see also \cite{Du2020}

Recalling that
$\int_{\mathbb{R}}J_1(x)dx = 1$ and $J_1$ is symmetric, we have
{\small
\be
\begin{array}{lll}
&&\int_{g(t)}^{h(t)}\mathcal{L}_1[u_1](t,x)dx\\[2mm]
&=&
\int_{g(t)}^{h(t)}\big[\int_{g(t)}^{h(t)}J_1(x-y)u_1(t,y)dy-u_1(t,x)\big]dx\\[2mm]
&=&\int_{g(t)}^{h(t)}\int_{g(t)}^{h(t)}J_1(x-y)[u_1(t,y)-u_1(t,x)]dy dx\\[2mm]
&&-\int_{g(t)}^{h(t)}\int_{h(t)}^{\infty}J_1(x-y)u_1(t,x)dy dx
-\int_{g(t)}^{h(t)}\int_{-\infty}^{g(t)}J_1(x-y)u_1(t,x)dydx\\[2mm]
&=&-\int_{g(t)}^{h(t)}\int_{h(t)}^{\infty}J_1(x-y)u_1(t,x)dydx
-\int_{g(t)}^{h(t)}\int_{-\infty}^{g(t)}J_1(x-y)u_1(t,x)dy dx.
\end{array}
\label{d19}\ee}
Similarly,
{\small
\be
\int_{g(t)}^{h(t)}\mathcal{L}_2[u_2](t,x)dx=-\int_{g(t)}^{h(t)}\int_{h(t)}^{\infty}J_2(x-y)u_2(t,x)dydx
-\int_{g(t)}^{h(t)}\int_{-\infty}^{g(t)}J_2(x-y)u_2(t,x)dydx.
\label{d20}\ee
}
Using the aforementioned equations and \eqref{a07} for $t\in (m\omega, m\omega+(1-\delta)\omega]$, we obtain
{\small
\be
\begin{array}{lll}
&&\frac{d}{dt}\int_{g(t)}^{h(t)}[u_1(t,x)+\frac{a_1e_1}{b_2}u_2(t,x)]dx\\[2mm]
&=&\int_{g(t)}^{h(t)}[u_{1t}(t,x)+\frac{a_1e_1}{b_2}u_{2t}(t,x)]dx\\[2mm]
&=&\int_{g(t)}^{h(t)}[d_1\mathcal{L}_1[u_1]+a_1(e_1-u_1)u_2-b_1u_1
+\frac{a_1e_1}{b_2}(d_2\mathcal{L}_2[u_2]+a_2(e_2-u_2)u_1-b_2u_2)]dx\\[2mm]
&=&\int_{g(t)}^{h(t)}(d_1\mathcal{L}_1[u_1]+\frac{a_1e_1d_2}{b_2}\mathcal{L}_2[u_2])dx
+\int_{g(t)}^{h(t)}[(\frac{a_1a_2e_1e_2}{b_1b_2}-1)b_1 u_1-(a_1+\frac{a_1a_2e_1}{b_2})u_1u_2]dx\\[2mm]
&\leq&\int_{g(t)}^{h(t)}(d_1\mathcal{L}_1[u_1]+\frac{a_1e_1d_2}{b_2}\mathcal{L}_2[u_2])dx\\[2mm]
&\leq&-\min\{d_1, \frac{a_1e_1d_2}{b_2}\}[\int_{g(t)}^{h(t)}\int_{h(t)}^{\infty}J_1(x-y)u_1(t,x)dydx
+\int_{g(t)}^{h(t)}\int_{-\infty}^{g(t)}J_1(x-y)u_1(t,x)dydx\\[2mm]
&&+\int_{g(t)}^{h(t)}\int_{h(t)}^{\infty}J_2(x-y)u_2(t,x)dydx
+\int_{g(t)}^{h(t)}\int_{-\infty}^{g(t)}J_2(x-y)u_2(t,x)dydx]\\[2mm]
&\leq&-\frac{\min\{d_1, \frac{a_1e_1d_2}{b_2}\}}{\max\{\mu_1, \mu_2\}}[h'(t)-g'(t)]:=-D[h'(t)-g'(t)]
\end{array}
\label{d09}
\ee}
by using \eqref{d19} and \eqref{d20}. Further, integrating  \eqref{d09} from $m\omega$ to $t$ yields
\bess
\begin{array}{lll}
&&[h(t)-g(t)]-[h(m\omega)-g(m\omega)]\\[2mm]
&\leq& -\int_{g(t)}^{h(t)}\frac{1}{D}[u_1+\frac{a_1e_1}{b_2}u_2](t,x)dx
+\int_{g(m\omega)}^{h(m\omega)}\frac{1}{D}[u_1
+\frac{a_1e_1}{b_2}u_2](m\omega,x)dx
\end{array}
\eess
for the warm season $t\in(m\omega, m\omega+(1-\delta)\omega]$.

Denote
$$S(t)=h(t)-g(t)\ \textrm{and}\  F(t, x)=\frac{1}{D}[u_1+\frac{a_1e_1}{b_2}u_2](t,x),$$
then the above inequation is rewritten as
\be
S(t)-S(m\omega)\leq -\int_{g(t)}^{h(t)}F(t, x)dx+\int_{g(m\omega)}^{h(m\omega)}F(m\omega, x)dx.
\label{d10}
\ee
For cold season $t\in ((m-1)\omega+(1-\delta)\omega, m\omega],$
\be
S(m\omega)-S((m-1)\omega+(1-\delta)\omega)=0.
\label{d13}\ee
Notice that
{\small
\be\begin{array}{lll}
&&\int_{g(m\omega)}^{h(m\omega)}F(m\omega, x)dx-\int_{g((m-1)\omega+(1-\delta)\omega)}^{h((m-1)\omega+(1-\delta)\omega)}
F\big((m-1)\omega+(1-\delta)\omega, x\big)dx\\[2mm]
&=&\frac{1}{D}\int_{g(m\omega)}^{h(m\omega)}[u_1
+\frac{a_1e_1}{b_2}u_2](m\omega,x)-[u_1
+\frac{a_1e_1}{b_2}u_2]((m-1)\omega+(1-\delta)\omega, x)dx\\[2mm]
&=&\frac{1}{D}\int_{(m-1)\omega+(1-\delta)\omega}^{m\omega}\int_{g(m\omega)}^{h(m\omega)}[u_{1t}
+\frac{a_1e_1}{b_2}u_{2t}](t, x)dxdt\\[2mm]
&=&\frac{1}{D}\int_{(m-1)\omega+(1-\delta)\omega}^{m\omega}
\int_{g(m\omega)}^{h(m\omega)}(d_1\mathcal{L}_1[u_1]-b_1u_1)(t,x)dxdt\\[2mm]
&&-\frac{ka_1e_1}{Db_2}\int_{(m-1)\omega+(1-\delta)\omega}^{m\omega}
\int_{g(m\omega)}^{h(m\omega)}u_2(t,x)dxdt\\[2mm]
&\leq&\frac{d_1}{D}\int_{(m-1)\omega+(1-\delta)\omega}^{m\omega}
\int_{g(m\omega)}^{h(m\omega)}\mathcal{L}_1[u_1](t,x)dxdt\\[2mm]
&=&\frac{d_1}{D}\int_{(m-1)\omega+(1-\delta)\omega}^{m\omega}\big\{\int_{g(m\omega)}^{h(m\omega)}
\int_{g(t)}^{h(t)}J_1(x-y)[u_1(t,y)-u_1(t,x)]dy dx\\[2mm]
&&-\int_{g(t)}^{h(t)}\int_{h(t)}^{\infty}J_1(x-y)u_1(t,x)dy dx
-\int_{g(t)}^{h(t)}\int_{-\infty}^{g(t)}J_1(x-y)u_1(t,x)dydx\big\}dt\\[2mm]
&\leq&\frac{d_1}{D}\int_{(m-1)\omega+(1-\delta)\omega}^{m\omega}\int_{g(m\omega)}^{h(m\omega)}
\int_{g(t)}^{h(t)}J_1(x-y)[u_1(t,y)-u_1(t,x)]dy dxdt\\[2mm]
&=&0
\end{array}
\label{d14}
\ee}
by using $h(t)=h(m\omega+(1-\delta)\omega)$ and $g(t)=g(m\omega+(1-\delta)\omega)$ for
$m\omega+(1-\delta)\omega<t\leq (m+1)\omega$. It follows from \eqref{d13} and \eqref{d14}
that{\small
\be\begin{array}{lll}
&&S(m\omega)-S((m-1)\omega+(1-\delta)\omega)\\[2mm]
&\leq& -\int_{g(m\omega)}^{h(m\omega)}F(m\omega, x)dx+\int_{g((m-1)\omega+(1-\delta)\omega)}^{h((m-1)\omega+(1-\delta)\omega)}
F\big((m-1)\omega+(1-\delta)\omega, x\big)dx.
\end{array}
\label{d15}
\ee}
Adding \eqref{d10} and \eqref{d15} yields
{\small
\be
\begin{array}{lll}
&&S(t)-S((m-1)\omega+(1-\delta)\omega)\\[2mm]
&\leq& -\int_{g(t)}^{h(t)}F(t, x)dx+\int_{g((m-1)\omega+(1-\delta)\omega)}^{h((m-1)\omega+(1-\delta)\omega)}
F\big((m-1)\omega+(1-\delta)\omega, x\big)dx
\end{array}
\ee}
for $t\in (m\omega, m\omega + (1- \delta)\omega].$

Repeat the above process, yields
\bess
\begin{array}{lll}
S(t)-S(0)
&\leq& -\int_{g(t)}^{h(t)}F(t, x)dx+\int_{g(0)}^{h(0)}F(0, x)dx\\[2mm]
&=&\int_{g(0)}^{h(0)}\frac{1}{D}[u_1(0,x)+\frac{a_1e_1}{b_2}u_2(0,x)]dx,
\end{array}
\eess
that is
\bess
h(t)-g(t)\leq 2h_0+\int_{-h_0}^{h_0}\frac{1}{D}[u_1(0,x)+\frac{a_1e_1}{b_2}u_2(0,x)]dx.
\eess
As a result, for all $t>0$ we have $h(t)-g(t)$ is bounded, which implies $h_{\infty}-g_{\infty}<\infty$ and vanishing happens.
\epf

\vspace{5mm}

\begin{lem}
\label{lem3.9}
If\, $0\leq \delta<1$ and $\lambda_1^P([g(t_0),h(t_0)])\leq0$ for some $t_0\geq 0$, then $g_{\infty}=-\infty, h_{\infty}=+\infty$ and
\be
\lim\limits_{n\rightarrow\infty}(u_1, u_2)(t+n\omega, x)=(U_1^{\vartriangle}, U_2^{\vartriangle})(t)
\label{d03}
\ee
 uniformly for $t\in [0, \omega]$ and locally uniformly for $x\in \mathbb{R}$, where $\lambda_1^P([g(t_0),h(t_0)])$ is the principal eigenvalue of \eqref{c01} with $[-L_1, L_2]$ replaced by $[g(t_0),$  $h(t_0)]$, and $(U_1^{\vartriangle}, U_2^{\vartriangle})(t)$ is the unique positive periodic solution to problem \eqref{c11}.
\end{lem}
\bpf
Notice that $[g(t_0), h(t_0)]\subsetneq[g(t_1), h(t_1)]$ for $t_1> t_0+\omega$, the monotonicity of $\lambda_1^P$ gives that
\be
\lambda_1^P([g(t_1), h(t_1)])<0 \quad {\rm for}\,\,t_1> t_0+\omega.
\label{d04}
\ee
Using the above inequality, we get by the similar proof of Lemma \ref{lem3.5} that $h_{\infty}=+\infty$ and $g_{\infty}=-\infty$.

Now, we prove \eqref{d03}. For big integer $n$ ($n > [\frac{t_0}{\omega}]+1$),  and let $(u_{1,n \omega}(t, x), u_{2, n \omega}(t, x))$ be the solution of
\eqref{c01} with $Q^m_{\textrm w}$ and $Q^m_{\textrm c}$ replaced by $(m\omega, m\omega+(1-\delta)\omega]\times(g(n \omega), h(n \omega))$ and $(m\omega+(1-\delta)\omega,$  $(m+1)\omega]\times(g(n \omega), h(n \omega))$ respectively, and initial functions $(u_{1,n \omega}(0,x), u_{2, n \omega}(0,x))=(u_1(n \omega, x), u_2(n \omega, x)).$ By comparison principle, we obtain
\bess
(u_{1,n \omega}, u_{2,n \omega})(t+m\omega, x)\leq (u_1, u_2)(t+m\omega+n \omega, x)
\eess
for $(t, x)\in [0,\omega]\times[g(n \omega), h(n \omega)], m=0, 1, \dots$. Letting $m\to \infty$ and using Lemma \ref{lem3.1} (ii)  yield
{\small
\bess
(0, 0)<(U_{1,n \omega}, U_{2,n \omega})(t, x)=\lim\limits_{m\rightarrow\infty}(u_{1,n \omega}, u_{2,n \omega})(t+m\omega, x)\leq \liminf\limits_{m\rightarrow\infty}(u_1, u_2)(t+m\omega+n \omega, x)
\eess}
with the convergence uniform for $t\in [0, \omega]$ and $x\in[g(n \omega), h(n \omega)]$, where $(U_{1,n \omega}, U_{2,n \omega})(t, x)$ is the positive periodic solution of \eqref{c02} with $[-L_1, L_2]$ replaced by $[g(n \omega), h(n \omega)]$.
Noticing the fact that $(g(n \omega), h(n \omega))\rightarrow(-\infty, \infty)$ as $n\rightarrow\infty$,  we obtain by Lemma \ref{lem3.2} that
\be
({U}_1^{\vartriangle}, {U}_2^{\vartriangle})(t)\leq \liminf\limits_{n\rightarrow\infty}(u_1, u_2)(t+n\omega, x)
\label{d05}
\ee
uniformly for $t\in [0,\omega]$ and locally uniformly for  $x\in \mathbb{R}$.
Moreover, it follows from Lemma \ref{lem3.7} that
\be
\limsup\limits_{n\rightarrow\infty}(u_1, u_2)(t+n\omega, x)\leq ({U}_1^{\vartriangle}, {U}_2^{\vartriangle})(t),
\label{d06}
\ee
where $(U_1^{\vartriangle}, U_2^{\vartriangle})$ is the unique solution of problem \eqref{c11}.
This completes the proof.
\epf

\vspace{5mm}
{\bf Proof of Theorem \ref{thm1.1}.} We first consider the solution in the warm season of the first year ($m=0$). Similarly as in \cite{Cao2019, Du2020},
denote {\small\bess
\begin{array}{lll}
&\mathbb{H}_{h_0,(1-\delta)\omega}:=\{h\in C^1([0, (1-\delta)\omega]): h(0)=h_0,\ h(t)\,\textrm{is strictly increasing}\}, \\[2mm]
&\mathbb{G}_{h_0,(1-\delta)\omega}:=\{g\in  C^1([0, (1-\delta)\omega]): -g\in \mathbb{H}_{h_0,(1-\delta)\omega}\}.
\end{array}
\eess}
Similarly as in [\cite{Du2020}, Lemma 4.2], for any $(g, h)\in \mathbb{G}_{h_0, (1-\delta)\omega}\times \mathbb{H}_{h_0, (1-\delta)\omega}$, the following problem
{\small
\be\left\{
\begin{array}{ll}
u_{1t}=d_1\mathcal{L}_1[u_1]+a_1(e_1-u_1)u_2 -b_1 u_1, &0<t\leq (1-\delta)\omega, g(t)<x< h(t), \\[2mm]
u_{2t}=d_2\mathcal{L}_2[u_2]+a_2(e_2-u_2) u_1 -b_2 u_2,  &0<t\leq (1-\delta)\omega, g(t)<x< h(t),\\[2mm]
u_1(t,x)=u_2(t,x)=0,  &0<t\leq (1-\delta)\omega, x\in\{g(t), h(t)\},\\[2mm]
h(0)=-g(0)=h_0,\\[2mm]
u_1(0,x)=u_{1,0}(x),\ u_2(0,x)=u_{2,0}(x),   &-h_0\leq x\leq h_0.
\end{array} \right.
\label{d07} \ee
}
has a unique positive solution for $t\in(0,s]$ by the contraction mapping theorem, we can further extend the unique solution defined over $(0,s]$ to $(0,(1-\delta)\omega]$ with respect to $t$. As a result, problem \eqref{d07} with the initial function pair $(u_{1,0}, u_{2,0})$ admits a unique solution $(u_1^{g,h}, u_2^{g,h})$ which satisfies
$$0<u_i^{g,h}\leq e_i\,\,  {\rm for }\, \, (t,x)\in (0, (1-\delta)\omega]\times(g(t), h(t)).$$

We now return to free boundary problem \eqref{a07}.
Then, we define a mapping $\mathcal{F}(h,g)=(\hat{h}, \hat{g})$ by
{\small
\bess\left\{
\begin{array}{ll}
\hat{h}(t)=h_0+\sum\limits_{i=1}^{2}\mu_i\int_{0}^{t}\int_{g(s)}^{h(s)}\int_{h(s)}^{+\infty}
J_i(x-y)u_i(s,x)dydxds, &0<t\leq (1-\delta)\omega,\\[2mm]
\hat{g}(t)=-h_0-\sum\limits_{i=1}^{2}\mu_i\int_{0}^{t}\int_{g(s)}^{h(s)}\int_{-\infty}^{g(s)}
J_i(x-y)u_i(s,x)dydxds, &0<t\leq (1-\delta)\omega.
\end{array} \right.
\eess
}
To prove that the existence and uniqueness of solution to \eqref{a07} for $t\in(0, (1-\delta)\omega]$, it suffices to show that $\mathcal{F}$ has a unique fixed point in  $\Sigma_{(1-\delta)\omega}
:=\mathbb{H}_{h_0,(1-\delta)\omega}\times \mathbb{G}_{h_0,(1-\delta)\omega}$.
To get the above conclusion, we firstly show that for small $0<\tau <(1-\delta)\omega$, $\mathcal{F}$ maps a suitable closed subset of $\Sigma_{\tau}$ into itself and is a contraction
mapping, and then showing that any fixed point of $\mathcal{F}$ in $\Sigma_{\tau}$ belongs to this closed subset. This
guarantees the existence of a unique solution of \eqref{a07} in the range $t \in(0, \tau]$. Then one shows
that the unique solution can be extended to $t \in(0, (1-\delta)\omega]$.
Since the proof is almost the same as that in [\cite{Cao2019}, Theorem 2.1], we omit the details here.

On the other hand, over the time interval $((1-\delta)\omega, \omega]$ (the cold season), we find that the boundaries are fixed, that is, $g(t)\equiv g((1-\delta)\omega)$ and $h(t)\equiv h((1-\delta)\omega)$ for $t\in ((1-\delta)\omega, \omega]$. Then fixed boundary problem
{\small
\be\left\{
\begin{array}{ll}
u_{1t}=d_1\mathcal{L}_1[u_1] -b_1 u_1, &(1-\delta)\omega<t\leq\omega,\, g((1-\delta)\omega)<x< h((1-\delta)\omega), \\[2mm]
u_1(t,x)=0,  &(1-\delta)\omega<t\leq\omega,\, x\in\{g((1-\delta)\omega), h((1-\delta)\omega)\}
\end{array} \right.
\label{d0711} \ee
}
admits a unique solution $u_1$ by the result of [\cite{Cao2019}, Lemma 2.3], moreover $u_1$ satisfies
$$0<u_1\leq  e_1\,\,  {\rm for }\, \, (t,x)\in ((1-\delta)\omega, \omega]\times(g((1-\delta)\omega), h((1-\delta)\omega)).$$
Again, solving the forth equation of \eqref{a07} gives that
 $$u_2(t, x)=u_{2}((1-\delta)\omega, x)\, e^{k[(1-\delta)\omega-t]}\leq e_2,$$
  for $(t, x)\in [(1-\delta)\omega, \omega]\times [g((1-\delta)\omega), h((1-\delta)\omega)]$.

By taking $m = 1, 2\cdots$, recursively, we therefore obtain the existence and uniqueness of the solution $(u_1, u_2, g, h)$ to \eqref{a07} for $t\in [0, \infty)$.  It also easily follows from  $(u_1, u_2)\in [C([0, \infty)\times(g(t), h(t))]^2$ and $(h, g)\in [C([0, \infty))\cup C^1([m\omega, m\omega+(1-\delta)\omega])]^2$
that $0<u_i \leq e_i(i=1,2)$ in $\{(t,x):\, t>0,\, g(t)\leq x\leq h(t)\}$, $-g^{'}, h^{'}>0$ in the warm season $(m\omega, m\omega+(1-\delta)\omega]$ and $g^{'}=h^{'}=0$ in the cold season $(m\omega+(1-\delta)\omega, (m+1)\omega]$, hence $g(t)$ is nonincreasing and $h(t)$ is nondecreasing in $(0, \infty)$. This completes the proof of  Theorem \ref{thm1.1}.
\epf

\vspace{5mm}
{\bf Proof of Theorem \ref{thm1.2}.}\,\,
$(1)$ If $\delta=1$, vanishing always happens by Remark \ref{vani}.

$(2)$ If $0\leq\delta<1$ and $\lambda_1^O<0$, by Lemma \ref{lem2.3} (ii)(iii),
we find $\lambda_1^P([g(t_0), h(t_0)]) < 0$ for some large $t_0 > 0$. According to Lemma \ref{lem3.9}, \eqref{d03} holds.

$(3)$ If $0\leq\delta<1$ and $\lambda_1^O\geq0$,  then \eqref{f04} holds by Lemma \ref{lem3.7}.

Therefore, spreading (case (2))-vanishing (cases (1) and (3)) dichotomy holds.
\epf

\begin{lem}
\label{thm4.111}
Suppose $\mathbf{(J)}$, $J_1(x)=J_2(x)$ and $0\leq\delta<1$ hold, and the initial functions satisfy \eqref{a08}. Let $\lambda_1^P(\Gamma, b_1-d_1\lambda_1^*, b_2-d_2\lambda_1^*, \delta)$ is given by \eqref{b01}, denote
$$\lambda_1^O(\Gamma, b_1, b_2, \delta):=\lim\limits_{L_1, L_2\rightarrow+\infty}
\lambda_1^P(\Gamma, b_1-d_1\lambda_1([-L_1, L_2]), b_2-d_2\lambda_1([-L_1, L_2]), \delta),$$
where $\Gamma=\{a_i, e_i, \omega, k\}$. The following conclusion is valid.

$(i)$\ If $\lambda_1^O\geq 0$, then vanishing always happens.

$(ii)$\ If $\lambda_1^O<0$ and  $\lambda_1^P([-h_0, h_0])\leq 0,$
then spreading always happens.

$(iii)$\ If $\lambda_1^O<0$ and $\lambda_1^P([-h_0, h_0])> 0,$
 then

 \quad$(a)$ \ For any given initial datum $(u_{1, 0}, u_{2, 0})$ satisfying \eqref{a08}, there exists $\mu^*\geq\mu_*>0$ such that vanishing happens for $0<\mu_1+\mu_2\leq \mu_*$ and spreading happens for $\mu_1+\mu_2>\mu^*$.

\quad $(b)$\ For fixed $\mu_1, \mu_2>0$, vanishing happens if initial datum $(u_{1, 0}, u_{2, 0})$ is small enough.
\end{lem}

\bpf
$(i)$ If $0\leq\delta<1$ and $\lambda_1^O\geq 0,$ then vanishing happens by Lemma \ref{lem3.7}.

$(ii)$ If $\lambda_1^O<0$ and  $\lambda_1^P([-h_0, h_0])\leq 0,$ we derive that $\lambda_1^P([g(t_0), h(t_0)])<0$ for some large $t_0 > 0$ by continuity in Lemma \ref{lem2.3}, hence
spreading occurs by Lemma \ref{lem3.9}.

$(iii)$ According to $\lambda_1^P([-h_0, h_0])>0$, then conclusion $(b)$ is directly proved by Lemma \ref{lem3.6}.  It remains to prove $(a)$.

We construct the same upper solution as in the proof Lemma \ref{lem3.6}, since that
$$ \lim\limits_{(\mu_1, \mu_2)\rightarrow (0,0)}M:=\lim\limits_{(\mu_1, \mu_2)\rightarrow(0, 0)}\gamma \epsilon_0\bigg(\max\limits_{0\leq t\leq\omega}\int_{-h_1}^{h_1}(\mu_1\phi+\mu_2\psi)dx\bigg)^{-1}=\infty,$$
therefore, for any given initial function pair $(u_{1,0}, u_{2,0})$ satisfying \eqref{a08}, there exists $\mu_*>0$ such that \eqref{d08} holds for $0<\mu_1+\mu_2\leq \mu_*$, which indicates that vanishing occurs for problem \eqref{a07}.

Next, if $0\leq \delta<1$, we claim that there exists $\mu^*>0$ such that spreading occurs when $\mu_1+\mu_2>\mu^*$.
To emphasize the dependence of solution on $\mu:=(\mu_1, \mu_2)$, let $(u_1^{\mu}, u_2^{\mu}; g^{\mu}, h^{\mu})$ denote the unique positive solution of \eqref{a07}.
First of all, we show that there exists $\mu_1^*>0$ such that for some large $t_1>0$,
\be
\lambda_1^P([g^{\mu_1^*}(t_1), h^{\mu_1^*}(t_1)])<0.
\label{a10}
\ee
If not, then we assume that for all $t>0$ and $\mu>0$,
\bess
\lambda_1^P([g^{\mu}(t), h^{\mu}(t)])\geq 0.
\eess
Note that $-g^{\mu}(t)$ and $h^{\mu}(t)$ are nondecreasing with respect to $t$, we obtain that $-g^{\mu}(t)$ and $h^{\mu}(t)$ are also nondecreasing in $\mu>0$ by Corollary \ref{cor4.2}. Denote
\bess
H_{\infty}:=\lim\limits_{t, \mu\rightarrow+\infty}h^{\mu}(t), \,\,G_{\infty}:=\lim\limits_{t, \mu\rightarrow+\infty}g^{\mu}(t).
\eess
And recall that $J_i(0)>0$, there exist $\epsilon_0, \delta_0$ such that $J_i(x)>\delta_0$ when $|x|<\epsilon_0$. Then for the above $\epsilon_0$, there exist large enough positive constants $\mu_0$ and $ m_0$ such that for $t_0\in(m_0\omega, m_0\omega+(1-\delta)\omega]$
$$h^{\mu}(t)>H_{\infty}-\frac{\epsilon_0}{4},\qquad  {\rm when\ }\mu\geq\mu_0, t\geq t_0.$$
Thus, for the sixth equation of \eqref{a07}, integrating both sides of the equation from $t_0$ to $m_0\omega+(1-\delta)\omega$ with respect to $t$, gives
{\small
\bess
\begin{array}{lll}
&&h^{\mu}(m_0\omega+(1-\delta)\omega)-h^{\mu}(t_0)\\[2mm]
&=&
\sum\limits_{i=1}^{2}\mu_i\int_{t_0}^{m_0\omega+(1-\delta)\omega}\int_{g^{\mu}(\tau)}^{h^{\mu}(\tau)}
\int_{h^{\mu}(\tau)}^{+\infty}J_i(x-y)u_i^{\mu}(\tau,x)dydxd\tau,\\[2mm]
&\geq& \sum\limits_{i=1}^{2}\mu_i\int_{t_0}^{m_0\omega+(1-\delta)\omega}\int_{g^{\mu_0}(\tau)}^{h^{\mu_0}(\tau)}
\int_{h^{\mu_0}(\tau)+\frac{\epsilon_0}{4}}^{+\infty}J_i(x-y)u_i^{\mu_0}(\tau,x)dydxd\tau,\\[2mm]
&\geq& \sum\limits_{i=1}^{2}\mu_i\int_{t_0}^{m_0\omega+(1-\delta)\omega}
\int_{h^{\mu_0}(\tau)-\frac{\epsilon_0}{2}}^{h^{\mu_0}(\tau)}
\int_{h^{\mu_0}(\tau)+\frac{\epsilon_0}{4}}^{h^{\mu_0}(\tau)
+\frac{\epsilon_0}{2}}J_i(x-y)u_i^{\mu_0}(\tau,x)dydxd\tau,\\[2mm]
&\geq& \frac{\epsilon_0}{4}\delta_0\sum\limits_{i=1}^{2}\mu_i\int_{t_0}^{m_0\omega+(1-\delta)\omega}
\int_{h^{\mu_0}(\tau)-\frac{\epsilon_0}{2}}^{h^{\mu_0}(\tau)}u_i^{\mu_0}(\tau,x)dxd\tau,\\[2mm]
&\geq& \frac{\epsilon_0}{4}\delta_0\sum\limits_{i=1}^{2}\mu_i\int_{t_0}^{m_0\omega+(1-\delta)\omega}
\int_{h^{\mu_0}(\tau)-\frac{\epsilon_0}{2}}^{h^{\mu_0}(\tau)}u_{\min}^{\mu_0}(\tau,x)dxd\tau,
\end{array}
\eess
}
where $u_{\min}^{\mu_0}(t,x)=\min\limits_{\substack{m_0\omega<t\leq m_0\omega+(1-\delta)\omega\\g(t)<x< h(t)}}\{u_1^{\mu_0}(t,x),\, u_2^{\mu_0}(t,x)\}$.
Therefore, we have
{\footnotesize
\bess\mu_1+\mu_2\leq \bigg( \frac{\epsilon_0}{4}\delta_0\int_{t_0}^{m_0\omega+(1-\delta)\omega}
\int_{h^{\mu_0}(\tau)-\frac{\epsilon_0}{2}}^{h^{\mu_0}(\tau)}u_{\min}^{\mu_0}(\tau,x)dxd\tau \bigg)^{-1}[h^{\mu}(m_0\omega+(1-\delta)\omega)-h^{\mu}(t_0)]<+\infty,
\eess}
which is a contradiction. Therefore, \eqref{a10} holds. We can conclude that spreading occurs when $\mu=\mu_1^*$ by Lemma \ref{lem3.9}.
\epf

\vspace{5mm}
{\bf Proof of Theorem \ref{thm1.3}.}\,\,
$(i)$\ Fix $a_i, e_i, b_2$ and $k$, then $b_1^*(:=a_1a_2e_1e_2/b_2)$ is determined. It follows from  Fig. \ref{tu3}
that $\lambda_1^O\geq 0$ if $b_1\in[b_1^*, +\infty)$.
While if $b_1\in(0, b_1^*)$, there exists $\delta^*\in(0, 1)$ such that $\lambda_1^O\geq 0$ for $\ \delta\in[\delta^*, 1]$.
Then vanishing always happens from Lemma \eqref{thm4.111} $(i)$.

$(ii)$\ If $b_1\in(0, b_1^*)$ and $\delta\in(0, \delta^*)$, then $\lambda_1^O<0$. Recalling that $\lambda_1^P([-L_1, L_2])$ is decreasing with respect to $L_1+L_2$, there exists $h_0^*>0$ such that $\lambda_1^P([-h_0, h_0])\leq 0$ for $h_0\geq h_0^*$. Therefore, for some $t_0\geq 0$ such that $g(t_0)\leq-h_0^*$ and $h(t_0)\geq h_0^*$, we have $\lambda_1^P([g(t_0), h(t_0)])\leq 0$. Hence spreading happens from Lemma \eqref{thm4.111} $(ii)$.

$(iii)$\ If $b_1\in(0, b_1^*)$, $\delta\in(0, \delta^*)$ and $h_0\in(0, h_0^*)$, then $\lambda_1^O<0$ and $\lambda_1^P([-h_0, h_0])>0$. The proof is finished by Lemma \eqref{thm4.111} $(iii)$ $(a)$.
\epf

\begin{lem}
\label{thm4.54}
Suppose $\mathbf{(J)}$ holds, and the initial functions satisfy \eqref{a08}.  There exists $0\leq \mu_\vartriangle\leq\mu^\vartriangle\leq +\infty$ such that vanishing happens for $(0, 0)<(\mu_1, \mu_2)< (\mu_\vartriangle, \mu_\vartriangle)$, and spreading happens for $(\mu_1, \mu_2)>(\mu^\vartriangle, \mu^\vartriangle)$.
Specially, if $\underline\lambda_1^O\geq 0$, then $\mu_\vartriangle=\mu^\vartriangle=+\infty$, which means that vanishing happens for any $\mu_1$ and $\mu_2$. While, if $\overline\lambda_1^O< 0$ and $\overline \lambda_1^P([-h_0, h_0])\leq 0$, then $\mu_\vartriangle=\mu^\vartriangle=0$, which means that spreading happens for any $\mu_1$ and $\mu_2$.
\end{lem}

\bpf
 Define
$$\mu_\vartriangle=\sup\{\,\mu^{\star}\in [0, +\infty): \ \textrm{vanishing occurs for}\ (\mu_1, \mu_2)\ \textrm{with}\ (\mu_1, \mu_2)\leq (\mu^{\star}, \mu^{\star})\, \}$$
and
$$\mu^\vartriangle=\inf\{\,\mu^{\star}\in (0, +\infty]: \ \textrm{spreading occurs for}\ (\mu_1, \mu_2)\ \textrm{with}\ (\mu_1, \mu_2)> (\mu^{\star}, \mu^{\star})\, \}.$$
It is easy to see that $\mu_\vartriangle$, $\mu^\vartriangle$ are well-defined, and $0\leq \mu_\vartriangle, \mu^\vartriangle\leq +\infty$. Owing to $h(t)-g(t)$ is nondecreasing in $t$, recalling that $h(\infty)-g(\infty)<\infty$ for the vanishing case and $h(\infty)-g(\infty)=\infty$ for the spreading case yields that $\mu_\vartriangle\leq \mu^\vartriangle$. Moreover,  vanishing happens for $(0, 0)<(\mu_1, \mu_2)< (\mu_\vartriangle, \mu_\vartriangle)$, and spreading happens for $(\mu_1, \mu_2)>(\mu^\vartriangle, \mu^\vartriangle)$.

If $\underline\lambda_1^O\geq 0$, similarly we can prove that $h_\infty-g_\infty<\infty$ by using the energy-type integral as in Lemma \ref{lem3.7} and
$$\lim\limits_{t\rightarrow\infty}\|u_1\|_{C([g(t),h(t)])}=
\lim\limits_{t\rightarrow\infty}\|u_2\|_{C([g(t),h(t)])}=0$$
by constructing an upper solution as in Lemma \ref{lem3.5}, so vanishing happens for any $(\mu_1, \mu_2)$, which implies that $\mu_\vartriangle=\mu^\vartriangle=+\infty$.

  If $\overline\lambda_1^O< 0$, we have that $\delta \neq 1$, and the assumptions $\overline\lambda_1^O< 0$ and $\overline \lambda_1^P([-h_0, h_0])\leq 0$ are equivalent to that $\lambda_1^O< 0$ and $\lambda_1^P([-h_0, h_0])\leq 0$, which immediately gives that spreading happens for any $(\mu_1, \mu_2)$ by Lemma \ref{lem3.9}, in this case, $\mu_\vartriangle=\mu^\vartriangle=0$.
\epf

\vspace{5mm}
{\bf Proof of Theorem \ref{thm1.4}}\,\,
If $\delta=1$, then vanishing always happens and $\mu_\vartriangle=\mu^\vartriangle=+\infty$.

If $0\leq \delta<1$ and  $\min\{b_1, k\}\delta>c_1(1-\delta)$, then $\lambda_1^O\geq 0$ from Corollary \ref{cor1.3}. According to Lemma \ref{thm4.54}, $\mu_\vartriangle=\mu^\vartriangle=+\infty$, which means that vanishing happens for any $\mu_1$ and $\mu_2$.

While if $0\leq\delta<1$ and $\max\{b_1, k\}\delta<c_1(1-\delta)$, then $\lambda_1^O<0$ from Corollary \ref{cor1.3}. Recalling that $\lambda_1^P([-L_1, L_2])$ is decreasing with respect to $L_1+L_2$, there exists a sufficiently large constant $h_0^{**}>0$ such that $\lambda_1^P([-h_0, h_0])\leq 0$ for $h_0\geq h_0^{**}$. It concluded that $\mu_\vartriangle=\mu^\vartriangle=0$, which means that spreading happens for any $\mu_1$ and $\mu_2$.
\epf

\section{Discussion}

As is known to all, WNv usually spreads among people by means of mosquitoes infected with the virus after biting the infected birds. On the one hand, mosquitoes are very climate sensitive and breed from June to November (warm season) each year with unknown changing habitat, whose fronts are described by free boundaries. However, for harsh living conditions in cold season, the population of mosquitoes is assumed to follow the Malthusian growth law i.e. decaying at an exponential rate, does not migrate and stays in a hibernating status. Mathematically, the density $u_2$ of infected mosquitoes satisfies a nonlocal diffusion with free boundary in the warm season and degenerate differential equation without diffusion in a fixed interval in the cold season.

On the other hand, the ecological pressure, being brought by the increased food-demand for active proliferation during the warm season, causes certain birds to fly to the northern glacier retreating place in summer and then return to the south again to spend the winter as soon as the glaciers recover. In other words, the density $u_1$ of infected birds always satisfies a reaction-diffusion equation both in the warm season and in the cold season.

In this paper, we study the WNv nonlocal model with free boundaries and seasonal succession, which can effectively describe the long-distance diffusion process, in fact, it is more suitable than local diffusion for models characterizing relatively dense population densities. We define the generalized principal eigenvalues of nonlocal diffusion operators with seasonal succession, and put forward a calculation method for the principal eigenvalue of the corresponding ODE system (Theorem \ref{thm2.2}). What's more, the monotonicity of (generalized) principal eigenvalues with respect to length $(L)$ of the interval (Lemmas \ref{lem2.7} and \ref{lem2.2}),  the duration ratio $(\delta)$ of the cold season  (Lemma \ref{lem2.5} and Corollary \ref{rmk1.2}) are given.

We prove that the spreading-vanishing dichotomy (Theorem \ref{thm1.2}) is valid,
and the criteria (Theorems \ref{thm1.3} and \ref{thm1.4}) that completely determine when spreading and vanishing can happen are deduced. Our theoretical results show that the longer the warm season lasts, the bigger the risk of infection gets and the less beneficial for the prevention and control of WNv (Remark \ref{rmk111}).

\end{document}